\documentclass{article}

\newcommand{\Z}{\mathbb{Z}}

\newcommand{\R}{\mathbb{R}}
\newcommand{\C}{\mathbb{C}}
\newcommand{\D}{\mathbb{D}}
\newcommand{\T}{\mathbb{T}}

\newcommand{\at}{{\mathcal A}_ T}
\newcommand{\iit}{{\mathcal I}_ T}
\newcommand{\LL}{L^{1}_{\omega_ T}}
\newcommand{\M}{{\mathcal M}_{\omega_ T}}
\newcommand{\ma}{{\mathcal A}}
\newcommand{\mb}{{\mathcal B}}

\newcommand{\qm}{{\mathcal Q}{\mathcal M}}
\newcommand{\llt}{-\lambda \Delta_T}
\newcommand{\sab}{S_{\alpha,\beta}}

\newcommand{\FB}{{\mathcal{FB}}}
\newcommand{\scab}{\overline {S}_{\alpha,\beta}}
\newcommand{\mga}{{\mathcal G}({\mathcal A})}

\newcommand{\scabp}{{\overline S}_{-{\pi\over 2}-\alpha,{\pi\over 2}-\beta}}

\newcommand{\uab}{{\mathcal U}_{\alpha,\beta}}

\newcommand{\dfi}{Dom(\FB(\phi))}

\newcommand{\fab}{{\mathcal F}_{\alpha, \beta}}

\newcommand{\gab}{{\mathcal G}_{a, b}}

\newcommand{\dab}{\Delta_{\alpha, \beta}}

\newcommand{\QM}{{\mathcal {QM}}}

\newcommand{\jjt}{{\mathcal I}_ {T,\mathcal B}}

\newcommand{\dt}{{\mathcal U}_ T}

\newcommand{\sia}{\widehat{\iit}}
\newcommand{\bbt}{{\mathcal B}_ T}
\newcommand{\mub}{{\bf \mu}}
\newcommand{\vab}{{\mathcal V}_{\alpha,\beta}}
\newcommand{\ssab}{S^*_{\alpha, \beta}}
\newcommand{\sscab}{\overline{S}^*_{\alpha, \beta}}
\newcommand{\usab}{{\mathcal U}^*_{\alpha,\beta}}
\newcommand{\vsab}{{\mathcal V}^*_{\alpha, \beta}}
\newcommand{\att}{A_ T}

\newcommand{\wiit}{\widehat{{\mathcal I}_T}}
\newcommand{\fb}{{\mathcal F}{\mathcal B}}
\newcommand{\0}{\{0\}}
\usepackage[ansinew]{inputenc}
\usepackage[T1]{fontenc}
\usepackage[english]{babel}
\usepackage{graphicx}
\usepackage{color}
\usepackage{amsfonts}
 \newtheorem{thm}{Theorem}[section]
 \newtheorem{cor}[thm]{Corollary}
 \newtheorem{lem}[thm]{Lemma}
 \newtheorem{prop}[thm]{Proposition}
  \newtheorem{defn}[thm]{Definition}
  
\begin{document}
%
%
%
%
%
%
%
%
%

\title
 {A holomorphic functional calculus for finite families of commuting semigroups}
\author{Jean Esterle}

\maketitle


\begin{abstract}

Let $\ma$ be a commutative Banach algebra such that $u\ma \neq \0$ for $u\in \ma \setminus \0$ which possesses dense principal ideals. The purpose of the paper is to give a general framework to define $F(-\lambda_1\Delta_{T_1}, \dots,-\lambda_k\Delta_{T_k})$ where $F$ belongs to a natural class of holomorphic functions defined on suitable open subsets of $\C^k$ containing the "Arveson spectrum" of $(-\lambda_1\Delta_{T_1}, \dots, -\lambda_k\Delta_{T_k})$, where $\Delta_{T_1}, \dots, \Delta_{T_k}$ are the infinitesimal generators of commuting one-parameter semigroups of multipliers on $\ma$ belonging to one of the following classes

(1) The class of strongly continous semigroups $T=(T(te^{ia})_{t>0}$ such that $\cup_{t>0}T(te^{ia})\ma$ is dense in $\ma,$ where $a \in \R.$

(2) The class of semigroups $T=(T(\zeta))_{\zeta \in S_{a, b}}$ holomorphic on an open sector $S_{a, b}$ such that $T(\zeta)\ma$ is dense in $\ma$ for some, or equivalently for all $\zeta \in S_{a, b}.$

We use the notion of quasimultiplier, introduced in 1981 by the author at the Long Beach Conference on Banach algebras: the generators of the semigroups under consideration will be defined as quasimultipliers on $\ma,$ and for $\zeta$ in the Arveson resolvent set $\sigma_{ar}(\Delta_T)$ the resolvent $(\Delta_T -\zeta I)^{-1}$ will be defined as a regular quasimultiplier on $\ma,$ i.e. a quasimultiplier $S$ on $\ma$ such that $\sup_{n\ge 1}\lambda^n\Vert S^n u\Vert <+\infty$ for some $\lambda >0$ and some $u$ generating a dense ideal of $\ma$ and belonging to the intersection of the domains of $S^n$, $n\ge 1.$ 

The first step consists in "normalizing" the Banach algebra $\ma,$ i.e. continuously embedding $\ma$ in a Banach algebra $\mathcal B$ having the same quasimultiplier algebra as $\ma$ but for which $\lim \sup_{t\to 0^+}\Vert T(te^{ia})\Vert_{{\mathcal M}(\mathcal B)}<+\infty$ if $T$ belongs to the class (1), and for which $\lim \sup _{\stackrel{\zeta \to 0}{_{\zeta \in S_{\alpha, \beta}}}}\Vert T(\zeta)\Vert <+\infty$ for all pairs $(\alpha, \beta)$ such that $a < \alpha < \beta < b$ if $T$ belongs to the class (2). Iterating this procedure  this allows to consider $(\lambda_j \Delta_{T_j}+\zeta I)^{-1}$ 
as an element of $\mathcal M(\mathcal B)$ for $\zeta \in Res_{ar}(-\lambda_j\Delta_{T_j}),$ the "Arveson resolvent set " of $-\lambda_j\Delta_{T_j},$ and to use the standard integral 'resolvent formula' even if the given semigroups are not bounded near the origin.

 A first approach to the functional calculus involves the dual $\mathcal G_{a,b}$ of an algebra of fast decreasing functions, described in Appendix 2.  Let $a=(a_1,\dots,a_k),$ $b=(b_1,\dots, b_k),$ with $a_j \le b_j \le a_j+\pi,$ and denote by $M_{a,b}$ the set of families  $(\alpha,\beta)=(\alpha_1, \beta_1), \dots, (\alpha_k, \beta_k)$ such that $\alpha_j=\beta_j=a_j$ if $a_j=b_j$ and such that $a_j <\alpha_j\le \beta_j <b_j$ if $a_j<b_j.$ Let $\uab$ denote the class of all functions $f: \zeta= (\zeta_1,\dots,\zeta_k)\to f(\zeta)$  continuous on the product space $\overline S_{\alpha, \beta}=\Pi_{j\le k}\overline S_{\alpha_j,\beta_j}$ and converging to $0$ at infinity such that the function $\sigma \to f(\zeta_1, \dots , \zeta_{j-1}, \eta, \zeta_{j+1},\dots, \zeta_k)$ is holomorphic on $S_{\alpha_j,\beta_j}$ whenever $\alpha_j<\beta_j.$ Elements of the dual $\uab'$ admit a "representing measure", and we describe in appendix 1
 some certainly well-known ways to implement the duality between $\uab$ and $\uab'$ and extend the action of elements of $\uab'$ to vector-valued analogs spaces $\uab(X)$ and their "bounded" counterparts $\vab(X)$ via representing measures and Cauchy and Fourier-Borel transforms. In appendix 2 we introduce a natural algebra of fast decreasing functions, which is the intersection for $(\alpha, \beta)\in M_{a,b}$  and $z \in \C^k$ of all $e_{-z}\uab,$ where $e_z(\zeta)=e^{z_1\zeta_1+\dots+z_k\zeta_k}$. The dual $\mathcal G_{a,b}$ of this algebra is an algebra with respect to convolution, this dual space is the union for $z\in \C^k$ and $(\alpha, \beta)\in M_{a,b}$ of the dual spaces $(e_{-z}\uab)',$ and elements
 of these dual spaces act on the vector-valued spaces $e_{-z}\uab(X).$ This action can also be implemented via representing measures, Cauchy transforms and Fourier Borel transforms as indicated in appendix 2. If $\lambda_jS_{\alpha_j,\beta_j}$ is contained in the domain of definition of $T_j$ for $j\le k,$ this allows to define the action of $\phi \in \mathcal G_{a,b}$ on $T_{(\lambda)}=(T_1(\lambda_1.), \dots, T_k(\lambda_k.)$ by using the formula
 
 $$<T_{(\lambda)},\phi>u=<T_1(\lambda_1\zeta_1)\dots T_{k}(\lambda_k\zeta_k)u, \phi_{\zeta_1,\dots,\zeta_k}> \ \ (u\in \mb),$$
 
 where $ \mathcal B$ denotes a normalization of the given commutative Banach algebra $\ma$  with respect to $T=(T_1,\dots,T_k),$ when $\phi \in (e_{-z}\uab)'$ and when 
 $\sup_{\zeta \in \overline S_{\alpha, \beta}}\left \Vert e^{z\zeta}T_1(\lambda_1\zeta_1)\dots T_k(\lambda_k\zeta_k)\right \Vert_{ \mathcal M(\mb)}<+\infty.$

 For $(\alpha, \beta)\in M_{a,b},$ set $S^*_{\alpha,\beta}=\Pi_{j\le k}S_{-{\pi\over 2}-\alpha_j,{\pi\over 2}-\beta_j}.$ An open set $U=\Pi_{\j \le k}U_j$ is said to be admissible with respect to $(\alpha, \beta)$ if for every $j \le k$ the boundary $\partial U_j$ is a piecewise ${\mathcal C}^1$-curve, if $U + \epsilon \subset U$ for every $\epsilon \in \overline S^*_{\alpha, \beta}$ and if $\overline S^*_{\alpha, \beta}\setminus U$ is compact. Standard properties of the class $H^{(1)}(U)$  of all holomorphic functions $F$ on $U$ such that
 $\Vert F\Vert_{H^{(1)}(U)}:=\sup_{\epsilon \in S^*_{\alpha, \beta}}\int_{\epsilon +\tilde \partial U}\Vert F(\sigma \vert d\sigma\vert<+\infty$ are given in appendix 3 (when $a_j=b_j$ for $j\le k,$ this space is the usual Hardy space $H^1$ on a product of open half-planes).
 
 The results of appendix 3 allow when an open set $U\subset \C^k$ admissible with respect to $(\alpha, \beta)\in M_{a,b}$ satisfies some more suitable admissibility conditions with respect to $T=(T_1,\dots,T_k)$ and $\lambda \in \cup_{(\gamma, \delta) \in M_{a-\alpha, b-\beta}}\overline S_{\gamma, \delta}$ to define $F(-\lambda_1\Delta_{T_1},\dots,-\lambda_k\Delta_{T_k})$ for $F\in H^{(1)}(U)$ by using the formula
 
 $$F(-\lambda_1\Delta_{T_1},\dots,-\lambda_k\Delta_{T_k})$$ $$={1\over (2i\pi)^k}\int_{\epsilon +\overline \partial U}F(\zeta_1,\dots,\zeta_k)(\lambda_1\Delta_{T_1} +\zeta_1 I)^{-1}\dots (\lambda_1\Delta_{T_k} +\zeta_k I)^{-1}d\zeta_1\dots d\zeta_k,$$

 where $\tilde \partial U$ denotes the "distinguised boundary of $U$ and where $\epsilon \in S^*_{\alpha, \beta}$ is choosen so that $\epsilon +U$ still satisfies the required admissibility conditions with respect to $T$ and $\lambda.$ Given $T$ and $\lambda$, this gives a family ${\mathcal W}_{T,\lambda}$ of open sets stable under finite intersection
 and an algebra homomorphism $F \to F(-\lambda_1\Delta_{T_1},\dots,-\lambda_k\Delta_{T_k})$ from $\cup_{U\in \mathcal W_{T,\lambda}}H^{(1)}(U)$ into the multiplier algebra
 $\mathcal M(\mb)\subset \qm_r(\ma).$ This homomorphism extends in a natural way to a bounded algebra homomorphism from $\cup_{U\in \mathcal W_{T,\lambda}}H^{\infty}(U)$
 into $\qm_r(\mb)=\qm_r(\ma),$ and we have, if  $\phi \in \mathcal F_{\alpha, \beta}$ for some $(\alpha, \beta)\in M_{a,b},$ and if $\lim_{\stackrel {\vert t \vert \to +\infty}{_{t\in \overline S_{\alpha,\beta}}}}\Vert e^{t_1z_1+\dots +t_tz_k}T_1(t_1\lambda_1)\dots T_k(t_k\lambda_k)\Vert=0$ for some $z$ in the domain of the Fourier-Borel transform of $\phi.$

 $$\mathcal{FB}(\phi)(-\lambda_1\Delta_{T_1},\dots,-\lambda_k\Delta_{T_k})=<T_{(\lambda)},\phi>,$$
 
 so that  $F(-\lambda_1\Delta_{T_1},\dots, -\lambda_k\Delta_{T_k})=T(\nu \lambda_j)$ if $F(\zeta)=e^{-\nu \zeta_j},$ when $\nu \lambda_j$ is in the domain of definition of $T_j.$

 A function $F \in H^{\infty}(U)$ will be said to be strongly outer if there exists a sequence $(F_n)_{n\ge 1}$ of invertible elements of $H^{\infty}(U)$ such that $\vert F(\zeta)\vert \le \vert F_n(\zeta)\vert$ and $\lim_{n\to +\infty}F(\zeta)F_n^{-1}(\zeta)=1$ for $\zeta \in U.$ Every bounded outer function on the open unit disc $\D$ is strongly outer, but the class of strongly outer functions on $\D^k$ is smaller than the usual class of bounded outer functions on $\D^k$ if $k\ge 2.$ We then define the Smirnov class $\mathcal S(U)$ to be the class of those holomorphic functions $F$ on $U$ such that $FG\in H^{\infty}(U)$ for some strongly outer function $G \in H^{\infty}(U).$ The bounded algebra homomorphism $F\to F(-\lambda_1\Delta_{T_1},\dots, -\lambda_k\Delta_{T_k})$ from $\cup_{U\in \mathcal W_{T,\lambda}}$ into $\qm_r(\mb)=\qm_r(\ma)$ extends to a bounded homomorphism from $\cup_{U\in \mathcal W_{T,\lambda}}\mathcal S(U)$ into $\qm(\mb)=\qm(\ma).$ If $F:\zeta \to \zeta_j$is the $j$-th coordinate projection, then of course $F(-\lambda_1\Delta_{T_1},\dots, -\lambda_k\Delta_{T_k})=\lambda_j\Delta_{T_j}.$

\end{abstract}



\section{Introduction}

The author observed in \cite{e2} that if  a  Banach algebra $A$ does not possess any nonzero idempotent then $\inf_{\stackrel{x\in A}{_{\Vert x \Vert\ge 1/2}}}\Vert x^2-x\Vert \ge 1/4.$ If $x$ is quasinilpotent, and if $\Vert x \Vert \ge 1/2,$ then $\Vert x\Vert >1/4.$ Concerning (nonzero) strongly continuous semigroups  ${\bf T} =(T(t))_{t>0}$ of bounded operators on a Banach space $X,$
these elementary considerations lead to the following results, obtained in 1987 by Mokhtari \cite{mo}

\begin{enumerate}

\item If $\lim \sup_{t\to 0^+}\Vert T(t)-T(2t)\Vert <1/4,$ then the generator of the semigroup is bounded, and so $\lim \sup \Vert T(t)-T(2t)\Vert =0.$

\item If the semigroup is quasinilpotent, then $\Vert T(t)-T(2t)\Vert >1/4$ when $t$ is sufficiently small.

\end{enumerate}

 If the semigroup is norm continuous, and if there exists a sequence $(t_n)_{n\ge 1}$ of positive real numbers such that $\lim_{n\to +\infty}t_n=0$ and $\Vert T(t_n) -T(2t_n)\Vert <1/4,$ then the closed subalgebra ${\mathcal A}_{\bf T}$ of  ${\mathcal B}(X)$ generated by the semigroup possesses an exhaustive sequence of idempotents, i.e. there exists a sequence $(P_n)_{n\ge 1}$ of idempotents of ${\mathcal A}_{\bf T}$ such that for every compact set $K \subset \widehat {\mathcal A}_{\bf T}$ there exists $n_K>0$ satisfying $\chi(P_n)=1$ for $\chi \in K,$ $n\ge n_K$.

More sophisticated arguments allowed A. Mokhtari and the author to obtain in 2002 in \cite{em} more general results valid for every integer $p\ge 1.$

These results led the author to consider in \cite{e3}  the behavior of the distance $\Vert T(s) -T(t)\Vert$ for $s > t$ near $0.$ The following results were obtained in \cite{e3}

\begin{enumerate}

\item If there exist for some $\delta >0$ two continuous functions $r\to t(r)$ and $r\to s(r)$ on $[0,\delta],$ such that $s(0)=0$ and such that $0<t(r)<s(r)$  and  $\Vert T(t(r))-T(s(r))\Vert < (s(r)-t(r))\frac{s(r)^{s(r)\over s(r)-t(r)}}{t(r)^{t(r)\over s(r)-t(r)}}$ for $r \in (0,\delta],$ then the generator of the semigroup is bounded, and so $\Vert T(t) -T(s)\Vert \to 0$ as $0<t<s, s\to 0^+.$

\item If the semigroup is quasinilpotent, there exists $\delta >0$ such that $\Vert T(t)-T(s)\Vert > (s-t)\frac{s^{s\over s-t}}{t^{t\over s-t}}$ for $0<t<s\le \delta.$

\item If the semigroup is norm continuous, and if there exists two sequences of positive real numbers such that $0<t_n<s_n,$ $\lim_{n\to +\infty}s_n=0,$ and such that $\Vert T(t_n)-T(s_n)\Vert < (s_n-t_n)\frac{s_n^{s_n\over s_n-t_n}}{t_n^{t_n\over s_n-t_n}},$ then the closed subalgebra ${\mathcal A}_{\bf T}$ of  ${\mathcal B}(X)$ generated by the semigroup possesses an exhaustive sequence of idempotents.
\end{enumerate}

The quantities appearing in these statements are not mysterious: consider the Hilbert space $L^2([0,1],$ and for $t>0$ define $T_0(t): L^2([0,1] \to L^2([0,1]$ by the formula $T_0(t)(f)(x)=x^tf(x) \ \ (0<x\le 1).$ Then $\Vert T_0(t) -T_0(s)\Vert = (s-t)\frac{s^{s\over s-t}}{t^{t\over s-t}}.$ This remark also shows that  assertions (1) and (3) in these statements are sharp, and examples show that  assertion $(2)$ is also sharp.

One can consider $T(t)$ as defined by the formula $\int_0^{+\infty}T(x)d\delta_t(x),$ where $\delta_t$ denotes the Dirac measure at $t.$ Heuristically, $T(t)=e^{t\Delta_T},$ where $\Delta_T$ denotes the generator of the semigroup, and since the Laplace transform of $\delta_t$ is defined by the formula ${\mathcal L}(\delta_t)(z)=\int_0^{+\infty}e^{-zx}d\delta_t(x)=e^{-zt},$ it is natural to write ${\mathcal L}(\delta_t)(-\Delta_T)=T(t).$ More generally, if an entire function $F$ has the form $F={\mathcal L}(\mu),$ where $\mu$ is a measure supported by $[a,b],$ with $0<a<b<+\infty,$ we can set

$$F(-\Delta_T)=\int_0^{+\infty}T(x)d\mu(x),$$

and  consider the behavior of the semigroup near 0 in this context.

I. Chalendar, J.R. Partington and the author used this point of view in \cite{cep1}. Denote by ${\mathcal M}_c(0,+\infty)$ the set of all measures $\mu$  supported by some interval $[a,b],$ where $0<a<b<+\infty.$ For the sake of simplicity we restrict attention to statements analogous to assertion 2. The following result is proved in \cite{cep1}

 Theorem: {\it  Let $\mu \in {\mathcal M}_c(0,+\infty)$ be a nontrivial  real measure such that $\int_{0}^{+\infty}d\mu(t)=0$ and let ${\bf T}=(T(t))_{t>0}$ be a quasinilpotent semigroup of bounded operators.Then there exists $\delta>0$ such that $\Vert F(-s\Delta_T)\Vert> \max_{x\ge 0}\vert F(x)\vert$ for $0 <s\le \delta.$}

When $\mu=\delta_1-\delta_2$ this gives assertion 3 of Mokhtari's result, and when $\mu =\delta_1-\delta_{p+1}$ this gives assertion 3 of the extension of Mokhtari's result given in \cite{em} (but several variables extensions of this functional calculus would be needed  in order to obtain extensions of the results of \cite{e3}). 

This theory applies, for example, to quantities of the form $\Vert T(t)-2T(2t) +T(3t)\Vert,$ or Bochner integrals $\Vert \int_1^2T(tx)dx-\int_2^3T(tx)dx\Vert,$ which are not accessible by the methods of \cite{mo} or \cite{em}. Preliminary results concerning semigroups holomorphic in a sector were obtained by I. Chalendar, J.R. Partington and the author in \cite{cep2}.

More generally it would be interesting to obtain lower estimates as $(\lambda_1, \dots, \lambda_k)\to (0,\dots,0)$ for quantities of the form $F(-\lambda_1\Delta_T,\dots,-\lambda_k\Delta_T)$ when the generator $\Delta_T$ of the semigroup is unbounded, and when $F$ is an analytic function of several complex variables  defined and satisfying natural growth conditions on a suitable neighbourhood of $\sigma_{ar}(\Delta_T)$, where $\sigma_{ar}(\Delta_T)$ denotes the "Arveson spectrum" of the infinitesimal generator $\Delta_{T}$ of $T.$ The purpose of the present paper is to pave the way to such a program by defining more generally $F(-\lambda_1\Delta_{T_1},\dots, -\lambda_k\Delta_{T_k})$ when $F$ belongs to a suitable class of holomorphic functions on some element of a family $\mathcal W_{T_1,\dots,T_k, \lambda}$ of open sets, and where $(T_1,\dots,T_k)$ denotes a finite family of commuting semigroups.

More precisely consider $a=(a_1,\dots,a_k)\in \R^k, b=(b_1,\dots,b_k)\in \R^k)$ satisfying $a_j\le b_j \le a_j+\pi$ for $j\le k,$ and consider a commutative Banach algebra $\ma$
such that $u\ma$ is dense in $\ma$ for some $u\in \ma$ and such that $u\ma \neq \{0\}$ for $u \in \ma \setminus \{0\}.$ This allows to consider the algebra $\qm(\ma)$ of all quasimultipliers on $\ma$ and the algebra $\qm_r(\ma)$ of all regular quasimultipliers on $\ma$ introduced by the author in \cite{e2}, see section 2, and the usual algebra ${\mathcal M}(\ma)$ of all multipliers on $\ma$ can be identified to the algebra of all quasimultipliers on $\ma$ of domain equal to the whole of $\ma.$ We will be interested here in finite families
$(T_1,\dots, T_k)$ of commuting semigroups of multipliers on $\ma$ satisfying the following conditions

\begin{itemize}

\item the semigroup $T_j$ is strongly continuous on $e^{ia_j}.(0,+\infty),$ and $\cup_{t>0}T_j(e^{ia_j}t)\ma$ is dense in $\ma$ if $a_j=b_j,$

\item the semigroup $T_j$  is holomorphic on the open sector $S_{a_j,b_j}:=\{ z \in \C \setminus \{0\} \ | \ a_j <arg(z)<b_j\}$ and $T_j(\zeta)\ma$ is dense in $\ma$ for some (or, equivalently, for all) $ \zeta \in S_{a_j,b_j}$ if $a_j<b_j.$

\end{itemize}

The first step of the construction consists in obtaining a "normalization" $\ma_T$ of the Banach algebra $\ma$ with respect to a strongly continuous one-parameter semigroup $(T(t))_{t>0}$ of multipliers on $\ma.$ The idea behind this normalization process goes back to Feller \cite{f1}, and we use for this the notion of "QM-homomorphism" between commutative Banach algebras introduced in section 2, which seems more appropriate than the related notion of "$s$-homomorphism" introduced by the author in \cite{e2}. Set $\omega_T=\Vert T(t)\Vert$ for $t>0.$ A slight improvement of a result proved by P. Koosis and the author in section 6 of \cite{e1} shows that the weighted convolution algebra $L^{1}(\R^+, \omega_T)$ possesses dense principal ideals, which allows to construct in section 3 a commutative Banach algebra $\ma_T\subset \qm_r(\ma)$ which contains $\ma$ as a dense subalgebra and has dense principal ideals such that the injection $\tilde j: \qm_(\ma)\to \qm(\ma_T)$ associated to the norm-decreasing inclusion map $j:\ma \to \ma_{T}$ is onto and such that $\tilde j ({\mathcal M}(\ma)) \subset  {\mathcal M}(\ma_T)$ for which $\lim \sup_{t\to 0^+}\Vert T(t) \Vert_{\mathcal M(A_{\bf T})}<+\infty.$ Set $\phi_T (f)=\int_0^{+\infty}f(t)T(t)dt$ for $f \in   L^{1}(\R^+, \omega_T)$, where the Bochner integral is computed with respect to the strong operator topology on $\mathcal M(\ma),$ and denote by $\mathcal I_T$ the closed subalgebra of $\mathcal M(\ma)$ generated by $\phi_T(L^{1}(\R^+, \omega_T)).$ In section 5 we give an interpretation of the generator $\Delta_T$ of the semigroup $T$ as a quasimultiplier on $\iit,$ and we define the "Arveson spectrum" $\sigma_{ar}(\Delta_T)$ to be the set $\{\tilde \chi(\Delta_T)\}_{\chi \in \widehat{\iit}},$ where $\tilde \chi$ denotes the unique extension to $\qm(\iit)$ of a character $\chi$ on $\iit,$ with the convention  $\sigma_{ar}(\Delta_T)=\emptyset$ if the "Arveson ideal" $\iit$ is radical. The quasimultiplier $\Delta_T-\lambda I$ is invertible in $\qm_r(\iit)$ and $(\Delta_T-\lambda I)^{-1}\in \mathcal M(\ma_T) \subset \qm_r(\ma)$ if $\lambda \in \C\setminus \sigma_{ar}(\Delta_T),$ and we observe in section 6 that we have, for $\zeta > \lim \sup_{t\to +\infty}{log \Vert T(t)\Vert\over t},$

$$(\Delta_T-\zeta I)^{-1}=-\int_0^{+\infty} e^{-\zeta t}T(t)dt \in \qm(\ma_T)\subset \qm_r(\ma),$$

which is the usual "resolvent formula" extended to strongly continuous semigroups not necessarily bounded near the origin.

In section 4 we construct a more sophisticated normalization of the Banach algebra $\ma$ with respect to a semigroup $T=(T(\zeta))_{\zeta \in S_{a, b}}$ which is holomorphic on an open sector $S_{a, b},$ where $a < b \le a +\pi.$ In this case the normalization $\ma_T$ of $\ma$ with respect to the semigroup $T$ satisfies two more conditions

\begin{itemize}

\item $T(\zeta)u\ma_T$ is dense in $\ma_T$ for $\zeta \in S_{a, b}$ if $u \ma$ is dense in $\ma,$ 

\item $\lim \sup_{\stackrel{\zeta \to 0}{{\alpha \le \arg(\zeta)\le \beta}}}\Vert T(\zeta)\Vert_{{\mathcal M}(\ma_T)}<+\infty$ for $a<\alpha < \beta<b.$

\end{itemize}

The generator of the holomorphic semigroup $T$ is interpreted as in \cite{cep0} as a quasimultiplier on the closed subalgebra of $\ma$ generated by the semigroup, which is equal to the Arveson ideal $\mathcal I_{T_0}$ where $T_0$ denotes the restriction of $T$ to the half-line $(0, e^{i{a+b\over 2}}.\infty),$ and the resolvent $\zeta \to (\Delta_T -\zeta I)^{-1}$, which is defined and holomorphic outside a closed sector of the form $z +\overline S_{-ie^{ia}, ie^{ib}}$ is studied in section 7.

Consider again $a=(a_1,\dots,a_k)\in \R^k, b=(b_1,\dots,b_k)\in \R^k)$ satisfying $a_j\le b_j \le a_j+\pi$ for $j\le k$ and a finite family $T=(T_1,\dots,T_k)$ of commuting semigroups of multipliers on $\ma$ satisfying the conditions given above. By iterating the normalization process of $\ma$  with respect to $T_1,\dots, T_k$ given in sections 3 and 4, we obtain a "normalization" of $\ma$ with respect to the family $T,$ see definition 8.1, which is a commutative Banach algebra $\mb\subset \qm_r(\ma)$ for which the injection $j:\ma \to \mb$ is norm-decreasing, has dense range and extends to a norm-decreasing homomorphism from $\mathcal M(\ma)$ into ${\mathcal M}(\mb)$, for which the natural embedding $\tilde j: \qm(\ma)\to \qm(\mb)$ is onto, and for which $\lim \sup_{t\to 0^+}\Vert T(te^{ia_j}\Vert_{\mathcal M(\mb)}<+\infty$ if $a_j=b_j,$ and for which $\lim \sup_{\stackrel{\zeta \to 0}{_{\alpha_j\le arg(\zeta) \le\beta_j}}}\Vert T(\zeta)\Vert<+\infty$ for $a_j<\alpha_j\le\beta_j<b_j$ if $a_j<b_j.$

Denote by $M_{a,b}$ the set of all pairs $(\alpha, \beta)\in \R^k\times \R^k$ such that $\alpha_j=\beta_j=a_j$ if $a_j=b_j$ and such that $a_j<\alpha_j\le \beta_j<b_j$ if $a_j<b_j.$
 Let $W_{a,b}$ be the algebra of continuous functions $f$ on $\cup_{(\alpha, \beta)\in M_{a,b}}\overline S_{\alpha, \beta}:=\Pi_{j\le k}\overline S_{\alpha_j,\beta_j}$such that $e_z(\zeta)f(\zeta) \to 0$ as $\vert \zeta \vert \to 0$ in $\Pi_{j\le k}\overline S_{\alpha_j,\beta_j}$ for every $z =(z_1,\dots, z_k) \in \C^k$ and every $(\alpha, \beta)\in M_{a,b},$ and such that the maps $\zeta \to f(\zeta_1, \zeta_{j-1},\zeta, \zeta_{j+1},\dots \zeta_k)$ are holomorphic on $S_{a_j,b_j}$  if $a_j<b_j.$ For every element $\phi$ of the dual space $\mathcal G_{a,b}=W_{a,b}'$ there exists $(\alpha, \beta)\in M_{a,b} ,$ $z \in \C^k$ and a measure $\nu$ of bounded variation on $\overline S_{\alpha, \beta}:=\Pi_{j\le k}\overline S_{\alpha_j,\beta_j}$ such that
 
 $$<f,\phi>= \int_{\overline S_{\alpha,\beta}}e^{z\zeta}f(\zeta)d\nu(\zeta) \ \ \ (f\in W_{a,b}),$$
 
 and this formula allows to extend the action of $\phi$ to $e_{-z}\mathcal V_{\alpha, \beta}(X)\supset e_{-z}\mathcal U_{\alpha, \beta}(X),$ where $X$ denotes a separable Banach space and where $\mathcal U_{\alpha,\beta}(X)$ (resp. $\mathcal V_{\alpha,\beta}(X)$) denotes the algebra of continous functions $f: \overline S_{\alpha, \beta}\to X$ which converge to 0  as $\zeta \to \infty$ (resp. bounded continuous functions $f: S_{\alpha, \beta} \to X$) such that the maps $\zeta \to f(\zeta_1, \zeta_{j-1}, \zeta, \zeta_{j+1}, \dots , \zeta)$ are holomorphic on $S_{\alpha_j, \beta_j}$ when $\alpha_j < \beta_j.$

 Set ${\mathcal U}_{\alpha, \beta}:= {\mathcal U}_{\alpha, \beta}(\C).$ We describe in appendix 1 some certainly well-known ways to implement the action of $\mathcal U_{\alpha, \beta}'$ on  $\mathcal V_{\alpha, \beta}(X)$ when $(\alpha, \beta) \in M_{a,b}$ by using Cauchy transforms and Fourier-Borel transforms, and these formulae are extended to the action of elements of $(e_{-z}{\mathcal U}_{\alpha, \beta})'$ to spaces $e_{-z}\mathcal V_{\alpha,\beta}(X)$ in appendix 2.
 
 If $\phi \in (\cap_{z \in \C^k}e_{-z}\mathcal U_{\alpha,\beta})',$ define the domain $Dom(\fb(\phi))$ of the Fourier-Borel transform $\fb(\phi)$ of $\phi$ to be the set of all $z \in \C^k$ such that $\phi \in (e_{-z}\mathcal U_{\alpha, \beta})',$ and set $\fb(\phi)(z)=<e_{-z}, \phi>$ for $z \in Dom(\fb(\phi)).$ One can also define in a natural way the Fourier-Borel transform of
 $f \in e^{-z}\mathcal V_{\alpha, \beta}(X).$ Let $\lambda \in \cup_{(\gamma, \delta)\in M_{a-\alpha, b-\beta}}\overline S_{\gamma, \delta},$ and set $T_{(\lambda)}(\zeta)=T(\lambda_1\zeta_1, \dots,\lambda_k \zeta_k)$  for $\zeta \in S_{\alpha, \beta}$.   If $\lim \sup_{\stackrel{\vert  \zeta \vert \to +\infty}{\zeta \in \tilde \partial \overline S_{\alpha, \beta}}} \vert e^{-z\zeta}\vert \Vert T(\lambda)(\zeta)\Vert<+\infty,$ where $\tilde \partial \scab$ denotes the "distinguished boundary" of $\sab,$ then
 $$\sup_{\zeta \in \scab}\left \Vert e^{z\zeta}T_1(\lambda_1\zeta_1)\dots T_k(\lambda_k\zeta_k)\right \Vert_{ \mathcal M(\mb)}<+\infty,$$
 
 and one can define the action of $\phi$ on $T_{(\lambda)}$ by using the formula

$$<T_{(\lambda)}, \phi>u=<T_{(\lambda)}u, \phi>=\int_{\overline S_{\alpha,\beta}}e^{z\zeta}T_{(\lambda)}(\zeta)ud\nu(\zeta) \ \ \ (u\in \mathcal B),$$

where $\nu$ is a representing measure for $\phi e_{-z}: f \to <e_{-z}f ,\phi>\ \ (f\in \mathcal U_{\alpha, \beta}).$

 Then $<T,\phi>\in \mathcal M(\mathcal B)\subset \qm_r(\ma).$ 

 The Fourier-Borel transform of $e_zT_{(\lambda)}$ takes values in $\mathcal M(\mathcal B)$ and extends analytically to $-Res_{ar}(\Delta_{T_{(\lambda)}}):=\Pi_{1\le j \le k}(\C \setminus \sigma_{ar}(-\lambda_j\Delta_{T_j}),$ which gives the formula

$$\fb(e_zT_{(\lambda)})(\zeta)=(-1)^k\Pi_{1\le j \le k}(\lambda_1\Delta_1+(z_1+\zeta_1) I)^{-1}\dots (\lambda_k\Delta_k-(z_k+\zeta_k I)^{-1}.$$

Set $S^*_{\alpha,\beta} =\Pi_{j\le k}S_{-{\pi\over 2}-\alpha_j, {\pi\over 2}-\beta_j}.$ and  set $W_n(\zeta)=\Pi_{1\le j \le k} {n^2\over\left ( n+\zeta_je^{i{\alpha_j+\beta_j\over 2}}\right )^2}$ for $n\ge 1, \zeta=(\zeta_1,\dots, \zeta_n)\in \overline S^*_{\alpha ,\beta}.$  The results of section 2 give for $u \in \mb,$ if $z \in Dom(\phi),$ and if $\lim \sup_{\stackrel{\vert  \zeta \vert \to +\infty}{\zeta \in \tilde \partial \overline S_{\alpha, \beta}}} \vert e^{-z\zeta}\vert \Vert T(\lambda)(\zeta)\Vert<+\infty,$ where  $\tilde \partial \overline S_{\alpha, \beta}$ denotes the ''distinguished boundary" of $\scab,$

$$<T_{(\lambda)},\phi>u$$ $$=\lim_{\stackrel{\epsilon \to 0}{_{\epsilon \in S^*_{\alpha,\beta}}}}\left (\lim_{n\to +\infty}{(-1)^k\over 2i\pi)^k}\int_{z+\tilde \partial S^*_{\alpha, \beta}}W_n(\sigma -z)\fb(\phi)(\sigma)((\sigma_1-\epsilon_1)I +\lambda_1\Delta_T)^{-1}\dots ((\sigma_k-\epsilon_k)I +\lambda_k \Delta_{T_k})^{-1}ud\sigma\right ).$$

where $\tilde \partial  S^*_{\alpha,\beta}:= \Pi_{1\le j \le k}\partial \overline S_{\alpha_j, \beta_j}$ denotes the "distinguished boundary" of $\overline S^*_{\alpha, \beta},$
and where $\partial  S^*_{\alpha_j, \beta_j}$ is oriented from $-ie^{i\alpha_j}.\infty$ to $ie^{i\beta_j}.\infty.$

If, further,  $\int_{z+\tilde \partial S^*_{\alpha, \beta}}\Vert \fb(\phi)(\sigma)\Vert \vert d\sigma \vert <+\infty,$ then we have, for $u\in \mb,$

$$<T_{\lambda)},\phi>u= \lim_{\stackrel{\epsilon \to (0,\dots,0)}{_{ \epsilon \in \ssab}}}<e_{-\epsilon}T_{(\lambda)},\phi>u=$$ $$\lim_{\stackrel{\epsilon \to (0,\dots,0)}{_{ \epsilon \in \ssab}}}{(-1)^k\over (2i\pi)^k}\int_{z+\tilde \partial  S^*_{\alpha, \beta}} \FB(\phi)(\sigma)((\sigma_1-\epsilon_1)I +\lambda_1\Delta_T)^{-1}\dots ((\sigma_k-\epsilon_k)I +\lambda_k \att)^{-1}ud\sigma.$$

Finally, if $z \in Dom(\phi),$ if  $\int_{z+\tilde \partial S^*_{\alpha, \beta}}\Vert \fb(\phi)(\sigma)\Vert \vert d\sigma \vert <+\infty,$ and if $\lim \sup_{\stackrel{\vert  \zeta \vert \to +\infty}{\zeta \in \tilde \partial \overline S_{\alpha, \beta}}} \vert e^{-z\zeta}\vert \Vert T(\lambda)(\zeta)\Vert=0,$ then we have, for $u\in \mb,$

$$<T_{\lambda)},\phi>u=  {(-1)^k\over (2i\pi)^k}\int_{z+\tilde \partial  S^*_{\alpha, \beta}} \FB(\phi)(\sigma)(\sigma_1I +\lambda_1\Delta_T)^{-1}\dots (\sigma_k I +\lambda_k \att)^{-1}ud\sigma.$$ 

The convolution product of two elements of $(e_{-z}\mathcal U_{\alpha, \beta})'$  may be defined in a natural way, and if $\lambda, \phi_1,\phi_2$ satisfies the conditions above we have

$$<T_{(\lambda)},\phi_1*\phi_2>=<T_{(\lambda)},\phi_1><T_{(\lambda)},\phi_2>,$$

but there is no direct extension of this formula to the convolution product of two arbitrary elements of $\mathcal G_{a,b},$ see the comments at the end of section 8.

In section 9 of the paper we introduce a class $\mathcal U$ of "admissible open sets" $U$, with piecewise $\mathcal C^1$-boundary,  of the form $(z+ S^*_{\alpha, \beta})\setminus K$, where $K$ is bounded and where $(\alpha, \beta)\in M_{a,b}.$ These open sets $U$ have the property that $U+\epsilon \subset U$ for $\epsilon \in \overline S^*_{\alpha, \beta}$ and that $\overline U + \epsilon
\subset  Res_{ar}(-\lambda \Delta_T)$ for some $\epsilon \in \overline S^*_{\alpha, \beta}.$ Also $R(-\lambda \Delta_T,.):=(-\lambda \Delta_{T_1} -.I)^{-1}\dots (-\lambda \Delta_{T_k} -.I)^{-1}$ is bounded on the distinguished boundary of $U+\epsilon$ for $\epsilon \in S^*_{\alpha, \beta}$ when $\vert \epsilon \vert$ is sufficiently small. Standard properties of the class $H^{(1)}(U)$  of all holomorphic functions $F$ on $U$ such that
 $\Vert F\Vert_{H^{(1)}(U)}:=\sup_{\epsilon \in S^*_{\alpha, \beta}}\int_{\epsilon +\tilde \partial U}\Vert F(\sigma \vert d\sigma\vert<+\infty$ are given in appendix 3 (when $a_j=b_j$ for $j\le k,$ this space is the usual Hardy space $H^1$ on a product of open half-planes). 
 
 The results of appendix 3 allow when an open set $U\subset \C^k$ admissible with respect to $(\alpha, \beta)\in M_{a,b}$ satisfies some more suitable admissibility conditions with respect to $T=(T_1,\dots,T_k)$ and $\lambda \in \overline S_{\gamma, \delta}$ for some $(\gamma, \delta) \in M_{a-\alpha, b-\beta}$ to define $F(-\lambda_1\Delta_{T_1},\dots,-\lambda_k\Delta_{T_k})\in \mathcal M(\mb)\subset \qm_r(\ma)$ for $F\in H^{(1)}(U)$ by using the formula
 
 $$F(-\lambda_1\Delta_{T_1},\dots,-\lambda_k\Delta_{T_k})$$ $$={1\over (2i\pi)^k}\int_{\epsilon +\overline \partial U}F(\zeta_1,\dots,\zeta_k)(\lambda_1\Delta_{T_1} +\zeta_1 I)^{-1}\dots (\lambda_1\Delta_{T_k} +\zeta_k I)^{-1}d\zeta_1\dots d\zeta_k,$$
 
 where $\tilde \partial U$ denotes the distinguished boundary of $U,$ where $\epsilon \in S^*_{\alpha, \beta}$ is choosen so that $\epsilon +U$ still satisfies the required admissibility conditions with respect to $T$ and $\lambda.$
 
 Given  $T$ and $\lambda\in \cup_{(\alpha,\beta)\in M_{a,b}}\overline S_{a-\alpha,b-\beta},$ denote by ${\mathcal W}_{T,\lambda}$ the family of all open sets $U\subset \C^k$ satsisfying these admissibility conditions with respect $T$ and $\lambda.$ Then ${\mathcal W}_{T,\lambda}$ is stable under finite intersections, $\cup_{U\in {\mathcal W}_{T, \lambda}}H^{(1)}(U)$ is stable under products and we have

$$(F_1F_2)(-\lambda_1\Delta_{T_1},\dots,-\lambda_k\Delta_{T_k})=F_1(-\lambda_1\Delta_{T_1},\dots,-\lambda_k\Delta_{T_k})F_2(-\lambda_1\Delta_{T_1},\dots,-\lambda_k\Delta_{T_k}).$$

This homomorphism extends in a natural way to a bounded algebra homomorphism from $\cup_{U\in \mathcal W_{T,\lambda}}H^{\infty}(U)$
 into $\qm_r(\mb)=\qm_r(\ma),$ and we have, if  $\phi \in \mathcal F_{\alpha, \beta}$ for some $(\alpha, \beta)\in M_{a,b}$ such that $\lambda \in \overline S_{\gamma, \delta}$ for some $(\gamma, \delta)\in M_{a-\alpha, b-\beta},$ and if $\lim_{\stackrel {\vert \zeta \vert \to +\infty}{_{\zeta \in \tilde \partial \overline S_{\alpha,\beta}}}}\Vert e^{-z\zeta}T_1(\lambda_1 \zeta_1)\dots T_k(\lambda_k \zeta_k)\Vert=0$ for some $z\in Dom(\fb(\phi)),$

 $$\mathcal{FB}(\phi)(-\lambda_1\Delta_{T_1},\dots,-\lambda_k\Delta_{T_k})=<T_{(\lambda)},\phi>,$$
 
 so that  $F(-\lambda_1\Delta_{T_1},\dots, -\lambda_k\Delta_{T_k})=T(\nu \lambda_j)$ if $F(\zeta)=e^{-\nu \zeta_j},$ where $\nu \lambda_j$ is in the domain of definition of $T_j.$

 A function $F \in H^{\infty}(U)$ will be said to be strongly outer if there exists a sequence $(F_n)_{n\ge 1}$ of invertible elements of $H^{(\infty)}(U)$ such that $\vert F(\zeta)\vert \le \vert F_n(\zeta)\vert$ and $\lim_{n\to +\infty}F(\zeta)F_n^{-1}(\zeta)=1$ for $\zeta \in U.$ If $U$ is admissible with respect to some $(\alpha, \beta)\in M_{a,b}$ then there is a conformal map $\theta$ from $\D^k$ onto $U$ and the map $F\to F\circ \theta$ is a bijection from the set of strongly outer bounded functions on $U$ onto the set of strongly outer bounded functions on $\D^k.$ Every bounded outer function on the open unit disc $\D$ is strongly outer, but the class of strongly outer bounded functions on $\D^k$ is smaller than the usual class of bounded outer functions on $\D^k$ if $k\ge 2.$ We then define the Smirnov class $\mathcal S(U)$ to be the class of those holomorphic functions $F$ on $U$ such that $FG\in H^{\infty}(U)$ for some strongly outer function $G \in H^{\infty}(U).$ The bounded algebra homomorphism $F\to F(-\lambda_1\Delta_{T_1},\dots, -\lambda_k\Delta_{T_k})$ from $\cup_{U\in \mathcal W_{T,\lambda}}$ into $\qm_r(\mb)=\qm_r(\ma)$ extends to a bounded homomorphism from $\cup_{U\in \mathcal W_{T,\lambda}}\mathcal S(U)$ into $\qm(\mb)=\qm(\ma).$ If $F:\zeta \to \zeta_j,$ is the $j$-th coordinate projection then of course $F(-\lambda_1\Delta_{T_1},\dots, -\lambda_k\Delta_{T_k})=\lambda_j\Delta_{T_j}.$

 The author wishes to thank Isabelle Chalendar and Jonathan Partington for valuable discussions during the preparation of this paper.

\section{Quasimultipliers on weakly cancellative commutative Banach algebras with dense principal ideals}

We will say that a  Banach algebra $\mathcal A$ is {weakly cancellative} if $u{\mathcal A}\neq \{0\}$ for every $u \in {\mathcal A}\setminus \{0\}.$ In the whole paper we will consider weakly cancellative commutative Banach algebras with dense principal ideals, i.e. weakly cancellative commutative Banach algebras  such that the set  $\Omega({\mathcal A}):=\{ u \in {\mathcal A} \ | \ \left [ u{\mathcal A}\right ]^-={\mathcal A}\}$ is not empty.

 A quasimultiplier on such an algebra ${\mathcal A}\neq \{0\}$ is a closed operator $S=S_{u/v}:{\mathcal D}_S\to {\mathcal A},$ where $u \in {\mathcal A},$ $v\in \Omega({\mathcal A}),$ where ${\mathcal D}_S:=\{ x \in {\mathcal A} \ | \ ux \in v{\mathcal A}\},$ and where $Sx$ is the unique $y \in {\mathcal A}$ such that $vy=ux$ for $x \in {\mathcal D}_{S}.$ Let $\qm ({\mathcal  A})$ be the algebra of all quasimultipliers on ${\mathcal A}.$ A set $U \subset\qm({\mathcal A})$ is said to be pseudobounded if $\sup_{S \in U}\Vert Su\Vert<+\infty$ for some $u \in \Omega({\mathcal A})\cap \left ( \cap_{S\in U}{\mathcal D}(S)\right )$, and a quasimultiplier $S \in \qm({\mathcal A})$ is said to be regular if the family $\{\lambda^nS^n\}_{n \ge 1}$ is pseudobounded for some $\lambda >0.$ The algebra of all regular quasimultipliers on ${\mathcal A}$ will be denoted by $\qm_r({\mathcal A}).$ A multiplier on $
{\mathcal A}$ is a bounded linear operator $S$ on ${\mathcal A}$ such that $S(uv)=(Su)v$ for $u \in {\mathcal A}, v\in {\mathcal A},$ and the multiplier algebra ${\mathcal M}({\mathcal A})$ of all multipliers on ${\mathcal A},$ which is a closed subalgebra of ${\mathcal B}({\mathcal A}),$ is also the algebra of all quasimultipliers on ${\mathcal A}$ such that ${\mathcal D}_S={\mathcal A},$ and  ${\mathcal M}({\mathcal A})\subset \qm_r({\mathcal A}).$ Also if $S=S_{u/v}\in \qm(\ma), w \in {\mathcal D}(S), R \in {\mathcal M}(\ma), $ then $u(Rw)=R( v (Sw))=v(R(Sw)),$ so $Rw\in {\mathcal D}(S),$ and we have

\begin{equation} R(Sw)=S(Rw).\end{equation}

If $\ma$ is unital then $\Omega(\ma)=\mga,$ where $\mga$ denotes the group of invertible elements of $\ma,$ and $\QM(\ma)={\mathcal M}(\ma).$

The following notion if slightly more flexible than the notion of $s$-homomorphism introduced by the author in \cite{e2}.

\begin{defn}: Let ${\mathcal A}$ be a weakly cancellative commutative Banach algebra with dense principal ideals, and let $\mb$ be a weakly cancellative Banach algebra. A homomorphism $\Phi: \ma \to \mb$ is said to be a $\qm$-homomorphism if the following conditions are satisfied

(i) $\Phi$ is one-to-one, and $\Phi(\ma)$ is dense in $\mb.$

(ii)  $\Phi(u)\mb\subset \Phi(\ma)$ for some $u\in \Omega(\ma).$

\end{defn}

If the conditions of definition 2.1 are satisfied, we will say that $\Phi$ is a $\qm$-homomorphism with respect to $u.$ Notice that $\Phi(u)\in \Omega(\mb),$ and so the existence of such an homomorphism implies that $\mb$ is a weakly cancellative commutative Banach algebra with dense principal ideals. Notice also that condition (ii) shows that $\mb$ may be identified to a subalgebra of $\qm(\Phi(\ma))\approx \qm(\ma).$  

\begin{prop} Let $\Phi:\ma \to \mb$ be a homomorphism between weakly cancellative commutative Banach algebras with dense principal ideals, and assume that $\Phi$ is a $\qm$-homomorphism with respect to some $u_0 \in \Omega(\ma).$

(i) There exists $M>0$  such that $\Vert \Phi^{-1}(\Phi(u_0)v)\Vert \le M\Vert v \Vert$ for $v\in \mb.$

(ii) $\Phi^{-1}(\Phi(u_0)v) \in \Omega(\ma)$ for $v \in \Omega(\mb).$

(iii) Set $\tilde \Phi(S_{u/v})=S_{\Phi(u)/ \Phi(v)}$ for $S_{u/v}\in \qm(\ma).$ Then $\tilde \Phi:\qm(\ma)\to \qm(\mb)$ is a pseudobounded isomorphism, and $\tilde \Phi^{-1}(S_{u/v})
=S_{\Phi^{-1}(\Phi(u_0)u)/\Phi^{-1}(\Phi(u_0)v)}$ for $S_{u/v}\in \qm(\mb).$ 

\end{prop}

Proof: (i) Set $\Psi(v)=\Phi^{-1}(\Phi(u_0)v)$ for $v \in \mb.$ If $\lim_{n\to +\infty}v_n=v \in \mb,$ and if $\lim_{n\to +\infty}\psi(v_n)=w\in \ma,$ then $\Phi(u_0)v=\Phi(w),$ so that $w=\Psi(v)$ and (i) follows from the closed graph theorem. 

(ii) Let  $v\in \Omega(\mb),$ and let $(w_n)_{n\ge 1}$ be a sequence of elements of $\ma$ such that $\lim_{n\to +\infty} v\Phi(w_n)=\Phi(u_0).$ Then $\lim_{n\to +\infty}\Phi^{-1}(\Phi(u_0)v)w_n=u_0^2\in \Omega(\ma),$ and so $\Phi^{-1}(\Phi(u_0)v)\in \Omega(\ma).$

(iii) Let $U \subset \qm(\ma)$ be a pseudobounded set, and let $w\in \Omega(\ma)\cap\left (\cap_{S\in U}{\mathcal D}(S)\right )$ be such that 
$\sup_{S\in U}\Vert Sw\Vert<+\infty.$ Then $\Phi(w)\in \Omega(\mb),$ and $\Phi(w)\in \cap_{S\in U}{\mathcal D}(\tilde \phi(S)).$ Since $\sup_{S\in U}\Vert \tilde \Phi(S)\Phi(w)\Vert \le \Vert \Phi \Vert \sup_{S\in U}\Vert Sw\Vert <+\infty,$ this shows that $\tilde \Phi: \QM(\ma)\to \QM(\mb)$ is pseudobounded.

Now set $\theta(S_{u/v})
=S_{\Phi^{-1}(\Phi(u_0)u)/\Phi^{-1}(\Phi(u_0)v)}=S_{\Psi(u)/\Psi(v)}$ for $S_{u/v}\in \qm(\mb).$ It follows from (ii) that $\theta: \qm(\mb)\to \qm(\ma)$ is well-defined. Let $U\subset \qm(\mb)$ be pseudobounded, and let $w \in \Omega(\mb)\cap\left ( \cap_{S\in U}{\mathcal D}(S)\right )$ be such that $\sup_{S\in U}\Vert Sw\Vert <+\infty.$ We have, for $S=S_{u/v}\in U,$

$$u(S_{u/v}w)=vw, u\Phi(u_0)(S_{u/v}w\Phi(u_0))=v\Phi(u_0)w\Phi(u_0),$$ $$ \Phi^{-1} (u\Phi(u_0))\left (\Phi^{-1}(S_{u/v}w\Phi(u_0) \right )=\Phi^{-1}(v\Phi(u_0))\Phi^{-1}(w\Phi(u_0)).$$
So $\Phi^{-1}(w\phi(u_0))\in {\mathcal D}(\theta(S)),$ and $\theta(S)\Phi^{-1}(w\Phi(u_0))=\phi^{-1}\left ( S_{u/v}w\Phi(u_0)\right ).$ Since $\sup_{S\in U}\Vert \left ( Sw\Phi(u_0)\right )\Vert \le M\sup_{S\in U}\Vert Sw\Vert <+\infty,$ this shows that $\theta: \qm(\mb)\to \qm(\ma)$ is pseudobounded. We have, for $S=S_{u/v}\in \qm(\ma),$

$$(\theta\circ \tilde \Phi)(S)=S_{\Phi^{-1}(\Phi(u)\Phi(u_0))/\Phi^{-1}(\Phi(v)\Phi(u_0))}=S_{uu_0/vu_0}=S_{u/v}=S.$$

We also have, for $S=S_{u/v}\in \qm(\ma),$

$$(\tilde \Phi\circ \theta)(S_{u/v})=S_{\Phi(\Phi^{-1}(u\Phi(u_0))/\Phi(\Phi^{-1}(v\Phi(u_0))}=S_{u\Phi(u_0)/v\Phi(u_0)}=S_{u/v}=S.$$

Hence $\tilde \Phi=\qm(\ma)\to \qm(\mb)$ is bijective, and $\tilde \Phi^{-1}=\theta:\qm(\mb)\to \qm(\ma)$ is pseudobounded. $\square$

The following result is a simplified version of theorem 7.11 of \cite{e2}, which will be used in the next two sections.

\begin{prop} Let $\ma$ be a weakly cancellative commutative Banach algebra with dense principal ideals, and let $U\subset \qm(\ma)$ be a pesudobounded set stable under products.
 Set ${\mathcal L}:=\{u \in \cap_{S\in U}{\mathcal D}(S) \ | \  \sup_{S\in { U}}\Vert Su\Vert <+\infty\},$ and set $\Vert u \Vert_{\mathcal L}:=\sup_{S\in U \cup\{I\}}\Vert Su\Vert$ for $u\in {\mathcal L}.$ Then $(\mathcal L, \Vert .\Vert_{\mathcal L})$ is a Banach algebra, $Ru\in {\mathcal L}$ and $\Vert Ru\Vert_{\mathcal L}\le \Vert R \Vert_{{\mathcal M}(\ma)}\Vert u \Vert _{\mathcal L}$ for $R\in {\mathcal M}(\ma),$ $u \in {\mathcal L},$ and if we denote by $\mathcal B$ the closure of $\ma$ in $({\mathcal M}(\mathcal L), \Vert .\Vert_{{\mathcal M}(\mathcal L)}),$ then the following properties hold

(i) $\mathcal B$ is a weakly cancellative commutative Banach algebra, and the inclusion map $j: \ma \to {\mathcal B}$ is a $\qm$-homomorphism with respect to $w$ for $w\in \mathcal(\ma)\cap \mathcal L.$

(ii) $\tilde j(\mathcal M(\ma)) \subset \mathcal M(\mb)$, and $\Vert \tilde j(R)\Vert_{{\mathcal M}(\mb)}\le \Vert_{\mathcal M(\ma)}$ for $R \in \mathcal M(\ma),$ where $\tilde j:\qm(\ma)\to \qm(\mb)$ is the pseudobounded isomorphism associated to $j$ in proposition 2.2(iii). 

(iii) $S \in \mathcal M(\mathcal B),$ and $\Vert S \Vert_{\mathcal M(\mathcal B)}\le \Vert S \Vert_{\mathcal M(\mathcal L)} \le 1$ for every $S \in U.$

\end{prop}

Proof: The fact that $(\mathcal L, \Vert .\Vert_{\mathcal L})$ is a Banach space follows from a standard argument given in the proof of theorem 7.11 of \cite{e2}. Clearly, ${\mathcal L}$ is an ideal of $\ma,$ and it follows from the definition of $\Vert .\Vert_{\mathcal L}$ that $\Vert u \Vert \le \Vert u\Vert _{\mathcal L}$ for $u\in \mathcal L.$ We have, for $u \in {\mathcal L}, v\in {\mathcal L},$

$$\Vert uv \Vert_{\mathcal L}\le \Vert u \Vert_{\mathcal L}\Vert v \Vert\le \Vert u \Vert_{\mathcal L}\Vert v \Vert_{\mathcal L},$$

and so $(\mathcal L, \Vert .\Vert_{\mathcal L})$ is a Banach algebra. If $R \in {\mathcal M}(\ma),$ $u\in \mathcal L,$ then it follows from (i) that $Ru\in \cap_{S\in U}{\mathcal D}(S),$ and that we have

$$\sup_{S\in U\cup\{I\}}\Vert S(Ru)\Vert=\sup_{S\in U\cup\{I\}}\Vert R(Su)\Vert\le \Vert R\Vert_{\mathcal M (\ma)}\Vert u \Vert_{\mathcal L}.$$

Hence $Ru \in \mathcal L,$ and $\Vert Ru\Vert_{\mathcal L}\le \Vert R \Vert_{{\mathcal M}(\ma)}\Vert u \Vert _{\mathcal L}.$

Now denote by $\mathcal B$ the closure of $\ma$ in $({\mathcal M}(\mathcal L), \Vert .\Vert_{{\mathcal M}(\mathcal L)}),$ and let $w \in \Omega(\ma)\cap \mathcal L.$ Since $\mathcal L \subset \mathcal B,$ $\mathcal B$ is weakly cancellative, and $\mathcal B$ is commutative and has dense principal ideals since $j(\ma)$ is dense in $\mathcal B.$ Since $w\mathcal B  \subset \mathcal L \subset \ma,$ we see that the inclusion map $j: \ma \to {\mathcal B}$ is a $\qm$-homomorphism with respect to $w,$ which proves (i).

Let $R\in \mathcal M(\ma),$ and denote by $R_1$ the restriction of $R$ to $\mathcal L.$ Then $R_1\in \mathcal M(\mathcal L),$ and $\Vert R_1\Vert_{\mathcal M(\mathcal L)}\le \Vert R \Vert_{\mathcal M(\ma)}.$ Set $R_2u=R_1u$ for $u\in \mb.$ Then $R_2u\in \mb,$ and $\Vert R_2u\Vert_{\mathcal M(\mb)} \le \Vert R_1\Vert_{\mathcal M(\mathcal L)} \Vert u\Vert_{\mathcal M(\mathcal L)}\le \Vert R\Vert_{\mathcal M(\mathcal \ma)}\Vert \Vert u\Vert_{\mathcal M(\mathcal L)}.$ Hence $R_2\in \mathcal M(\mb),$ and $\Vert R_2\Vert_{\mathcal M(\mb)}\le 
\Vert R\Vert_{\mathcal M(\ma)}.$

Now let $S_0\in U.$ we have, for $u \in \mathcal L,$ since $U$ is stable under products,

$$\Vert S_0u\Vert_{\mathcal L}=\sup_{S\in U}\Vert S_0Su\Vert\le \sup_{S\in U}\Vert Su\Vert=\Vert u\Vert_{\mathcal L},$$

and so $S_0u\in \mathcal L,$ and $\Vert S_0\Vert_{{\mathcal M}(\mathcal L)}\le 1.$ This implies that $S_0(\mathcal B)\subset \mathcal B,$ so that $S_0\in \mathcal M(\mathcal B),$
and $\Vert S \Vert_{\mathcal M(\mathcal B)}\le \Vert S_0\Vert_{{\mathcal M}(\mathcal L)}\le 1,$ which proves (iii). $\square$

We have the following very easy observation.

\begin{prop}: Let $\ma_0,\ma_1$ and $\ma_2$ be weakly cancellative commutative Banach algebras, and assume that $\Phi_0:\ma_0\to \ma_1$ is a $\qm$-homomorphism with respect to $u_0\in \Omega(\ma_0)$ and that $\Phi_1:\ma_1\to \ma_2$ is a $\qm$-homomorphism with respect to $u_1\in \Omega(\ma_1).$ Then $\Phi_1\circ \Phi_0: \ma_0\to \ma_2$
is a $\qm$-homomorphism with respect to $\Phi_0^{-1}(\Phi_0(u_0)u_1).$

\end{prop}

Proof: The homomorphism $\Phi_1\circ \Phi_0$ is one-to-one and has dense range, and it follows from proposition 2.2 (ii) that $\Phi_0^{-1}(\Phi_0(u_0)u_1) \in \Omega(\ma_0).$ Let $u\in \ma_2,$ let $v\in \ma_2$ be such that $\Phi_1(v)= \Phi_1(u_1)u,$ and let $w\in \ma_0$ be such that $\Phi_0(w)=\Phi_0(u_0)v.$ Then $(\Phi_1\circ \Phi_0)\left (\Phi_0^{-1}(\Phi_0(u_0)u_1) \right)u=(\Phi_1\circ \Phi_0)(u_0)\Phi_1(u_1)u=\Phi_1(\Phi_0(u_0)v)=(\Phi_1\circ \Phi_0 )(w)\subset (\Phi_1\circ \Phi_0)(\ma_0),$ and so $\Phi_1\circ \Phi_0$ is a $\qm$-homomorphism with respect to $\Phi_0^{-1}(\Phi_0(u_0)u_1).$ $\square$

 We will denote by $\widehat A$ the space of all characters on a commutative Banach algebra $\ma,$ equipped with the Gelfand topology. Recall that $\ma$ is said to be radical when $\widehat A=\emptyset.$ 

\begin{defn}: Let $\ma$ be a weakly cancellative commutative Banach algebra with dense principal ideals.  For $\chi \in \widehat {\ma},$ $S=S_{u/v}\in \qm(\ma),$ set $\tilde \chi(S)={\chi(u)\over \chi(v)},$ and for $S\in \qm(\ma),$ set $\sigma_{\ma}(S):=\{\tilde \chi(S)\}_{\chi \in \widehat {\ma}},$ with the convention $\sigma_{\ma}(S)=\emptyset$ if $\ma$ is radical.

\end{defn}

Clearly, $\tilde \chi$ is a character on $\qm(\ma)$ for $\chi \in \widehat{\ma},$ and the map $\chi \to \tilde \chi(S)$ is continuous on $\widehat \ma$ for $S\in \qm(\ma).$

We will use the following result in the study of the resolvent of semigroups.

\begin{prop} Let $\ma$ be a weakly cancellative commutative Banach algebra with dense principal ideals, and let $S\in \qm(\ma).$ If  $\lambda_0 -S$  has an inverse $(\lambda_0 I-S)^{-1}$ in $\qm(\ma)$ which belongs to $\ma$ for some $\lambda_0 \in \C,$ where $I$ denotes the unit element of ${\mathcal M}(\ma),$ then $\sigma_{\ma}(S)$ is closed, $\lambda I -S$ has an inverse $(\lambda I-S)^{-1}$ in $\qm(\ma)$ which belongs to $\ma$ for every $\lambda \in  \C \setminus
\sigma_{\ma}(S),$ and the $\ma$-valued map $\lambda \to (\lambda I-S)^{-1}$ is holomorphic on $ \C \setminus \sigma_{\ma}(S).$
\end{prop}

Proof: If $\ma$ is unital, then $\qm(\ma)=\ma,$ and there is nothing to prove. So assume that $\ma$ is not unital, and set $\tilde \ma:=\ma\oplus \C I.$ Then $\widehat {\tilde \ma}=\{{\chi_0}\}\cup\{\tilde \chi_{_{|{\tilde \ma}}}\}_{\chi \in \widehat{\ma}},$ where $\chi_0(a+\lambda I)=\lambda$ for $a \in \ma, \lambda \in \C.$


Set $a=(\lambda_0 I-S)^{-1}\in \ma.$ Then $\chi(a) =\tilde \chi (a)={1\over \lambda_0 -\tilde \chi(S)},$ so that $\tilde \chi(S)=\lambda_0-{1\over \chi(a)}$ for $\chi \in \widehat{\ma}.$ Since $\sigma_{\ma}(a)\cup\{0\}=\sigma_{\tilde {\ma}}$ is a compact subset of $\C,$ this shows that $\sigma_{\ma}(S)$ is closed.  We have, for $\lambda \in \C,$

$$spec_{\tilde{\ma}}(I+(\lambda -\lambda_0)a)=\{1 \}\cup \left \{ \frac{\lambda -\tilde \chi(S)}{\lambda_0 -\tilde \chi(S)}\right \}_{\chi \in \widehat{\ma}}.$$

So $I+(\lambda  -\lambda_0)a$ is invertible in $\tilde{\ma}$ for $\lambda \in \C \setminus \sigma_{\ma}(S),$ and the map $\lambda \to a(I+(\lambda -\lambda_0)a)^{-1}\in \ma$ is holomorphic on $\C \setminus \sigma_{\ma}(S).$ We have, for $\lambda \in \C \setminus \sigma_{\ma}(S),$

$$(\lambda I -S)a(I+(\lambda-\lambda_0)a)^{-1}=((\lambda-\lambda_0)I +(\lambda_0 -S))a(I+(\lambda-\lambda_0)a)^{-1}$$ $$=(I +(\lambda-\lambda_0)a)(I+(\lambda-\lambda_0)a)^{-1}=I.$$

Hence $\lambda I-S$ has an inverse $(\lambda I -S)^{-1}\in \ma$ for $\lambda \in \C -\sigma_{\ma}(S),$ and the map $\lambda \to (\lambda I-S)^{-1}=a(I+(\lambda-\lambda_0)a)^{-1}$ is holomorphic on $\C \setminus \sigma_{\ma}(S)$. $\square$

\section{Normalization of a commutative Banach algebra with respect to a strongly continuous semigroup of multipliers}

A semigroup $T=(T(t))_{t>0}$ of multipliers on a commutative Banach algebra $\ma$ is said to be strongly continuous if the map $t \to T(t)u$ is continuous on $(0,+\infty)$ for every $u\in \ma.$ This implies that $\sup_{\alpha \le t \le \beta}\Vert T(t)\Vert <+\infty$ for $0<\alpha\le \beta <+\infty,$ and so $\Vert T(t)\Vert ^{1\over t}$ has a limit $\rho_T$ as $t\to +\infty,$ and $\rho_T=\lim_{n\to +\infty}\Vert T(n)\Vert ^{1\over n}.$ 
In the remainder of the section $T=(T(t))_{t>0}$ will denote a strongly continuous group of multipliers on a weakly cancellative commutative Banach algebra $\ma$ such that $\cup_{t>0}T(t)\ma$ is dense in $\ma.$ Hence if $u\in \cap_{t>0}Ker(T(t)),$ then $uv=0$ for every $v \in \cup_{t>0}T(t)\ma,$ so $u\ma=\{0\}$ and $u=0.$ 

Notice that in this situation if $\ma$ has a unit element $1$ then if we set $\tilde T(t)=T(t).1$ then $\tilde T:=(\tilde T(t))_{t>0}$ is a norm-continuous semigroup of elements of $\ma.$
Since $\cup_{t>0}T(t)\ma$ is dense in $\ma,$ $T(t_0).1$ is invertible in $\ma$ for some $t_0>0,$ and so $\lim_{t\to 0^+}\Vert T(t).1-1\Vert=0,$ which implies that the generator of $\tilde T$ is bounded. So there exist $R\in {\mathcal M}(\ma)\approx \ma$ such that $T(t)=e^{tR}$ for $t >0$ if $\ma$ is unital.

Let $\omega$ be a positive measurable weight on $(0,+\infty).$ Recall that if $\omega(s+t)\le \omega(s)\omega(t)$ for $s>0,t>0,$ then the space $L^1_{\omega}(\R^+)$ of all (classes of) measurable functions on $[0,+\infty)$ satisfying $\Vert f \Vert_{\omega}:= \int_0^{+\infty}\Vert f(t)\vert \omega(t)dt<+\infty,$ equipped with the norm $\Vert .\Vert_{\omega},$ is a Banach algebra with respect to convolution. Set $\omega_T(t)=\Vert T(t)\Vert$ for $t>0.$ For $f \in L^1_{\omega_T}(\R^+),$ define $\Phi_T(f)\in {\mathcal M}(\ma)$ by the formula

\begin{equation} \Phi_T(f)u=\int_{0}^{+\infty}f(t)T(t)udt \ \ \ (u \in \ma).\end{equation}

Denote by $\iit$ the closure of $\Phi_T(L^1_{\omega_T}(\R^+))$ in ${\mathcal M}(\ma).$ Let $(f_n)_{n\ge 1}$ be a Dirac sequence, i.e. a sequence $(f_n)$ of nonnegative integrable functions on $\R^+$  such that $\int_0^{+\infty}f_n(t)dt=1$ and such that $f_n(t)=0$ a.e. on $(\alpha_n,+\infty)$ with $\lim_{n\to +\infty}\alpha_n=0.$ Since the semigroup $T$ is strongly continuous on $\ma,$ a standard argument shows that $\lim_{n\to +\infty}\Vert \Phi_T(f_n*\delta_t)u-T(t)u\Vert=0$ for every $t>0$ and every $u\in \ma.$ This shows that if $v\in \Omega(\iit),$ and if $w \in \Omega(\ma),$ then $vw \in \Omega(\ma).$

The following result is then a consequence of theorem 6.8 of \cite{e1} and of proposition 5.4 of \cite{e4}.

\begin{lem} There exists $w\in \Omega(\ma)$ such that $\lim \sup_{t\to 0^+}\Vert T(t)w\Vert <+\infty.$

\end{lem}

Proof:  Let $\lambda > log\left (\rho_{\bf T}\right )$,  and set $\omega_{\lambda}(t)=e^{\lambda t }\sup_{s\ge t}e^{-\lambda s}\Vert T(s)\Vert$ for $t>0.$ An extension to lower semicontinuous weights of theorem 6.8 of \cite{e1} given in \cite{e4}  shows that $\Omega(L^1_{\omega_{\lambda}}(\R^+))\neq \emptyset.$ It follows also from proposition 5.4 of \cite{e4}  that $\Vert \Phi_{ T}(g)T(t)\Vert \le e^{\lambda t}\Vert g \Vert_{L^1_{\omega_{\lambda}}}$ for every $g \in L^1_{\omega_{\lambda}}(\R^+)$ and every $t>0,$ and that $\Phi_{ T}(g) \subset \Omega(\iit)$ for every $g \in \Omega \left (L^1_{\omega_{\lambda}}(\R^+)\right ).$ Hence if $g\in \Omega \left (L^1_{\omega_{\lambda}}(\R^+)\right )$ and if $v \in \Omega(\ma),$ then $\Phi_{T}(g)v \in \Omega(\ma),$ and $\lim \sup_{t\to 0^+}\Vert T(t)\Phi_{T}(g)v \Vert <+\infty.$ $\square$
 
The following result is a version specific to one-parameter semigroups of proposition 2.3.
\begin{prop} Let $T:=(T(t))_{t>0}$ be a strongly continuous of multipliers on a weakly cancellative commutative Banach algebra $\ma$ with dense principal ideals such that $\cup_{t>0}T(t)\ma$ is dense in $\ma.$  Let ${\mathcal L}_{T}:=\{ u \in \ma \ | \ \lim \sup_{t\to 0^+}\Vert T(t)u\Vert <+\infty\}\supset \cup_{t>0}T(t)\ma$, choose $\lambda > log(\rho_{T}),$ set $\Vert u \Vert_{\lambda}:=sup_{s\ge 0}e^{-\lambda s}\Vert T(s)u\Vert$ for $u \in {\mathcal L}_{T},$ with the convention $T(0)=I,$ and set  $\Vert R \Vert_{\lambda,op} := \sup\{\Vert Ru\Vert _{\lambda} \ | \ u \in {\mathcal L}_{T}, \Vert u \Vert_{\lambda}\le 1\}=\Vert R \Vert_{{\mathcal M}({\mathcal L}_{T})}$ for $R\in {\mathcal M}({\mathcal L}_{T}).$ Denote by  $\at$ the closure of $\ma$ in $( {\mathcal M}(\mathcal L_T), \Vert .\Vert_{\lambda, op}).$    

Then $({\mathcal L}_{T}, \Vert .\Vert_{\lambda})$ is a Banach algebra, the norm topology on ${\mathcal L}_{T}$ does not depend on the choice of $\lambda,$  and the following properties hold
  
(i) The inclusion map $j:\ma \to \ma_T$ is a $\qm$-homomorphism with respect to $w$ for every $\omega \in \Omega(\ma)\cap {\mathcal L}_T,$ the tautological map $\tilde j: S_{u/v}\to S_{u/v}$ is a pseudobounded isomorphism from $\qm(\ma)$ onto $\qm(\ma_T)$ and if $w\in \Omega(\ma)\cap {\mathcal L}_T,$ then $\tilde j^{-1}(S)=S_{uw/vw}$ for $S=S_{u/v}\in \qm(\ma_T).$

(ii) $\tilde j(\mathcal M(\ma))\subset \mathcal M(\ma_T),$  and $\Vert R\Vert_{\mathcal M(\at)}\le \Vert R \Vert_{\mathcal M(\ma)}$ for $R\in \mathcal M(\ma).$
 
(iii) $\Vert T(t)\Vert_{\mathcal M(\ma_T}\le \Vert T(t)\Vert_{\lambda, op}\le e^{\lambda t}$ for $t>0,$ and $\lim \sup_{t\to 0^+}\Vert T(t)u-u\Vert_{\lambda, op}=0$ for every $u\in \ma_T.$
  
\end{prop}
  
 Proof:  It follows from lemma 3.1 that the family $U=\{e^{-\lambda t}T(t)\}_{t>0}$ is pseudobounded for $\lambda >log(\rho_T).$ The fact that $(\mathcal L_T, \Vert .\Vert_{\lambda})$ is a Banach algebra, and assertions (i) and (ii) follow from proposition 2.2 and proposition 2.3 applied to $U$, and an elementary argument given in the proof of theorem 7.1 of \cite{e4} shows that there exists $k>0$ and $K>0$ such that $k\Vert u\Vert_{\lambda}\le \sup_{0\le t \le 1}\Vert T(t)u\Vert \le K\Vert u\Vert_{\lambda}$ for $u \in {\mathcal L}_T,$ which shows that the norm topology on $\mathcal L_T$ does not depend on the choice of $\lambda.$ 
 
 It follows also from proposition 2.3 applied to $U$ that $\Vert T(t)\Vert_{\mathcal M(\ma_T}\le \Vert T(t)\Vert_{\lambda, op}\le e^{\lambda t}$ for $t>0,$ and $\lim_{t\to 0^+}\Vert T(t)u-u\Vert_{\lambda, op}=0$ for every $u \in \cup_{t>0}T(t)\ma.$ Since $\cup_{t>0}T(t)\ma_T$ is dense in $\ma_T,$ a standard density argument shows that $\lim_{t\to 0^+}\Vert T(t)u-u\Vert_{\lambda, op}=0$ for every $u \in \ma_T.$
$\square$

  This suggests the following definition
  
  \begin{defn}  Let $\ma$ be a weakly cancellative commutative Banach algebra with dense principal ideals, let $T=((T(t))_{t>0}$ be a strongly continuous semigroup of multipliers on $\ma$ such that $T(t)\ma$ is dense in $\ma$ for $t>0.$  A normalization $\mb$ of $\ma$ with respect to $T$ is a subalgebra $\mb$ of $\qm(\ma)$ which is a Banach algebra with respect to a norm $\Vert .\Vert_{\mb}$ and satisfies the following conditions

(i) The inclusion map $j: \ma \to \mb$ is a $\qm$-homomorphism, and $\Vert j(u)\Vert_{\mb}\le \Vert u \Vert_{\ma}$ for $u \in \ma.$

(ii) $\tilde j (R) \subset \mathcal M(\mb),$ and $\Vert \tilde j(R)\Vert _{\mathcal M(\mb)}\le  \Vert R\Vert _{\mathcal M(\ma)}$ for $R\in \mathcal M(\ma),$ where $\tilde j: \qm(\ma)\to \qm(\mb)$ is the pseudobounded isomorphism associated to $j$ introduced in proposition 2.2 (ii).

(iii) $\lim \sup_{t\to 0^+}\Vert \tilde j(T(t))\Vert_{\mathcal M(\mb)} <+\infty.$

\end{defn}

For example the algebra $\ma_T$ constructed in proposition 3.2 is a normalization of the given Banach algebra $\ma$ with respect to the semigroup $T.$ Notice that if $\mb$ is a normalization of $\ma$ with respect to $T,$ the same density argument as above shows that $\lim_{t\to 0^+}\Vert T(t)u-u\Vert_{\mb}=0$ for every $u \in \mb.$

\section{Normalization of a commutative Banach algebra with respect to a holomorphic semigroup of multipliers}

For $a<b \le a+\pi,$ denote by $S_{a,b}$ the open sector $\{ z \in \C \setminus \{0\} \ | \ a < arg( z )< b\},$ with the convention $\overline{S}_{a,a}=\{re^{ia}\ | \ r\ge 0\}.$ In this section we consider again a weakly cancellative commutative Banach algebra $\ma$ with dense principal ideals and we consider a semigroup $T=(T(\zeta))_{\zeta\in S_{a,b}}$ of multipliers on $\ma$ such that $\cup_{t\in S_{a,b}}T(t)\ma$ is dense in $\ma$ which is holomorphic  on  $S_{a,b}$, which implies that $T(\zeta)\ma$ is dense in $\ma$ for every $\zeta \in S_{a,b}.$ So $T(\zeta)u \in \Omega(\ma)$ for every $\zeta\in S_{a,b}$
and every $u \in \Omega(\ma).$ We state as a lemma the following easy observations.

\begin{lem}  Let $u \in {\mathcal M}({\mathcal A})$ such that  $\lim \sup_{\stackrel{\zeta\to 0}{_{\zeta \in \overline{\mathcal S}_{\alpha,\beta}}}}\Vert T(\zeta)u\Vert <+\infty,$ where $a<\alpha \le \beta <b.$

(i) If $\lambda  >{1\over cos\left ({\beta -\alpha\over 2}\right )}\left (\lim_{r \to +\infty} \frac{log (max(\Vert T(r e^{ i\alpha}),\Vert T(re^{i\beta}\Vert)}{r}\right )$  then \\ $\sup_{\zeta\in \scab}\Vert e^{-\lambda \zeta e^{-i{\alpha +\beta\over 2}} }T(\zeta)u\Vert <+\infty.$

 $$(ii) \lim \sup_{\lambda \to +\infty}\left [ \sup_{\zeta\in \scab}\Vert e^{-\lambda \zeta e^{-i{\alpha +\beta\over 2}} }T(\zeta)u\Vert \right ]=\lim \sup_{\stackrel{\zeta\to 0}{_{\zeta \in \scab}}}\Vert T(\zeta)u\Vert.\ \ \ \ \ \ \ \ \ \ \ \ \ \ \ \ \ \ \ \ \  \ \ \ $$
\end{lem}

Proof: Set $\gamma:={\beta -\alpha\over 2}.$ We have, for $r_1\ge 0, r_2\ge 0,$

$$\left \Vert e^{-\lambda e^{-i{\alpha +\beta\over 2}}(r_1e^{i\alpha}+r_2e^{i\beta})}T(r_1e^{i\alpha}+r_2e^{i\beta})u \right \Vert $$ $$\le \inf \left [e^{-{\lambda cos(\gamma)} r_1}\Vert T(r_1e^{i\alpha})u \Vert e^{-{\lambda cos(\gamma)} r_2}\Vert T(r_2e^{i\beta})\Vert, e^{-{\lambda cos(\gamma)} r_1}\Vert T(r_1e^{i\alpha}) \Vert e^{-{\lambda cos(\gamma)} r_2}\Vert T(r_2e^{i\beta})u\Vert \right ].$$

Set $m:=\lim \sup_{\stackrel{\zeta\to 0}{_{\zeta \in \scab}}}\Vert T(t)u\Vert,$ and let $\epsilon >0.$ There exists $\delta >0$ such that $\Vert T(\zeta)u\Vert\le m+\epsilon$ for $0\le \vert \zeta \vert \le \delta, \zeta \in \scab.$

We obtain, considering the cases $\sup(r_1,r_2)\le {\delta/2}, \inf(r_1,r_2)\le \delta/2 \ \mbox{and} \ \sup(r_1,r_2)\ge{ \delta\over 2},$ and the case where $\inf(r_1,r_2)\ge {\delta\over 2},$

$$ \sup_{\zeta \in \scab}\Vert e^{-\lambda \zeta e^{i{\alpha +\beta\over 2}}}T(\zeta)u\Vert $$ $$\le \max \left [ (m+\epsilon)\sup_{\stackrel{\vert \zeta\vert  \le \delta}{_{\zeta \in \overline S_{\alpha, \beta}}}}\vert e^{-\lambda \zeta}\vert \left (1+ \sup_{r \ge \delta/2}e^{-\lambda r cos (\gamma)}\max(\Vert T(re^{i\alpha})\Vert, \Vert T(re^{i\beta})\Vert\right ), \right .$$ $$\left . \Vert u \Vert \left ( \sup_{r \ge \delta/2}e^{-\lambda r cos (\gamma)}\Vert T(re^{i\alpha})\Vert \right ) \left (\sup_{r \ge \delta/2}e^{-\lambda rcos (\gamma)}\Vert T(re^{i\beta})\Vert\right )\right ],$$ and (i) and (ii) follow from this inequality. $\square$

We now use a construction proposed by I. Chalendar, J.R. Partington and the author in proposition 3.6 of \cite{cep0} to associate to construct a $\qm$- homomorphism from $\ma$ into a weakly cancellative commutative Banach algebra $\mb$ such that $\sup_{\stackrel{t \in \overline{\mathcal S}_{\alpha,\beta}}{_{\vert t \vert \le 1}}}\Vert T(t)\Vert<+\infty$ for every $\alpha, \beta$ satisfying $a <\alpha \le \beta < b.$

The following result is more general than proposition 3.6 of \cite{cep0}.

\begin{prop}  Let $\ma$ be a weakly cancellative commutative Banach algebra with dense principal ideals, and let $T=((T(\zeta))_{\zeta\in S_{a,b}}$ be a holomorphic semigroup of multipliers on $\ma$ such that $T(\zeta)\ma$ is dense in $\ma$ for $\zeta\in S_{a,b}.$ Set $\alpha_n={a+{b-a\over 2(n+1)}}$, $\beta_n:=b -{b-a\over 2(n+1)}$ for $n\ge 0,$ and let  ${\bf \mu}= (\mu_n)_{n\ge 0}$ be a nondecreasing sequence  of positive real numbers satisfying the following conditions

\begin{equation}\mu_n >{1\over cos\left ( {\beta_n-\alpha_n\over 2}\right )}\lim_{r \to +\infty} \frac{log (max(\Vert T(r e^{ i\alpha_n})\Vert,\Vert (T(re^{i\beta_n})\Vert)}{r} \ \ (n\ge 0).\end{equation}

\begin{equation}\sup_{\zeta \in \overline {\mathcal S}_{\alpha_n,\beta_n}}\Vert e^{-\mu_n \zeta e^{-i{a+b\over 2}}}T(\zeta +2^{-n}e^{i{a+b\over 2}})\Vert_{op, \mu_0} \le e^{2^{-n}(\mu_0+1)} \ (n\ge 1).\end{equation}

where $\Vert .\Vert_{op,\mu_0}$ is the norm on the normalization $\ma_{T_0}$ of $\mathcal A$ with respect to the semigroup $ T_0:=(T(te^{i{(\alpha +\beta)\over 2}}))_{t>0}$
associated to $\mu_0$ in theorem 2.2.

For $n\ge 1$ set $V_{n}:=\{ e^{-\mu_n \zeta e^{-i{(\alpha_n+\beta_n)\over 2}}}T(\zeta)\}_{\zeta \in \overline{\mathcal S}_{\alpha_n, \beta_n}},$ set $W_{ n}:=V_{1}\dots V_{n}$ and set $W=\cup_{n\ge 1}W_{ n},$ so that $W$ is stable under products. Set ${\mathcal L}_{\mub}:=\{ u \in \ma_{T_0} \ | \ \sup_{w\in W}\Vert wu\Vert_{\mub_0,op} <+\infty,$ and set $\Vert u\Vert_{\mub}= \sup_{w\in W\cup\{I\}}\Vert wu\Vert_{\mub_0,op}$ for $u\in {\mathcal L}_{\mub}. $ Denote by by ${\mathcal A}_{\mub, T}$ the closure of ${\mathcal A}_{T_0}$ in ${\mathcal M}({\mathcal L}_{\mub}).$Then $T(\zeta)u\in {\mathcal L}_{\mub}$  for $\zeta \in S_{a,b}, u \in \ma_{T_0},$ $({\mathcal L}_{\mub}, \Vert . \Vert _{\mub})$ is a Banach algebra, the inclusion map $j: \ma \to \ma_{\mu, T}$ is a $\qm$ homomorphism with respect to $T(\zeta)w^2$ for $w\in \Omega(\ma),$ and we have the following properties

(i) The tautological map $\tilde j: S_{u/v}\to S_{u/v}$ is a pseudobounded isomorphism from $\qm(\ma)$ onto $\qm(\ma_{\mu,T})$ and if $ \zeta \in S_{a,b}$ then $\tilde j^{-1}(S)=S_{T(\zeta )w^2u/T(\zeta) w^2v}$ for $S=S_{u/v}\in \qm(\ma_{\mu, T}),$ $w\in \Omega(\ma).$

(ii)  $\tilde j(\mathcal M(\ma))\subset \mathcal M(\ma_{\mu,T})$  and $\Vert R\Vert_{\mathcal M(\ma_{\mu, T})}\le \Vert R \Vert_{\mathcal M(\ma)}$ for $R\in \mathcal M(\ma).$

(iii)$\Vert T(\zeta) \Vert_{{\mathcal M}({\ma}_{\mu,T})} \le \Vert T(\zeta) \Vert_{{\mathcal M}({\mathcal L}_{\mub})}  \le e^{\mu_nRe(\zeta e^{-i{a+b\over 2}})}$ for $\zeta \in \overline{\mathcal S}_{\alpha_n,\beta_n}, n\ge 1,$

(iv)  If $a <\alpha  < \beta <b,$ then we have, for  $v\in {\mathcal A}_{\mub, T},$
   $$\lim \sup_{\stackrel{\zeta \to 0}{_{\zeta \in S_{\alpha,\beta}}}}\Vert T(\zeta)v-v\Vert _{\mathcal M(\ma_{\mu, T})}\le \lim \sup_{\stackrel{\zeta \to 0}{_{\zeta \in S_{\alpha,\beta}}}}\Vert T(\zeta)v-v\Vert_{{\mathcal M}({\mathcal L}_{\mu})}=0.$$

\end{prop}

Proof: Since $\Vert T_0(2^{-n})\Vert_{op, \mu_0}\le e^{2^{-n}\mu_0},$ the existence of a sequence $(\mu_n)_{n\ge 1}$ satisfying the required conditions follows from the lemma. 

Let $t>0,$ let $n_0\ge 2$ be such that $2^{-n_0+1}<t,$ and let $v\in W.$ Since $1\in V_n$ for $n\ge 1,$ we can assume that $v=v_1\dots v_n,$ where $n \ge n_0$ and where $v_j \in V_j$ for $1\le j \le n.$ Then $\Vert T_0(2^{-j})v_j \Vert_{op, \mu_0}\le  e^{2^{-j}(\mu_0+1)}$ for $j\ge n_0,$ and so $\Vert T_0(2^{-n_0+1}-2^{-n})v_{n_0}\dots v_n\Vert_{op,\mu_0}\le e^{(2^{-n_0}-2^{-n})(\mu_0+1)}.$ We obtain

$$\Vert T_0(2^{-n_0+1})v_{n_0}\dots v_n\Vert_{op,\mu_0} $$ $$\le \Vert T_0(2^{-n})\Vert_{op, \mu_0} \Vert T((2^{-n_0+1}-2^{-n}))v_{n_0}\dots v_n\Vert_{op,\mu_0}$$ $$\le e^{2^{-n_0+1}(\mu_0+1)-2^{-n}}\le e^{2^{-n_0+1}(\mu_0+1)}.$$

Set $r={t-2^{-n_0+1}\over n_0}.$ It follows from lemma 4.1(i) that for every $j\le n_0-1$ there exists $k_j>0$ such that $sup_{v\in V_j}\Vert T(r)v\Vert_{op, \mu_0}\le k_j.$ This gives

$$\Vert T_0(t)v\Vert_{op,\mu_0}\le \Vert T_0(r)v_1\Vert_{op,\mu_0}\dots\Vert T_0(r)v_{n_0-1}\Vert_{op,\mu_0} \Vert T_0(2^{-n_0+1})v_{n_0}\dots v_n\Vert_{op,\mu_0}$$
$$\le k_1\dots k_{n_0-1}e^{2^{-n_0+1}(\mu_0+1)},$$

and so $T_0(t)\ma_{T_0} \subset {\mathcal L}_{\mu}.$ Now let $\zeta \in S_{a,b}.$ Then $\zeta-te^{i{a+b\over 2}}\in S_{a,b}$ for some $t>0,$ and so $T(\zeta)\ma_{T_0}=T(\zeta- te^{i{(\alpha+\beta)\over 2}})T_0(t)\ma_{T_0}\subset T(\zeta- te^{i{(\alpha+\beta)\over 2}}){\mathcal L}_{\mu}\subset {\mathcal L}_{\mu}.$

Since  $W$ is stable under products, the fact that $\mathcal L_{\mu}$ is a Banach algebra  follows from proposition 2.3, which also implies (ii) and (iii). Let $\zeta_0\in S_{a,b}.$ Since (iv) holds for $u\in T(\zeta_0)\ma_{\mu,T},$ and since $T(\zeta_0)\ma_{\mu,T})$ is dense in $(\ma_{\mu, T}, \Vert .\Vert_{{\mathcal L}_{\mu}})$ (iv) follows from (iii) by a standard density argument.

Let $\zeta \in S_{a,b},$ and let $w \in \Omega(\ma) \subset \Omega(\ma_{T_0}).$ Since $T(\zeta/2)w\in \Omega(\ma),$ and since   $\lim \sup_{t\to 0^+}\Vert T_0(t)T(\zeta)w\Vert <+\infty,$ it follows from proposition 2.3 that the inclusion map $j_0: \ma \to \ma_{T_0}$ is a $\qm$-homomorphism with respect to $T(\zeta/2)w.$ Since $T(\zeta/2)w\subset {\mathcal M}_{\mu}\cap \Omega(\ma_{T_0}),$ it follows also from proposition 2.3 that the inclusion map $j_1:\ma_{T_0}\to \ma_{\mu, T}$ is a $\qm$-homomorphism with respect to $T(\zeta/2)w,$ 
and it follows then from proposition 2.4 that the inclusion map $j=j_1\circ j_0: \ma \to \ma_{\mu, T}$ is a $\qm$ homomorphism with respect to $T(\zeta/2)wT(\zeta/2)w=T(\zeta)w^2.$
It follows then from proposition 2.2 that the tautological map $\tilde j: S_{u/v}\to S_{u/v}$ is a pseudobounded homomorphism from $\qm(\ma)$ onto $\qm(\ma_{\mu, T})$ and that $\tilde j^{-1}(S_{u/v})=S_{T(\zeta)uw^2/T(\zeta)vw^2}$ for $S=S_{u/v}\in \qm(\ma_{\mu, T}).$ Since $T(\zeta)\ma_{\mu,T}\subset T(\zeta/2)\left [T(\zeta/2)\ma_{\mu,T}\right ]\subset T(\zeta/2)\ma_{T_0}\subset \ma,$ and since $T(\zeta)v\ma$ contains $T(\zeta)vw^2\ma$ which is dense in $\ma$ for $v\in \Omega(\ma_{\mu,T}),$ we have $T(\zeta)v\in \Omega(\ma)$ for $\zeta \in S_{a,b}, v \in \Omega(\ma_{\mu, T}),$ and so $\tilde j^{-1}(S)=S_{T(\zeta)uw^2/T(\zeta)vw^2}$ for $S=S_{u/v} \in \qm(\ma_{\mu, T}).$ $\square$

We will use the following notion.

\begin{defn} Let $\ma$ be a weakly commutative Banach algebra with dense principal ideals, and let $T:=(T(\zeta))_{\zeta\in S_{a,b}}$ be an analytic semigroup of multipliers on $\ma$ such that  $T(\zeta)\ma$ is dense in $\ma$ for $\zeta \in S_{a,b}.$ A normalization of the algebra $\ma$ with respect to the semigroup $T$ is a subalgebra $\mb$ of $\qm(\ma)$ which is a Banach algebra with respect to a norm $\Vert .\Vert_{\mb}$ and satisfies the following conditions

(i) There exists $u \in \Omega(\ma)$ such that the inclusion map $j: \ma \to \mb$ is a $\qm$-homomorphism with respect to $T(\zeta)u$ for every $\zeta \in S_{a,b}$, and $\Vert j(u)\Vert_{\mb}\le \Vert u \Vert_{\ma}$ for $u \in \ma.$

(ii) $\tilde j (R) \subset \mathcal M(\mb),$ and $\Vert \tilde j(R)\Vert _{\mathcal M(\mb)}\le  \Vert R\Vert _{\mathcal M(\ma)}$ for $R\in \mathcal M(\ma),$ where $\tilde j: \qm(\ma)\to \qm(\mb)$ is the pseudobounded isomorphism associated to $j$ introduced in proposition 2.2 (ii).

(iii) $\lim \sup_{\stackrel{t\to 0}{_{\zeta\in S_{\alpha,\beta}}}}\Vert \tilde j(T(t))\Vert_{\mathcal M(\mb)} <+\infty$ for $a<\alpha <\beta < b.$

\end{defn}

If $\mb$ is a normalization of $\ma$ with respect to the holomorphic semigroup $T=(T(\zeta))_{\zeta \in S_{a,b}},$ a standard density argument shows that  if $a<\alpha <\beta < b$ then $\lim \sup_{\stackrel{\zeta\to 0}{_{\zeta in S_{\alpha, \beta}}}}\Vert T(\zeta)u-u\Vert_{\mb}=0$ for every $u\in \mb$

Notice that the algebra ${\mathcal A}_{0, T}={\mathcal A}_{\mu_0,\bf T}$ and its norm topology associated to the norm $\Vert .\Vert_{\mu_0, op}$ discussed above do not depend on the choice of $\mu_0.$ This is no longer the case  for the Banach algebra ${\mathcal A}_{\mu, T}$ and its norm topology, which may depend on the choice of the sequence $\mu.$ In order to get a more intrinsic renormalization one could consider the Fr\'echet algebra ${\mathcal L}:=\cap_{n\ge 1}\{ u \in {\mathcal A_ T \ | \ \sup_{\stackrel{\zeta \in {\overline S}_{\alpha_n,\beta_n}}{_{\Vert u \Vert_{\mu_0,op}\le 1}}}}\Vert T(\zeta )u\Vert _{\mu_0,op}<+\infty\},$ then consider the closed subalgebra ${\mathcal U}$ of $\mathcal L$ generated by the semigroup and introduce an intrinsic normalization of $\at$ to be the closure of ${\mathcal U}$ in ${\mathcal M}(\mathcal U)$ with respect to the Mackey-convergence associated to a suitable notion of boundedness on subsets of ${\mathcal M}(\mathcal U),$ but it seems more convenient to adopt the point of view used in proposition 4.2.

\section{Generator of a strongly continuous semigroup of multipliers and Arveson spectrum}

In this section we consider again  a weakly cancellative commutative Banach algebra ${\ma}$ with dense principal ideals and a strongly continuous semigroup $T=(T(t))_{t>0}$ of multipliers on $\ma$ such that $\cup_{t>0}T(t)\ma$ is dense in $\ma.$ We set again $\omega_{T}(t)=\Vert T(t)\Vert$ for $t>0.$ Denote by ${\mathcal M}_{\omega_T}$ the space of all measures $\nu$ on $(0,+\infty)$ such that $\int_0^{+\infty}\omega_T(t)d\vert \nu\vert (t)<+\infty,$ and for $\nu \in {\mathcal M}_{\omega_T}$ define a $\Phi_{T}: {\mathcal M}_{\omega_T}\to \mathcal M(\ma)$ by the formula 

$$\Phi_T(\nu)u:=\int_0^{+\infty}T(t)ud \nu(t) \ \ (u\in \ma). \ \ \ \ \ \ \ \  \ \ \ \ \ \ \ \ \ \ \ \ \ \ \ \  \ \ \ \ \ \ \ \ \ \ \ \ \ \ \ \  \ \ \ \ \ \ \ \ \ \ \ \ \ \ \ \  \ \ \ \ \ \ \ \ (2^{'})$$

The Bochner integral is well-defined since the semigroup is strongly continuous, ${\mathcal M}_{\omega_T}$ is a Banach algebra with respect to convolution of measures on the half-line, and we will identify again the space $L^1_{\omega_T}$ of (classes of) measurable functions on $[0,+\infty)$ satisfying $\int_0^{+\infty}\vert f(t)\vert \Vert T(t)\Vert dt <+\infty$ to the ideal of all $\nu \in {\mathcal M}_{\omega_T}$ which are absolutely continuous with respect to Lebesgue measure. Denote by ${\mathcal B}_T$ the closure of $\phi_T(\mathcal M_{\omega})$ in ${\mathcal M}(\ma),$ and denote again by $\iit$ the closure of $\Phi_T(L^1_{\omega_T})$ in ${\mathcal M}(\ma),$ so that the "Arveson ideal" $\iit$ is a closed ideal of ${\mathcal B}_{\bf T}.$

The idea of considering the generator of a semigroup as a quasimultiplier on some suitable Banach algebra goes back to \cite{ga} and \cite{gm} for groups of bounded operators and, more generally, for groups of regular quasimultipliers. An obvious such interpretation was given by I. Chalendar, J. R. Partington and the author in \cite{cep0} for analytic semigroups, and the author interpreted in section 8 of \cite{e4} the generator of a semigroup of bounded operators which is weakly continuous in the sense of Arveson \cite{ar} as a quasimultiplier on the corresponding Arveson ideal $\iit.$ Since in the present context $\Phi_T\left (\Omega(L^1_{\omega_T})\right )\subset \Omega(\iit)$ and $uv \in \Omega(\ma)$ for $u\in \Omega(\iit), v \in \Omega(\ma),$ the map $j_{T}: S_{u/v}\to S_{uw/vw}$ is a pseudobounded homomorphism from $\qm(\iit)$ into $\qm(\ma)$ for every $w\in \Omega(\ma),$ and the definition of $j_T$ does not depend on the choice of $w.$

The generator $\Delta_{T,\iit}$ of $T$ considered as a strongly continuous semigroup of multipliers on $\iit$ has been defined in \cite{e4}, def. 8.1 by the formula

\begin{equation}\Delta_{T, \iit}=S_ {-{\Phi_{T}(f_0')/\Phi_{T}(f_0)}},\end{equation}

where $f_0 \in {\mathcal C}^1([0,+\infty))\cap \Omega \left (L^1_{\omega_{T}}\right )$ satisfies $f_0=0,$ $f_0' \in L^1_{\omega_{T}},$

and an easy verification given in \cite{e4} shows that this definition does not depend on the choice of $f_0.$ This suggests the following definition

   \begin{defn}:  The infinitesimal generator $\Delta_{T, \ma}$ of  ${T}$ is the quasimultiplier on $\ma$ defined by the formula
 
$$\Delta_{T, \ma}=j_T(\Delta_{T,\iit})=S_ {-{\Phi_{T}(f_0')u_0/\Phi_{T}(f_0)u_0}},$$

where $f_0 \in {\mathcal C}^1([0,+\infty))\cap \Omega \left (L^1_{\omega_{T}}\right )$ satisfies $f_0=0,$ $f_0' \in L^1_{\omega_{T}},$ and where $u_0\in \Omega(\ma).$

\end{defn}
   
Assume that $f_1$ and $u_1$ also satisfy the conditions of the definition, and set $f_2=f_0*f_1, u_2=u_0u_1$ Since $\Omega \left (L^1_{\omega_{T}}\right )$ is stable under convolution, $f_2 \in  \Omega \left (L^1_{\omega_{T}}\right ),$ and $f'_2=f'_0*f_1=f_0*f'_1\in L^1_{\omega_{T}}$ is continuous. Also $f_2(0)=0,$ and we have

$$\Phi_{ T}(f'_2)u_2\Phi_{T}(f_0)u_0 =\Phi_{T}(f'_0)\Phi_{T}(f_1)u_2\Phi_T(f_0)u_0=\Phi_T(f'_0)u_0\Phi_T(f_2)u_2$$

and similarly $\Phi_{ T}(f'_2)u_2\Phi_{T}(f_1)u_1 =  \Phi_T(f'_1)u_1\Phi_T(f_2)u_2,$ which shows that the definition of $\Delta_{\bf T, \ma}$ does not depend on the choice of $f_0$ and $u_0.$

\begin{prop} Let $T_1$ and $T_2$ be two semigroups satisfying the conditions of definition 5.1with respect to $\ma.$ If $T_1(t)T_2(t)=T_2(t)T_1(t)$ for $t>0,$ then $\Delta_{T_1T_2,\ma}= \Delta_{T_1,\ma}+\Delta_{T_2,\ma}.$

\end{prop}

Proof: Let $f_1$ and $f_2$ be two functions on $[0,+\infty)$ satisfying the conditions of definition 5.1 with respect to $T_1$ and $T_2.$ Since $\omega_{T_1T2}(t)\le \omega_{T_1}(t)\omega_{T_1}(t)$ for $t>0,$ $f_1f_2$ satisfies the conditions of definition 5.1 with respect to $T_1T_2,$ and it follows from Leibnitz rule and Fubini's theorem that we have

$$\Phi_{T_1T_2}(f_1f_2)= \Phi_{T_1}(f_1)\Phi_{T_2}(f_2), \Phi_{T_1T_2}((f_1f_2)')=\Phi_{T_1}(f'_1)\Phi_{T_2}(f_2)+\Phi_{T_1}(f_1)\Phi_{T_2}(f'_2),$$

and the results follows. $\square$
 
 We now give a link between the quasimultiplier approach and the classical approach based on the study of ${T(t)u-u\over t}$ as $t \to 0^+.$ A proof of the following folklore result  
  is given for example in \cite{e4}, lemma 8.4.
  
  \begin{lem} Let $\omega$ be a lower semicontinuous submultiplicative weight on $(0,+\infty),$ and let $f \in {\mathcal C}^1([0,+\infty))\cap L^1_{\omega}.$ If 
  $f(0)=0,$ and if $f' \in L^1_{\omega},$ then the Bochner integral $\int_t^{+\infty}(f'*\delta_s)ds$ is well-defined in $L^1_{\omega}$ for $t\ge 0,$ and we have
  
  $$f*\delta_t  -f =-\int_0^{t}(f'*\delta_s)ds, \ \ \mbox{and} \lim_{t\to 0^+}\left \Vert \frac{f*\delta_t-f}{t}+f'\right \Vert _{L^1_{\omega}}=0.$$
  
  \end{lem}
  
It follows from the lemma that we have if $f\in {\mathcal C}^1([0,+\infty))\cap L^1_{\omega_T},$ and if $f'\in  L^1_{\omega_T},$

\begin{equation} T(t)\Phi_T(f)-\Phi_T(f)=-\int_0^tT(s)\Phi_T(f')ds  \ \ (t\ge 0). \end{equation}
  
  \begin{prop} (i) Let $u \in \ma.$ If $\lim_{t\to 0^+}\Vert \frac{T(t)u-u}{t}-v\Vert =0$ for some $v \in \ma,$ then  $u \in {\mathcal D}_{\Delta_{T, \ma}},$ and $\Delta_{T, \ma}u=v.$
  
  (ii) Conversely if $\lim \sup_{t\to 0^+}\Vert T(t) \Vert <+\infty,$ then $\lim \sup_{t\to 0^+} \Vert {T(t) u -u \over t} -\Delta_{T,\ma}u\Vert =0$ for every $u \in {\mathcal D}_{\Delta_{T, \ma}}.$
  
  \end{prop}
  
  Proof: (i) If $u \in \ma,$ and if $\lim_{t\to 0^+}\Vert \frac{T(t)u-u}{t}-v\Vert =0$ for some $v \in \dt,$ let $f_0 \in {\mathcal C}^1([0,+\infty))\cap \Omega \left (L^1_{\omega_{\bf T}}\right )$ satisfiying $f_0=0,$ $f_0' \in L^1_{\omega_{T}},$ and let $u_0 \in \Omega(\ma).$ It follows from the lemma that we have, with respect to the norm topology on $\ma,$
  
  $$-\Phi_{\bf T}(f'_0)u_0u=\left [\lim_{t\to 0^+}\frac {T(t)\Phi_{\bf T}(f_0)-\Phi_{\bf T}(f_0)}{t}\right ]u_0u= \Phi_{\bf T}(f_0)u_0\left [\lim _{t\to 0^+} {T(t)u-u\over t}\right ]$$ $$= \Phi_{\bf T}(f_0)u_0v,$$
  
  and so $u \in  {\mathcal D}_{\Delta_{ T, \ma}},$ and $\Delta_{T, \ma}u=v.$
  
 (ii)  Conversely assume that $\lim \sup_{t\to 0^+}\Vert T(t)\Vert <+\infty,$ let $u\in {\mathcal D}_{\Delta_{T, \ma}},$ let $f_0 \in {\mathcal C}^1([0,+\infty))\cap \Omega \left (L^1_{ \omega_{\bf T}}\right )$ such that $f'_0\in L^1_{ \omega_{\bf T}},$ and let $u_0\in \Omega(\ma).$ Let  $v=\Delta_{{\bf T},\jjt}u.$  It follows from (6) that we have, for $t\ge 0,$
  
  $$ \Phi _{\bf T}(f_0)u_0\int_0^tT(s)vds= \int_0^tT(s) \Phi_{\bf T}(f_0)u_0vds=- \int_0^tT(s) \phi_{\bf T}(f'_0)u_0uds$$ $$=-\left [ \int_0^tT(s)\Phi_{\bf T}(f'_0)ds\right ]u_0u=\left [T(t) \Phi_{\bf T}(f_0)u_0- \Phi_{\bf T}(f_0)u_0\right ]u=\Phi_{\bf T}(f_0)u_0(T(t)u-u).$$
  
  Since $ \Phi _{\bf T}(f_0)u_0\in \Omega(\ma),$ this shows that $T(t)u-u=\int_0^tT(s)vds,$ and so $\lim_{t\to 0^+}\left \Vert \frac{T(t)u-u}{t}-v\right \Vert=0.$ $\square$

  We now consider a normalization $\ma$ with respect to $T$, see definition 3.3.
  
\begin{prop} Let $\mb$ be a normalization of $\ma$ with respect to $T.$ Set $v_{\lambda}(t)=te^{-\lambda t}$ for $\lambda \in \R, t\ge 0,$ and let $\lambda > log(\rho_T).$

(i) If $u \in \Omega(\mb),$ then 

$$\Delta_{T,\mb}=-S_{\int_0^{+\infty}v'_{\lambda}(t)T(t)udt/\int_{0}^{+\infty}v_\lambda(t)T(t)udt} .$$

(ii) Let $\tilde j: \qm(\ma)\to \qm(\mb)$ be the pseudobounded isomorphism given in proposition 3.2 (i). Then $\tilde j^{-1}(\Delta_{T,\mb}) =\Delta_{T,\ma}.$ So if $u \in \Omega(\ma)$, and if $\lim \sup_{t\to 0^+}\Vert T(t)u\Vert
 <+\infty,$ then
 
 $$\Delta_{T, \ma}=-S_{\int_0^{+\infty}v'_{\lambda}(t)T(t)udt/\int_{0}^{+\infty}v_\lambda(t)T(t)udt} .$$
  
  \end{prop}
  
  Proof: (i) Set $e_\lambda(t)=e^{\lambda t}$ for $\lambda \in \R, t\ge 0.$ If $ \lambda >\mu > log(\rho_{\bf T}),$ then $v_{\lambda} \in L^1_{e_{\mu}}\subset L^1_{\tilde \omega_T},$ where $\tilde \omega_T(t)=\Vert T(t)\Vert_{{\mathcal M}(\ma_T)},$ and $L^1_{e_{\mu}}$ is dense in  $L^1_{\tilde \omega_T}$ since it contains the characteristic function of $[\alpha, \beta]$ for $0<\alpha < \beta <+\infty.$ It follows for example from Nyman's theorem \cite{ny} about closed ideals of $L^1(\R^+)$ that $v_{\lambda}e_{\mu}\in \Omega(L^1(\R^+))$ and so $v_{\lambda}\in \Omega(L^1_{e_{\mu}})\subset \Omega(L^1_{\tilde \omega_T}).$ So $v_{\lambda}$ and $u$ satisfy the conditions of definition 5.1 with respect to $T$ and $\mb,$
  and (i) holds.
  
  (ii) The map $\tilde j$ is the tautological map $S_{u/v}\to S_{u/v},$ where $u \in \ma \subset \mb $ and $v \in \Omega(\ma)\subset \Omega(\mb).$ Now let $f_0\in L^1_{\omega_T}\cap{\mathcal C}^1([0,+\infty))$ satisfying definition 5.1 with respect to $T$ and $\ma$ and let $u_0\in \Omega(\ma).$ Since $\Omega(L^1_{\omega_T})\subset \Omega(L^1_{\tilde \omega_T}),$ and since $\Omega(\ma)\subset \Omega(\mb),$ it follows from definition 5.1 that $\tilde j (\Delta_{T, \ma})=\Delta_{T,\mb},$ and so $\tilde j^{-1}(\Delta_{T,\mb}) =\Delta_{T,\ma}.$

  Let $u \in \Omega(\ma)\subset \Omega(\mb)$, and assume that $\lim \sup_{t\to 0^+}\Vert T(t)u\Vert
 <+\infty.$ Let $w\in \Omega(\ma)$ be such that $w\mb \subset \ma.$ Since $v_{\lambda}\in \Omega(L^1_{\omega_{\lambda}}),$ we see as in the proof of lemma 3.1 that  $\int_0^{+\infty}v_{\lambda}(t)T(t)udt\in \Omega(\mb),$ and it follows from proposition 2.2 that $w\int_{0}^{+\infty}v_{\lambda}(t)T(t)udt \in \Omega(\ma).$ Using the characterization  of $\tilde j^{-1}$ given in proposition 2.2, we obtain
 
  $$\Delta_{T,\mb}=-S_{\int_0^{+\infty}v'_{\lambda}(t)T(t)udt/\int_{0}^{+\infty}v_\lambda(t)T(t)udt},$$ and 
 $$\Delta_{T,\ma}=\tilde j^{-1}(\Delta_{T,\mb})= -S_{w\int_0^{+\infty}v'_{\lambda}(t)T(t)udt/w\int_{0}^{+\infty}v_{\lambda}(t)T(t)udt}$$ $$=-S_{\int_0^{+\infty}v'_{\lambda}(t)T(t)udt/\int_{0}^{+\infty}v_\lambda(t)T(t)udt }.$$ $\square$

 We will  denote by $\sia$ the space of characters of $\iit,$ equipped with the usual Gelfand topology. Notice that if $\chi \in \sia$ then there exists a unique character
$\tilde \chi$ on ${\mathcal {QM}}(\iit)$ such that $\tilde \chi _{|_{\iit}}=\chi,$ which is defined by the formula $\tilde \chi(S_{u/v})=\frac{\chi(u)}{\chi(v)}$ for $u\in \iit, v \in \Omega(\iit).$ 

\begin{defn} Assume that $\iit$ is not radical, and let $S \in \qm(\iit).$ The Arveson spectrum $\sigma_{ar}(S)$ is defined by the formula

$$\sigma_{ar}(S)=\{\lambda = \tilde \chi(S) \ : \ \chi \in \sia\}.$$

\end{defn}

If $\nu$ is a measure on $[0,+\infty),$ the Laplace tranform of $\nu$  is defined by the usual formula ${\mathcal L}(\nu)(z)=\int_0^{+\infty}e^{-zt}d\nu(t)$
when $\int_0^{+\infty}e^{-Re(z)t}d\vert \nu \vert(t)<+\infty.$

We have the following easy observation.

\begin{prop}  Let $\nu \in \M.$ Then we have, for $\chi \in \sia,$ 

\begin{equation} \tilde \chi \left ( \int_{0}^{+\infty}T(t)d\nu(t)\right)={\mathcal L}(\nu)(-\tilde \chi(\Delta_{{T},\iit})).\end{equation}

Similarly we have, for $\nu \in {\mathcal M}_{\tilde \omega_{T}}, \chi \in \sia,$

\begin{equation} \tilde \chi(\Phi_T(\nu))=\tilde \chi \left ( \int_{0}^{+\infty}T(t)d\nu(t)\right)={\mathcal L}(\nu)(-\tilde \chi(\Delta_{{T},\iit}))={\mathcal L}(\nu)(-\tilde \chi(\Delta_{{\bf T},\jjt})) .\end{equation}

In particular $\tilde \chi (T(t))=e^{\tilde \chi(\Delta_{{T},\iit})t}$ for $t>0.$

\end{prop}

Proof: If $\chi \in \sia,$ then $\tilde \chi_{|_{\at}}$ is a character on $\at,$ the map $t \to \tilde \chi (T(t))$ is continuous on $(0,+\infty)$ and so there exists $\lambda \in \C$ such that $\tilde \chi (T(t))=e^{-\lambda t}$ for $t>0,$ and $\vert e^{-\lambda t}\vert\le \Vert T(t)\Vert,$ which shows that $Re(\lambda)\ge -log(\rho_{\bf T}).$

Let $u\in \Omega(\iit),$ and let $\nu \in \M.$ We have

$$\chi (u)\tilde \chi \left (\int_0^{+\infty}T(t)d\nu(t)\right )= \chi \left (u\int_0^{+\infty}T(t)d\nu(t)\right )
=\chi \left (\int_0^{+\infty}T(t)ud\nu(t)\right )$$
$$=\int_0^{+\infty}\chi(T(t)u)d\nu(t)=\chi(u)\int_0^{+\infty}e^{-\lambda t}d\nu(t)=\chi(u){\mathcal L}(\nu)(\lambda),$$

and so $\tilde \chi(\Phi_T(\nu))={\mathcal L}(\nu)(\lambda).$

Let $f_0\in {\mathcal C}^1((0,+\infty))\cap \Omega(\iit)$ such that $f_0(0)=0.$ We have

$$\lambda {\mathcal L}(f_0)(\lambda)={\mathcal L}(f'_0)(\lambda)=\chi (\Phi_{\bf T}(f'_0))=-\tilde \chi\left (\Delta_{{ T},\iit}\Phi_{\bf T}(f_0)\right )$$ $$=-\tilde \chi(\Delta_{{ T},\iit})\chi (\phi_{\bf T}(f_0))= -\tilde \chi(\Delta_{{ T},\iit}){\mathcal L}(f_0)(\lambda),$$

and so $\lambda =-\tilde \chi(\Delta_{{T},\iit}),$ which proves (7), and  formula (8) follows from a similar argument. In particular $\chi(T(t))={\mathcal L}(\Delta_t)(-\tilde \chi(\Delta_{{\bf T},\iit}))= e^{\tilde \chi(\Delta_{{T},\iit})t}$ for $t>0.$ $\square$

  The following consequence of  proposition 5.7 pertains to folklore.

\begin{cor} Assume that  $\iit$ is not radical. Then the map $\chi \to \tilde \chi(\Delta_{{ T},\iit})$ is a homeomorphism from $\widehat{\iit}$ onto $\sigma_{ar}(\Delta_{{T},\iit}),$ and the set
$\Lambda_{t}:=\{ \lambda \in \sigma_{ar}(\Delta_{{T},\iit}) \ | \ Re(\lambda)\le t\}$ is compact for every $t \in \R,$ so that $\sigma_{ar}(\Delta_{{T},\iit})$ is closed.
\end{cor}

Proof: Let $f_0\in {\mathcal C}^1((0,+\infty))\cap \Omega(\iit)$ such that $f_0(0)=0.$ We have $\chi(\Phi_{{\bf T}}(f_0))\neq 0$ and $\tilde \chi(\Delta_{{ T},\iit})=-\frac{\chi(\Phi_{\bf T}(f'_0))}{\chi(\Phi_{\bf T}(f_0))}$ for $\chi \in \widehat{\iit},$ 
and so  the map $\chi \to \tilde \chi(\Delta_{T})$ is continuous with respect to the Gelfand topology on $\widehat{\iit}.$

Conversely let $f \in \LL.$ It follows from proposition 5.6 that we have, for $\chi \in \widehat{\iit},$

$$\chi(\Phi_{\bf T}(f))={\mathcal L}(f)(-\tilde \chi(\Delta_{{T},\iit})).$$

Since the set $\{ u= \Phi_{\bf T}(f) : f \in \LL \}$ is dense in $\iit,$ this shows that the map $\chi \to \tilde {\chi}(\Delta_{{T},\iit})$ is one-to-one on  $\widehat{\iit},$ and that the inverse map $\sigma_{ar}(\Delta_{{T},\iit})\to \widehat{\iit}$ is continuous with respect to the Gelfand topology.

Now let $t \in \R,$ and set $U_t:=\{ \chi \in \widehat{\iit} :  Re(\chi(\Delta_{{\bf T},\iit}))\le t\}.$ Then $\vert \tilde {\chi}(T(1))\vert\ge e^{-t}$ for $\chi \in U_t,$ and so $0$ does not belong to the closure of $U_t$ with respect to the weak$^*$ topology on the unit ball of the dual of $\iit.$ Since $\widehat{\iit}\cup \{0\}$ is compact with respect to this topology, $U_t$ is a compact subset of $\widehat{\iit}$, and so the set $\Lambda_t$ is compact, which implies that $\sigma_{ar}(\Delta_{{T},\iit})=\cup_{n\ge 1}\Lambda _n$ is closed. $\square$

\section{The resolvent}
We now wish to discuss the resolvent of the generator of a strongly continuous semigroup ${ T}=(T(t))_{t>0}$ of multipliers on $\ma,$ where $\ma$  is a weakly cancellative commutative Banach algebra with dense principal ideals, and where $\cup_{t>0}T(t)\ma$ is dense in $\ma.$ From now on we will write $\Delta_T=\Delta_{T,\ma}$ and we will denote  by $\mathcal D_{\Delta_{T},\ma}$  the domain of $\Delta_{\bf T}$ considered as a quasimultiplier on $\ma$. The Arveson ideal $\iit$ is as above the closed subalgebra of ${\mathcal M}(\ma)$ generated by $\Phi_{\bf T}(L_{\omega_T}^1).$ 

The Arveson resolvent set is defined by the formula $Res_{ar}(\Delta_{T,\iit})=\C \setminus \sigma_{ar}(\Delta_{T,\iit}),$ with the convention $\sigma_{ar}(\Delta_{T,\iit})=\emptyset$ if $\iit$ is radical. The usual "resolvent formula," interpreted in terms of quasimultipliers, shows that $\lambda I -\Delta_{T,\iit}\in \qm(\iit)$ is invertible in $\qm(\iit)$ and that its inverse $(\lambda I -\Delta_{T,\iit})^{-1}$ belongs to the
Banach algebra $\jjt \subset \qm_r(\iit)$ obtained by applying theorem 2.2 to $\iit$ with respect to the semigroup $T,$ and that we have, for $Re(\lambda)> log(\rho_{\bf T}),$

$$(\lambda I -\Delta_{T,\iit})^{-1}=\int_0^{+\infty}e^{-\lambda s}T(s)ds \in \jjt,$$

where the Bochner integral is computed with respect to the strong operator topology on ${\mathcal M}(\jjt).$ Also the $\jjt$-valued map $\lambda \to (\lambda I -\Delta_{T,\iit})^{-1}$ is holomorphic on $Res_{ar}(T,\iit).$ The details of the adaptation to the context of quasimultipliers of this classical part of semigroup theory are given in \cite{e4}, proposition 10.2. 

We now give a slightly more general version of this result, which applies in particular to the case where $\mb$ is the normalization $\ma_T$ of $\ma$ with respect to the semigroup $T$ introduced in proposition 2.2.



In the following we will identify the algebras $\qm(\ma)$ and $\qm(\mb)$ using the isomorphism $\tilde j$ intoduced in proposition 2.2 (iii) if $\mb$ is a normalization of $\ma$ with respect to $T.$ We set $\Phi_{T,\mb}(\nu)u=\int_0^{+\infty}T(t)ud\nu(t)$ for $u\in \mb, \nu\in {\mathcal M}_{\omega_{T,\mb}},$ where $\omega_{T,\mb}(t)=\Vert T(t)\Vert_{\mathcal M(\mb)}$ for $t>0,$ and we denote by $\mathcal I_{T,\mb}$ the closure of $\Phi_{T, \mb}( L^1_{\omega_{T,\mb}})$ in $\mathcal M(\mb).$ 

\begin{prop} Let $\ma$ be a weakly cancellative commutative Banach algebra with dense principal ideals, let $T=((T(t))_{t>0}$ be a strongly continuous semigroup of multipliers on $\ma$ such that $T(t)\ma$ is dense in $\ma$ for $t>0,$ and let $\mb$ be a normalization of $\ma$ with respect to $T.$ Set $Res_{ar}(\Delta_T)=Res_{ar}(\Delta_{T,\mathcal I_T})=\C \setminus \sigma_{ar}(\Delta_{T,{\mathcal I}_T}).$

The quasimultiplier $\lambda I-\Delta_{ T} \in {\mathcal QM}(\ma)$ admits an inverse $(\lambda I-\Delta_{T})^{-1}\in {\mathcal I}_{T,\mb}\subset {\mathcal M}(\mb) \subset {\mathcal QM}_r(\ma)$ for $\lambda \in Res_{ar}(\Delta_{T}),$ and the map $\lambda \to (\lambda I -\Delta_{T})^{-1}$ is an holomorphic map from $Res_{ar}(\Delta_{T})$ into $\mathcal I_{T,\mb}.$ Moreover we have, for $Re(\lambda) > log(\rho_{\bf T}),$

$$(\lambda I -\Delta_{ T})^{-1}=\int_0^{+\infty}e^{-\lambda s}T(s)ds \in \mathcal I_{T,\mb},$$

where the Bochner integral is computed with respect to the strong operator topology on ${\mathcal M}(\mb),$ and $\Vert (\lambda I -\Delta_{T})^{-1}\Vert_{{\mathcal M}(\mb)}\le \int_0^{+\infty}e^{-Re(\lambda)t}\Vert T(t)\Vert _{{\mathcal M}(\mb)}dt.$

\end{prop}

Proof: We could deduce this version of the resolvent formula from proposition 10.2 of \cite{e4}, but we give a proof for the sake of completeness. Set again $e_\lambda(t)=e^{\lambda t}$ for $t\ge 0, \lambda \in \C.$ Assume that $Re(\lambda)>log(\rho_{\bf T})\ge \lim_{t\to +\infty}{log\Vert T(t)\Vert_{\mathcal M (\mb)}\over t},$ let $v \in \mathcal I_{T,\mb},$ and set $a=\Phi_{T,\mb}(e_{-\lambda}).$  We have

$$av=\int_0^{+\infty}e^{-\lambda s}T(s)vds, T(t)av-av=\int_0^{+\infty}e^{-\lambda s}T(s+t)vds-\int_0^{+\infty}e^{-\lambda s}T(s)vds$$ $$=e^{\lambda t}\int_t^{+\infty}e^{-\lambda s}T(s)vds -\int_0^{+\infty}e^{-\lambda s}T(s)vds= (e^{\lambda t}-1)av -e^{\lambda t}\int_0^te^{-\lambda s}T(s)vds.$$

Since $\lim_{t\to 0^+}\Vert T(t)v-v\Vert_{\jjt}=0,$ we obtain

$$\lim_{t\to 0^+}\left \Vert \frac{T(t)av-av}{t} -\lambda av +v\right \Vert_{\jjt}=0,$$

and so $av \in \mathcal D_{\Delta_{\bf T},\jjt},$ and $\Delta_{T,\jjt}(av)=\lambda av -v.$ This shows that $a\jjt \subset \mathcal D_{\Delta_{ T,\jjt}},$ and that $(\lambda I -\Delta_{ T,\jjt})av=v$ for every
$v \in \jjt.$ We have $\lambda  I -\Delta_{ T,\jjt} =S_{u/v},$ where $u \in \jjt, v\in \Omega(\jjt),$ and we see that $ua=v.$ Hence $u\in \Omega(\jjt), \lambda  I -\Delta_{T,\jjt}$ is invertible in $\qm(\jjt),$ and $(\lambda  I -\Delta_{T,\mb})^{-1}=a= \Phi_{T,\mb}(e_{-\lambda})=\int_0^{+\infty}e^{-\lambda t}T(t)dt \in \jjt,$ where the Bochner integral is computed with respect to the strong operator topology on ${\mathcal M}(\jjt).$

Let $\chi \in \widehat \jjt.$ Then $\chi\circ \tilde j \in \widehat \iit,$ and so $\sigma_{ar}(\Delta_{T,\jjt})\subset \sigma_{ar}(\Delta_{T,\iit}).$ It follows then from proposition 2.6 that $\lambda I -\Delta_{T,\mb}$ has an inverse $(\lambda I -\Delta_{T,\mb})^{-1}\in \jjt$ in $\qm_r(\jjt)$ for $\lambda \in Res_{ar}(\Delta_{T,\iit})$ and that the $\jjt$-valued map $\lambda \to (\lambda I-\Delta_{T,\mb})^{-1}$ is holomorphic on $Res_{ar}(\Delta_{T,\mb}).$

Fix $u_0 \in \Omega(\ma)\subset \Omega(\mb),$ and  set $j_T(S)=S_{uu_0/vu_0}$  for $S=S_{u/v}\in \qm(\jjt).$ Then $j_T:\qm(\jjt)\to \qm(\mb)$ is a pseudobounded homomorphism,
and $j_T(\Delta_{T,\jjt})=\Delta_{T,\mb}.$ Identifying $\jjt$ to a subset of $\qm(\jjt)$ as above in the obvious way, we see that the restriction of $j_T$ to $\jjt$ is the identity map, and so $\lambda I-\Delta_{T}$ is invertible in $\qm(\mb)$ for $\lambda \in Res_{ar}(\Delta_T),$ $(\lambda I-\Delta_T)^{-1}=(\lambda I -\Delta_{T,\jjt})^{-1}$ and the $\jjt$-valued  map $\lambda \to (\lambda I-\Delta_{T,\mb})^{-1}$ is holomorphic on $Res_{ar}(\Delta_T).$ 

If $Re(\lambda)>log(\rho_{\bf T})\ge \lim_{t\to +\infty}{log\Vert T(t)\Vert_{\mathcal M (\mb)}\over t},$ then if $u \in \jjt, v \in \mb$, we have $$(\lambda I -\Delta_T)^{-1}uv=((\lambda I -\Delta_{T,\mb})^{-1}u)v=\int_0^{+\infty}e^{-\lambda t}T(t)uvdt.$$ Since $uv\in \Omega(\mb)$ for $u\in \Omega(\jjt),v\in \Omega(\mb),$ $u\mb$ is dense in $\mb$ for $u\in \Omega(\jjt),$
and we obtain $(\lambda I - \Delta_T)^{-1}=\int_0^{+\infty}e^{-\lambda t}T(t)dt \in \jjt,$ where the Bochner integral is computed with respect to the strong operator topology on $\mathcal M(\mb),$ so that $\Vert (\lambda I - \Delta_T)^{-1}\Vert_{\mathcal M(\mb)}\le \int_0^{+\infty}e^{-Re(\lambda)t}\Vert T(t)\Vert _{{\mathcal M}(\mb)}dt.$

$\square$

If we consider $\Delta_{T}$ as a quasimultiplier on $\mb,$ the fact that $(\lambda I -\Delta_{T})^{-1}\in {\mathcal M}(\mb)$ is the inverse of $\lambda I -\Delta_{T}$ for $\lambda \in Res(\Delta_{T})$ means that $(\lambda I -\Delta_{T})^{-1}v \in {\mathcal D}_{\Delta_{{T},\mb}}$ and that $(\lambda I-\Delta_{T})\left ((\lambda I -\Delta_{ T})^{-1}v\right )=v$ for every $v \in \mb,$ and that if $w \in {\mathcal D}_{\Delta_{{T},\mb}},$ then 
$(\lambda I -\Delta_{T})^{-1}\left ( (\lambda I -\Delta_{\bf T})w\right )=w.$ The situation is slightly more complicated if we consider $\Delta_{ T}$ as a quasimultiplier on $\ma$ when $\lim \sup_{t \to 0^+}\Vert T(t)\Vert =+\infty.$ In this case the domain ${\mathcal D}_{(\lambda I -\Delta_{T})^{-1},\ma}$ of $(\lambda I -\Delta_{T})^{-1}\in {\mathcal QM}(\ma)$ is a proper subspace of $\ma$ containing $\mathcal L_T \supset \cup_{t>0}T(t)\ma,$ and we have $(\lambda I -\Delta_{T})^{-1}v \in {\mathcal D}_{\Delta_{{T},\ma}}$ and $(\lambda I-\Delta_{T})\left ((\lambda I -\Delta_{T})^{-1}v\right )=v$ for every $v \in {\mathcal D}_{(\lambda I -\Delta_{\bf T})^{-1},\ma}.$ Also if $w \in {\mathcal D}_{\Delta_{{\bf T},\ma}},$ then $(\lambda I -\Delta_{\bf T})w\in {\mathcal D}_{(\lambda I -\Delta_{\bf T})^{-1},\ma},$ and we have 
$(\lambda I -\Delta_{ T})^{-1}\left ( (\lambda I -\Delta_{T})w\right )=w.$

In order to interpret $(\lambda I -\Delta_{\bf T})^{-1}$ as a partially defined operator on $\ma$ for $Re(\lambda) >log(\rho_{\bf T}),$ we can use the formula

\begin{equation}(\lambda I -\Delta_{\bf T})^{-1}v=\int_0^{+\infty}e^{-\lambda t}T(t)vdt \ \ \ \ (v \in \mathcal L_T),\end{equation} 

which defines a quasimultiplier on $\ma$ if we apply it to some $v \in \Omega(\ma)$ such that $\lim \sup_{t\to 0^+}\Vert T(t)u\Vert<+\infty.$ The fact that this quasimultiplier is regular is not completely obvious but follows from the previous discussion since $(\lambda I -\Delta_{\bf T})^{-1}\in {\mathcal M }(\mb)\subset {\mathcal QM}_r(\ma).$ Notice that since $\cup_{t>0}T(t)\ma$ is dense in $(\dt, \Vert . \Vert_{\dt}),$ $(\lambda I -\Delta_{T})^{-1}$ is characterized by the simpler formula

\begin{equation}(\lambda I -\Delta_{\bf T})^{-1}T(s)v=e^{\lambda s}\int_s^{+\infty}e^{-\lambda t}T(t)vdt \ \ \ \ (s>0, v \in \ma).\end{equation}

  \section{The generator of a holomorphic semigroup and its resolvent}
  
  Let $a<b\le a+\pi.$ In this section we consider a holomorphic semigroup $T=(T(\zeta))_{\zeta \in S_{a,b}}$ of multipliers on a weakly cancellative commutative Banach algebra ${\mathcal A}$ having dense principal ideals such that  $T(\zeta)\mathcal A$ is dense in ${\mathcal A}$ for some, or, equivalently, for every $\zeta\in S_{a,b}.$
  
  Denote by $\mathcal I_T$ the closed span of $\{T(\zeta)\}_{\zeta \in S_{a,b}}$, which is equal to the closed span of $\{T_{\zeta}(t)\}_{t>0}$ for $\zeta \in S_{a,b}.$ For $\zeta \in S_{a,b},$ set $T_{\zeta}=(T(t\zeta))_{t>0},$ let $\Phi_{T_{\zeta}}: \mathcal M_{\omega_{T_\zeta}}\to \mathcal M(\ma)$ be the homomorphism defined by (2). Set $\omega_T(\zeta)=\Vert T(\zeta)\Vert$ for $\zeta \in S_{a,b},$ denote by $\mathcal M_{\omega_T}(S_{a,b})$ the space of all measures $\mu$ on $S_{a,b}$ such that $\Vert \mu\Vert_{\omega_T}:=\int_{S_{a,b}}\omega_T(\zeta)d\vert \mu \vert(\zeta)<+\infty,$ which is a Banach algebra with respect to convolution. The convolution algebra $L^1_{\omega_T}(S_{a,b})$ is defined in a similar way and will be identified to the closed ideal of $\mathcal M_{\omega_T}$ consisting of measures which are absolutely continuous with respect to Lebesgue measure. Define $\Phi_T: \mathcal M_{\omega_T}\to \mathcal I_T \subset \mathcal M(\ma)$ by the formula
  
  $$\Phi_T(\mu)=\int_{S_{a,b}}T(\zeta)d\mu(\zeta),$$
  
  which is well-defined since the map $\zeta \to T(\zeta)$ is continuous with respect to the norm topology on $\mathcal M(\ma)$ and since $\mathcal I_T$ is separable.

Let $\zeta \in S_{a,b}.$ Since the semigroup $T_\zeta$ is continuous with respect to the norm topology on $\mathcal M(\ma),$ a standard argument shows that we have, for every Dirac sequence $(f_n)_{n\ge 1},$

  $$\lim \sup_{n\to +\infty}\left \Vert \int_{0}^{+\infty}(f_n*\delta_s)(t)T(t\zeta)dt -T(s\zeta)\right \Vert=0,$$ 
and so $T(s\zeta)\in {\mathcal I}_{T_{\zeta}}$ for every $s>0,$ which implies that ${\mathcal I}_{T_{\zeta}}=\mathcal I_T$, and a similar argument shows that $\mathcal I_T$ is the closure in $(\mathcal M(\ma), \Vert .\Vert_{\mathcal M(\ma)})$ of $\Phi_T(\mathcal M_{\omega_T}(S_{a,b}))$, as well as the closure of $\Phi_{T}(L^1_{\omega_T}(S_{a,b}))$ and the closure of $\Phi_{T_{\zeta}}(L^1_{\omega_{T_\zeta}})$ in $(\mathcal M(\ma), \Vert .\Vert_{\mathcal M(\ma)})$, and the notation $\mathcal I_T$ is consistent with the notation used to denote the Arveson ideal associated to a strongly continuous semigroup of multipliers on the half-line.
  
The following interpretation of the generator of a holomorphic semigroup as a quasimultiplier follows the interpretation given in  \cite{cep0} in the case where $\mathcal A=\iit.$ 

\begin{prop} Set

$$\Delta_ {T, \ma}:=S_{T'(\zeta_0)u_0/T(\zeta_0)u_0}\in \qm(\ma),$$

where $\zeta_0 \in S_{a,b},$ $u_0 \in \Omega(\ma).$

Then this definition does not depend on the choice of $\zeta_0$ and $u_0,$ and we have, for $\zeta \in S_{a,b},$

\begin{equation} \Delta_{T_{\zeta}, \ma}=\zeta \Delta_{T, \ma},\end{equation}

where the generator $\Delta_{T_{\zeta}, \ma}$ of the semigroup $T_{\zeta}$ is the quasimultiplier on $\ma$ introduced  in definition 4.1.

Moreover If $T_1=(T_1(\zeta))_{\zeta \in S_{a,b}}$ and $T_2=(T_2(\zeta))_{\zeta \in S_{a,b}}$ are two holomorphic semigroups of multipliers on $\ma$ such that $T_1(\zeta)\ma$ and $T_2(\zeta)\ma$ are dense in $\ma$ and such that $T_1(\zeta)T_2(\zeta)=T_2(\zeta)T_1(\zeta)$ for $\zeta \in \sab,$ then we have

$$\Delta_{T_1T_2,\ma}=\Delta_{T_1,\ma}+\Delta_{T_2,\ma}.$$

\end{prop}

Proof: We have, for $\zeta \in S_{a,b},$ 

$$T'(\zeta_0)T(\zeta)=T'(\zeta_0+\zeta)=T'(\zeta)T(\zeta_0),$$

and so the definition of $\Delta_ {T, \ma}$ does not depend on the choice of $\zeta_0,$ and an easy argument given in the comments following definition 4.1 shows that this definition does not depend on the choice of $u_0\in \Omega(\ma)$ either.

Now let $\zeta_0 \in S_{a,b}$, and le $f \in {\mathcal C}^1([0,+\infty))\cap \Omega(L^1_{\omega_{T_{\zeta_0}}})$ such that $\int_0^{+\infty} \vert f(t)\vert \Vert T(t\zeta_0)\Vert dt<\infty$ and $\int_0^{+\infty} \vert f'(t)\vert \Vert T(t\zeta_0)\Vert dt<\infty$ satisfying $f(0)=0.$ We have, integrating by parts, since $\lim_{p\to +\infty}\vert f(n_p)\vert \Vert T(n_p\zeta_0)\Vert=0$ for some strictly increasing sequence $(n_p)_{p\ge 1}$ of integers,

$$T(\zeta_0)\int_0^{+\infty} f'(t)T(t\zeta_0)dt =\lim_{p\to +\infty}\int_0^{n_p} f'(t)T(\zeta_0+t\zeta_0)dt$$ $$=\lim_{p\to +\infty}\left ( \left [f(t)T(\zeta_0+t\zeta_0)\right ]_0^{n_p} -\zeta\int_0^{n_p}f(t)T'(\zeta_0+t\zeta_0 )dt \right )$$ $$=-\zeta T'(\zeta_0)\int_0^{+\infty}f(t)T(t\zeta_0)dt,$$

and formula (11) follows since  $\left (\int_0^{+\infty}f(t)T(t\zeta_0)dt\right )u=\phi_{T_{\zeta_0}}(f)u\in \Omega(\ma)$ for $u \in \Omega(\ma).$ The last assertion follows immediately from the Leibnitz rule . $\square$


The following corollary follows then from proposition 5.4.

  \begin{cor} (i) Let $u \in \ma,$ and let $\zeta \in S_{a,b}.$ If $\lim_{t\to 0^+}\left \Vert \frac{T(t\zeta)u-u}{t}-v\right \Vert =0$ for some $v \in \ma,$ then  $u \in {\mathcal D}_{\Delta_{T, \ma}},$ and $\zeta \Delta_{T, \ma}u=v.$
  
  (ii) Conversely if $\lim \sup_{t\to 0^+}\left \Vert T(t\zeta) \right \Vert <+\infty,$ then $\lim \sup_{t\to 0^+} \left \Vert {T(t) u -u \over t} -\zeta \Delta_{T,\ma}u\right \Vert =0$ for every $u \in {\mathcal D}_{\Delta_{T, \ma}}.$
  
  \end{cor}

In the remaining of the section we will denote by $\mb$ a normalization of $\ma$ with respect to the semigroup $T,$ see definition 4.3. Since $\qm(\ma)$ is isomorphic to $\qm(\mb)$, we can consider the generator $\Delta_{T,\ma}$ as a quasimultiplier on $\mb$, and it follows immediately from definition 7.1 that this quasimultiplier on $\mb$ is the generator of the semigroup $T$ considered as a semigroup of multipliers on $\mb.$ From now on we will thus set $\Delta_T=\Delta_{T,\ma}=\Delta_{T,\mb}.$ Applying corollary 7.2 to $T$ and $\mb,$ we obtain

  \begin{cor} (i) Let $u \in \mb$ Then the following conditions imply each other
  
  (i) There exists $\zeta_0 \in S_{a,b}$ and $v \in  \mb$ such that $\lim_{t\to 0^+}\left \Vert \frac{ T(\zeta_0 t)u -u}{t}-v\right \Vert_{\mb}=0,$
  
  (ii) $u \in {\mathcal D}_{\Delta_T,\mb},$
  
 and in this situation $\lim_{t\to 0^+}\left \Vert \frac{ T(\zeta t)u -u}{t}-\zeta \Delta_Tu\right \Vert_{\mb}=0$ for every $\zeta \in S_{a,b}.$

  \end{cor}

Denote by  $\wiit$ the space of characters on $\iit,$ equipped with the usual Gelfand topology. If $\chi \in \wiit,$ the map $\zeta \to \chi(T(\zeta))$ is holomorphic on $S_{a,b},$ and so there exists a unique complex number $c_{\chi}$ such that $\chi(T(\zeta))=e^{\zeta c_{\chi}}$ for $\zeta \in S_{a,b}.$ We see as in section 5 that there exists a unique character $\tilde \chi$ on $\qm(\iit)$ such that $\tilde \chi_{|_{\iit}}=\chi,$ and since $\Delta_{T_{\zeta},\iit}=\zeta\Delta_{T,\iit}$ it follows from proposition 5.7 and proposition 7.1 that $\chi(T(t\zeta))=e^{t\tilde \chi(\Delta_{T_{\zeta},\iit})}= e^{t\zeta \tilde \chi(\Delta_{T,\iit})}$ for $\zeta \in S_{a,b}, t>0,$ and so $c_{\chi}=\tilde \chi(\Delta_{T,\iit}).$

Since $\Delta_{T_{\zeta},\iit}=\zeta\Delta_{T,\iit}$ for $\zeta \in S_{a,b},$ we deduce from corollary 5.8 and proposition 6.1 the following result.

\begin{prop} Let $T=(T(\zeta)_{\zeta \in S_{a,b}}\subset {\mathcal M}(\ma)$ be a holomorphic semigroup. Set $\sigma_{ar}(\Delta_{T,\iit})=\{ \tilde \chi (\Delta_{T,\iit})\}_{\chi \in \wiit},$ with the convention $\sigma_{ar}(\Delta_{T,\iit})=\emptyset$ if the semigroup is quasinilpotent, and set $Res_{ar}(\Delta_T)= Res_{ar}(\Delta_{T,\iit})=\C \setminus \sigma_{ar}(\Delta_{T,\iit}).$ Let $\mb$ be a normalization of $\ma$ with respect to the holomorphic semigroup $T,$ and let ${\mathcal I}_{T,\mb}$ be the closed subalgebra of $\mathcal M(\mb)$ generated by the semigroup.

(i) The set $\Lambda_{t, \zeta}:=\{ \lambda \in \sigma_{ar}(\Delta_T,\iit) \ | \ Re(\lambda \zeta) \le t\}$ is compact for $\zeta \in S_{a,b}, t\in \R.$

(ii) The quasimultiplier $\lambda I -\Delta_T$ has an inverse $(\lambda I -\Delta_T)^{-1} \in {\mathcal I}_{T,\mb}\subset {\mathcal M}(\mb)\subset \qm_r(\ma)$ for $\lambda \in Res_{ar}(\Delta_T),$ and
the map $\lambda \to (\lambda I -\Delta_T)^{-1}$ is a holomorphic map from  $Res_{ar}(\Delta_T)$ into ${\mathcal I}_{T,\mb}.$ 

(iii) If $\zeta \in S_{a,b},$ then $\lambda \in Res_{ar}(\Delta_T)$ for $Re(\lambda \zeta) > \lim_{t\to +\infty}\frac{log\left ( \left \Vert T(t\zeta)\right \Vert \right )}{t},$ and we have

\begin{equation} (\lambda I -\Delta_T)^{-1}=\int_0^{\zeta.\infty}e^{-s \lambda}T(s)ds,\end{equation}

so that

\begin{equation} \Vert  (\lambda I -\Delta_T)^{-1}\Vert_{{\mathcal M}(\mb)}\le \vert \zeta \vert \int_0^{+\infty}e^{-tRe(\lambda \zeta)}\Vert T(t\zeta)\Vert _{{\mathcal M}(\mb)}dt. \end{equation}
\end{prop}

Proof: (i) Let $\zeta \in S_{a,b},$ $t>0,$ and set $V=\{ \lambda \in\zeta  \sigma_{ar}(\Delta_{T,\bbt}) \ | \ Re(\lambda) \le t\}=\{ \lambda \in \sigma_{ar}(\Delta_{T_{\zeta},\iit}) \ | \ Re(\lambda) \le t\}.$ It follows from corollary 4.7 that $V$ is compact, and so $\Lambda_{t, \zeta}=\zeta^{-1}V$ is compact.

(ii) Fix $\zeta_0 \in S_{a,b}.$ We have $\lambda I -\Delta_{T}=\lambda I -\zeta_0^{-1}\Delta_{T_{\zeta_0}}=\zeta_0^{-1}\left ( \lambda \zeta_0I -\Delta_{T_{\zeta_0}}\right ).$ If $\lambda \in Res_{ar}(\Delta_T),$ then $\lambda I-\Delta_T$ is invertible in $\qm(\ma),$ and $(\lambda I-T)^{-1}=\zeta_0^{-1}(\lambda \zeta_0 I-\Delta_{T_{\zeta_0}})^{-1}\in {\mathcal I}_{T,\mb}\subset \mathcal M(\mb)\subset \qm_r(\ma),$ since in this situation $\lambda \zeta_0\in Res(\Delta_{T_{\zeta_0}}),$ and it follows also from proposition 6.1 that the ${\mathcal I}_{T,\mb}$-valued map $\lambda
\to (\lambda I-T)^{-1}=\zeta_0^{-1}(\lambda \zeta_0 I-\Delta_{T_{\zeta_0}})^{-1}$ is holomorphic on $Res_{ar}(\Delta_T).$

(iii) This follows from proposition 6.1 applied to $\lambda \zeta$ and $T_{\zeta}.$ $\square$

\section{Multivariable functional calculus for holomorphic semigroups associated to linear functionals}

In the following definition, we write by convention $T_j(0)=I$ for $1\le j \le k.$ Set  $\sigma \zeta=\sigma_1 \zeta_1 + \dots + \sigma_k \zeta_k$ for $\sigma=(\sigma_1,\dots, \sigma_k), \zeta=(\zeta_1,\dots, \zeta_k)\in \C^k.$

Let $a=(a_1,\dots, a_k)\in \R^k, b=(b_1,\dots, b_k)\in \R^k$ such that $a_j\le b_j \le a_j+\pi$ for $1\le j \le k.$ As in appendix 2, we set $M_{a,b}= \{(\alpha,\beta)\in \R^k\times \R^k \ | \ a_j<\alpha_j\le \beta_j <b_j \ \mbox{if} \ a_j < b_j,  \alpha_j=\beta_j=a_j \ \mbox{if} \ a_j=b_j\}.$ 

\begin{defn}: Let $a=(a_1,\dots,a_p) \in \R^k, b=(b_1,\dots,b_p) \in \R^k$ such that $a_j \le b_j \le a_j+\pi$ for $j \le k,$ let $\ma$ be a weakly cancellative commutative Banach algebra with dense principal ideals, and let $T=({T_1,\dots,T_k})$ be a family of semigroups of multipliers on $\ma$ which possesses the following properties

$$\left \{ \begin{array}{l} T_j=(T_j(\zeta))_{\zeta \in (0, a_j.\infty)} \ \mbox{is strongly continuous on} \  (0, e^{ia_j}.\infty), \ \mbox{and} \ \cup_{t>0}T(te^{ia_j})\ma \\ \mbox{is dense}  \mbox{ in} \ \ma \ \mbox{if} \ a_j=b_j,\cr 
 T_j=(T(\zeta))_{\zeta \in S_{a_j,b_j}} \ \mbox{is holomorphic on} \ S_{a_j,b_j}, \ \mbox{and} \ T(\zeta)\ma \ \mbox{is dense in} \ \ma \\
 \mbox{for every} \ \zeta \in S_{{a_j,b_j}} \ \mbox{if }\ a_j<b_j.\end{array} \right .$$
 
 For $\zeta=(\zeta_1,\dots,\zeta_k)\in \cup_{(\alpha, \beta)\in M_{a,b}}\overline S_{\alpha,\beta}$ set 
 
 $$T(\zeta)=T_1(\zeta_1)\dots T_k(\zeta_k).$$

A subalgebra $\mb$ of $\qm(\ma)$ is said to be a normalization of $\ma$ with respect to $T$ if the following conditions are satisfied

(a) $(\mb, \Vert . \Vert_{\mb})$ is a Banach algebra with respect to a norm $\Vert .\Vert_{\mb}$ satisfying $\Vert u\Vert_{\mb}\le \Vert u \Vert_{\ma}$ for $u \in \ma,$ and there exists a family $(w_1,\dots,w_k)$ of elements of $\Omega(\ma)$ such that the inclusion map $j:\ma \to \mb$ is a $\qm$-homomorphism with respect  to $T_1(\zeta_1)\dots T_k(\zeta_k)w_1\dots w_k$ for every family  $(\zeta_1,\dots,\zeta_k)$ of complex numbers such that $\zeta_j\in S_{a_j,b_j}$ if $a_j<b_j$ and such that $\zeta_j=0,$  if $a_j=b_j.$

(b) $\tilde j(\mathcal M(\ma))\subset {\mathcal M}(\mb)),$ and $\Vert \tilde j(R)\Vert_{\mathcal M (\mb)}\le \Vert R\Vert_{\mathcal M (\ma)}$ for every $R\in \mathcal M(\ma),$ where $\tilde j: \qm(\ma)\to \qm(\mb)$ is the pseudobounded isomorphism associated to $j$ in proposition 2.2 (ii).

(c) $\lim \sup_{\stackrel{\zeta \to 0}{_{\zeta \in S_{\gamma,\delta}}}}\Vert T(\zeta)\Vert_{\mathcal M(\mb)}<+\infty$ for $a_j<\gamma<\delta<b_j$ if $a_j<b_j,$ and $\lim \sup_{t\to 0^+}\Vert T(te^{ia_j})\Vert_{\mathcal M(\mb)}<+\infty$ if $a_j=b_j.$

\end{defn}

It follows from proposition 3.2 and proposition 4.2 that there exists a normalization $\mb_1$ of $\ma$ with respect to $T_1.$ Also if $\mb_m$ is a normalization of $\ma$ with respect to $(T_1,\dots,T_m)$ and if $\mb_{m+1}$ is a normalization of $\mb_m$ with respect to $T_{m+1},$ it follows  from proposition 2.4 and definitions 3.3 and 4.3 that $\mb_{m+1}$ is a normalization of $\ma$ with respect to $(T_1,\dots,T_{m+1}).$ It is thus immediate to construct a normalization of $\ma$ with respect to $T$ by a finite induction. Notice that if $\mb$ is a normalization of $\ma$ with respect to $T,$ then $\mb$ is a normalization of $\ma$ with respect to $T_\sigma:=(T(t\sigma))_{t>0}$ for every $\sigma \in \cup_{(\alpha, \beta)\in M_{a,b}}\overline S_{\alpha,\beta}$.

 Since $\cup_{t>0}T(te^{ia_j})\mb$ is dense in $\mb,$ when $a_j=b_j,$ and since $T(\zeta)\ma$ is dense in $\ma$ for $\zeta \in S_{a_j,b_j}$ if $a_j<b_j,$ it follows from condition (c) of definition 10.1 that the map $\zeta \to T(\zeta)u_1\dots u_k$ is continuous on $\overline S_{\alpha,\beta}$ for $(\alpha,\beta) \in M_{a,b},$ $u_1,\dots, u_k \in \mb.$ Since $u_1\dots u_k \in \Omega(\mb)$ for $u_1, \dots, u_k\in \Omega(\mb),$ it follows again from condition (c) of definition 10.1 that the map $\zeta \to T(\zeta)u$ is continuous on $\overline S_{\alpha,\beta}$ for $(\alpha,\beta) \in M_{a,b}$ for every $u\in \mb.$ Let $(\alpha, \beta)\in M_{a,b}$ and assume that $a_j<b_j.$ Since the semigroup $T_j$ is holomorphic on $S_{a_j,b_j}$  the map $\eta \to T(\zeta_1, \zeta_2, \dots, \zeta_{j-1}, \eta, \zeta_{j+1}, \zeta_k)u$ is holomorphic on $S_{\alpha_j,\beta_j}$ for every $(\zeta_1, \dots, \zeta_{j-1},\zeta_{j+1}, \zeta_k)\in \Pi_{\stackrel{1\le s \le k}{_{ s \neq j}}}\overline {S}_{\alpha_s, \beta_s}.$ 

Notice that if $u\in \mb,$ where $\mb$ is a normalization of $\ma$ with respect to $T,$ then the closed subspace $\mb_{T,u}$ spanned by the set $\{T(\zeta)u \ | \ \zeta \in S_{a,b}\}$ is separable, and so the function $\zeta \to T(\zeta)u$ takes its values in a closed separable subspace of $\mb.$





With the convention $T_j(0)=I$ for $1\le j \le k,$ we see that if $(\alpha, \beta)\in M_{a,b}$ and if $\lambda \in \cup_{(\gamma,\delta)\in M_{a-\alpha, b-\beta}}\overline S_{\gamma, \delta}$ then $T_{(\lambda)}: \zeta \to T(\lambda \zeta)=(T_1(\lambda_1\zeta_1), \dots,T_k(\lambda_k\zeta_k))$ is well-defined for $\zeta \in \overline S_{\alpha,\beta}. $
\begin{prop} 
Let $(\alpha, \beta)\in M_{a,b}.$ For  $\lambda\in \cup_{(\gamma, \delta) \in M_{a-\alpha, b -\beta}}\overline S_{\gamma, \delta},$ denote by $N(T, \lambda, \alpha, \beta)$ the set of all $z \in \C^k$ such that

$$\lim \sup_{t\to +\infty}\vert e^{tz_j e^{i\omega}}\vert \Vert T_j(t\lambda_je^{i\omega})\Vert <+\infty \ \ \mbox{for} \ \alpha_j \le \omega \le \beta_j, 1 \le j \le k $$,

and denote by $N_0(T, \lambda, \alpha, \beta)$ the set of all $z \in \C^k$ such that

$$\lim_{t\to +\infty}\vert e^{tz_j e^{i\omega}}\vert \Vert T_j(t\lambda_je^{i\omega})\Vert =0 \ \ \mbox{for} \ \alpha_j \le \omega\le  \beta_j, 1 \le j \le k.$$

Then $z \in N(T, \lambda, \alpha, \beta)$ if and only if $\lim \sup_{t\to +\infty}\vert e^{tz_j e^{i\alpha_j}}\vert \Vert T_j(t\lambda_je^{i\alpha_j})\Vert <+\infty$ and
$\lim \sup_{t\to +\infty}\vert e^{tz_j e^{i\beta_j}}\vert \Vert T_j(t\lambda_je^{i\beta_j})\Vert <+\infty$ for $1\le j \le k,$ 

Also $z \in N_0(T, \lambda, \alpha, \beta)$ if and only if $Re(z_j e^{i\alpha_j})<
-\lim_{t\to +\infty}{log\Vert T(t\lambda_je^{i\alpha_j})\Vert \over t}$ and $Re(z_j e^{i\beta_j})<
-\lim_{t\to +\infty}{log\Vert T(t\lambda_je^{i\beta_j})\Vert \over t}$ for $1\le j \le k.$
\end{prop}

Proof:  Let $j\le k$ such that $\alpha_j <\beta_j.$ If $\alpha_j \le \omega \le \beta_j,$ there exists $r_0>0$ and $s_0>0$ such that $e^{i\omega} =r_0e^{i\alpha_j}+s_0e^{i\beta_j},$ and we have, for $z_j \in \C,$

\begin{equation}\vert e^{tz_je^{i\omega}}\vert \Vert T_j(t\lambda_j e^{i\omega}\Vert \le \vert e^{r_0tz_je^{i\alpha j}}\vert \Vert T_j(r_0 t\lambda_j e^{i\alpha_j}\Vert \vert e^{s_0tz_je^{i\beta j}}\vert \Vert T_j(s_0 t\lambda_j e^{i\beta_j}\Vert,\end{equation}

and we see that $z \in N(T, \lambda, \alpha, \beta)$ if and only if $\lim \sup_{t\to +\infty}\vert e^{tz_j e^{i\alpha_j}}\vert \Vert T_j(t\lambda_je^{i\alpha_j})\Vert <+\infty$ and
$\lim \sup_{t\to +\infty}\vert e^{tz_j e^{i\beta_j}}\vert \Vert T_j(t\lambda_je^{i\beta_j})\Vert <+\infty$ for $1\le j \le k,$ which implies that $Re( z_j e^{i\alpha_j})\le
-\lim_{t\to +\infty}{log\Vert T(t\lambda_je^{i\alpha_j}\Vert \over t}$ and $Re( z_j e^{i\beta_j})\le
-\lim_{t\to +\infty}{log\Vert T(t\lambda_je^{i\beta_j}\Vert \over t}$ for $1\le j \le k.$ A similar argument shows that $z \in N_0(T, \lambda, \alpha, \beta)$ if and only if $Re(z_j e^{i\alpha_j})<
-\lim_{t\to +\infty}{log\Vert T(t\lambda_je^{i\alpha_j}\Vert \over t}$ and $Re(z_j e^{i\beta_j})<
-\lim_{t\to +\infty}{log\Vert T(t\lambda_je^{i\beta_j}\Vert \over t}$,  which implies that $Re( z_j e^{i\omega})<
-\lim_{t\to +\infty}{log\Vert T(t\lambda_je^{i\omega}\Vert \over t}$ for $\alpha_j\le \omega \le \beta_j,$ $1\le j \le k.$ $\square$

Notice that it follows from (16) and (17) that $N(T,\lambda, \alpha, \beta)-{\overline S}^*_{\alpha,\beta}\subset N(T,\lambda, \alpha, \beta)$ and that $N(T,\lambda,\alpha, \beta)-S^*_{\alpha,\beta}\subset N_0(T,\lambda,\alpha, \beta).$

Set  again $e_z(\zeta)=e^{z\zeta}$ for $z\in \C^k, \zeta \in \C^k.$ If $\mb$ is a normalization of $\ma$ with respect to $T,$ then $\sup_{\stackrel{\vert \zeta \vert \le 1}{\zeta \in \scab}}\Vert T(\zeta)\Vert_{\mathcal M(\mb)}<+\infty$ for $(\alpha, \beta) \in M_{a,b},$ and it follows from (42) that $\sup_{\zeta \in S_{\alpha,\beta}}\vert  e_{z}(\zeta)\vert \Vert T(\lambda \zeta)\Vert_{\mathcal M(\mb)}<+\infty$ for $z \in N(T, \alpha, \beta, \lambda)$ and $\lim_{\stackrel{\vert \zeta \vert \to +\infty}{\zeta \in \scab}}\vert  e_{z}(\zeta)\vert \Vert T(\lambda \zeta)\Vert_{\mathcal M(\mb)}=0$ if $z \in N_0(T, \lambda,\alpha, \beta)$ when $\lambda \in \cup_{(\gamma, \delta) \in M_{a-\alpha, b -\delta}}\overline S_{\gamma, \delta}.$ With the notations of appendices 1 and 2, we obtain the following result.

\begin{prop} Let $(\alpha, \beta)\in M_{a,b},$ and let  $\lambda\in \cup_{(\gamma, \delta) \in M_{a-\alpha, b -\beta}}\overline S_{\gamma, \delta}.$

\smallskip

(i) If $z \in N(T, \lambda, \alpha, \beta),$  then $ e_{z}T(\lambda .)u_{|_{\scab}} \in {\mathcal V}_{\alpha, \beta}(\mb),$ $\zeta_j-z_j\in Res_{ar}(\lambda_j\Delta_{T_j})$ 
 for $\zeta \in S^*_{a,b},$ $u\in \mb,$ $1 \le j \le k,$ and we have
 
 $$\FB( e_{z}T(\lambda .)u_{|_{\scab}})(\zeta)=(-1)^k((z_1-\zeta_1)I+\lambda_1 \Delta_{T_1})^{-1}\dots ((z_k-\zeta_k)I+\lambda_k \Delta_{T_k})^{-1}u.$$ 

\smallskip

(ii)If  $z \in N_0(T,\lambda, \alpha, \beta)$ then $ e_{z}T(\lambda .)u_{|_{\scab}} \in {\mathcal U}_{\alpha, \beta}(\mb),$ $\FB( e_{z}T(\lambda .)u_{|_{\scab}}$ has a continuous extension to $\overline S^*_{\alpha, \beta},$ $z_j+\zeta_j\in Res_{ar}(\lambda_j\Delta_j)$ for $1 \le j \le k,$
and we have,  for $\zeta \in \overline {S}^*_{a,b},$ $u\in \mb,$ $$\FB( e_{z}T(\lambda .)u_{|_{\scab}})(\zeta)=(-1)^k((z_1-\zeta_1)I+\lambda_1 \Delta_{T_1})^{-1}\dots ((z_k-\zeta_k)I+\lambda_k \Delta_{T_k})^{-1}u.$$

\end{prop}

Proof: It follows from the discussion above that $ e_{z}T(\lambda .)u_{|_{\scab}} \in {\mathcal V}_{\alpha, \beta}(\mb)$ if $z \in N(T, \lambda, \alpha, \beta), $ and that $ e_{z}T(\lambda .)u_{|_{\scab}} \in {\mathcal U}_{\alpha, \beta}(\mb)$ if $z  \in N_0(T,\lambda,\alpha, \beta).$ Let $z \in N(T, \lambda, \alpha, \beta),$ and let $u\in \mb.$ It follows from definition 10.3 (iii) that we have, for $\zeta \in S^*_{a,b},$

$$\FB( e_{z}T(\lambda .)u_{|_{\scab}})(\zeta) $$ $$
=\int_0^{e^{i\omega_1}.\infty}\dots \int_0^{e^{i\omega_k}.\infty}e^{(z_1-\zeta_1)\sigma_1+ \dots +(z_k-\zeta_k)\sigma_k}T_1(\lambda_1\sigma_1)\dots T_k(\lambda_k\sigma_k)ud\sigma_1\dots d\sigma_k,$$ 

where $\alpha _j \le \omega_j \le \beta_j$ and where $Re(\zeta_je^{i\omega_j})>0$ for $1\le j \le k.$

Since $Re((\zeta_j -z_j)e^{i\omega_j})>\lim_{t\to +\infty} {log(\Vert T(t\lambda_j \omega_j\Vert\over t},$ it follows from proposition 6.1 and proposition 7.4 that $\zeta_j-z_j\in Res_{ar}(\lambda_j\Delta_{T_j})$ for $j\le k,$ and that we have, for $v \in \mb,$

$$\int_{0}^{e^{i\omega_ j}.\infty}e^{(z_j-\zeta_j)\sigma_j}T(\lambda_j \sigma_j)vd\sigma_j=-((z_j-\zeta_j)I +\lambda_j\Delta_{T_j})^{-1}v.$$

Using Fubini's theorem, we obtain

$$\FB( e_zT(\lambda .)u_{|_{\scab}})(\zeta) $$ $$=-((z_1-\zeta_1)I+\lambda_1\Delta_{T_1})^{-1}\int_0^{e^{i\omega_2}.\infty}\dots \int_0^{e^{i\omega_k}.\infty}e^{(z_2-\zeta_2)\sigma_2- \dots +(z_k-\zeta_k)\sigma_k)}T_2(\lambda_2\sigma_2)\dots T_k(\lambda_k\sigma_k)ud\sigma_1\dots d\sigma_k$$
$$=\dots=(-1)^k((z_1-\zeta_1)I+\lambda_1 \Delta_{T_1})^{-1}\dots ((z_k-\zeta_k)I+\lambda_k \Delta_{T_k})^{-1}u.$$

 Since $N_0(T,\lambda, \alpha, \beta) \subset  N(T,\lambda, \alpha, \beta) -S^*_{\alpha, \beta},$ (ii) follows then from (i). $\square$

Recall that $\fab=(\cap_{z\in \C^k} e_{-z}\uab)'=\cup_{z \in \C^k}( e_{-z}\uab)'.$ If $\phi \in \fab,$ then $Dom(\FB(\phi))$ is the set of all $z \in \C^k$ such that $\phi \in ( e_{-z}\uab)'.$
In the following definition the action of $\phi \in \fab$ on an element $f$ of $ e_{-z}\vab(\mathcal B)$ taking values in a closed separable subspace of $\mb,$ where $z \in Dom(\FB(\phi),$ is defined according to definition
11.3. by the formula

$$<f,\phi>= < e_{z}f, \phi  e_{-z}>,$$

where $<g,\phi  e_{-z}>=< e_{-z}g,\phi>$ for $g \in \uab({\mathcal B}).$ It follows from the remarks following definition 11.3 that the above definition does not depend on the choice of $z.$

\begin{defn}  Let $(\alpha, \beta)\in M_{a,b},$ let $\lambda \in  \cup_{(\gamma, \delta) \in M_{a-\alpha, b -\beta}}\overline S_{\gamma, \delta},$ let $\phi \in \fab,$ and let $\mb$ be a normalization of $\ma$ with respect to $T.$ For $\zeta \in \overline S_{\alpha, \beta},$ set $T_{(\lambda)}(\zeta)=T(\lambda_1\zeta_1,\dots, \lambda_k\zeta_k)=T_1(\lambda_1\zeta_1)\dots T_k(\lambda_k\zeta_k),$ with the convention $T_j(0)=I.$

If $N(T, \lambda, \alpha, \beta)\cap Dom(\FB(\phi))\neq \emptyset,$ set, for $u \in \mb,$

$$<T_{(\lambda)},\phi>u= <T(\lambda .)u_{|_{\overline S_{\alpha, \beta}}},\phi>.$$

\end{defn}

For $(\alpha, \beta)\in M_{a,b}$, $z^{(1)}\in \C^k, z^{(2)}\in \C^k,$ we define as in definition 11.1 $\sup(z^{(1)},z^{(2)})$ to be the set of all $z\in \C^k$ such that $z+\overline S^*_{\alpha, \beta}=(z^{(1)}+\overline S^*_{\alpha, \beta})\cap (z^{(2)}+\overline S^*_{\alpha, \beta}),$ so that $\sup(z^{(1)},z^{(2)})$ is a singleton if $a_j<b_j$ for $j\le k.$ 

\begin{lem} If $\phi_1 \in \fab,$ $\phi_2\in  \fab,$ and if $N(T,\lambda, \alpha,\beta)\cap Dom(\fb(\phi_1))\neq \emptyset, N(T,\lambda, \alpha,\beta)\cap Dom(\fb(\phi_2))\neq \emptyset,$ then $\sup(z^{(1)}, (z^{(2)})\in  N(T,\lambda, \alpha,\beta)\cap Dom(\fb(\phi_1))\cap Dom(\fb(\phi_2))\subset N(T,\lambda, \alpha,\beta)\cap Dom(\fb(\phi_1*\phi_2))$ for $z^{(1)}\in 
N(T,\lambda, \alpha,\beta)\cap Dom(\fb(\phi_1)), z^{(2)}\in 
N(T,\lambda, \alpha,\beta)\cap Dom(\fb(\phi_2)),$ and the same property holds for $N_0(T,\lambda, \alpha,\beta).$

\end{lem}

Proof: Let $z^{(1)}\in 
N(T,\lambda, \alpha,\beta)\cap Dom(\fb(\phi_1)),$ let  $z^{(2)}\in 
N(T,\lambda, \alpha,\beta)\cap Dom(\fb(\phi_2)),$ let $z \in \sup(z^{(1)},z^{(2)}),$ and let $j \le k.$ There exists $s_1\in \{1,2\}$  and $s_2\in \{1,2\}$ such that $z_j \in \left (z^{(s_1)}_j+[0, e^{(-{\pi\over 2}-\alpha_j)i}.\infty)\right ) \cap\left ([z^{(s_2)}_j+[0,e^{({\pi\over 2}-\beta_j)i}.\infty)\right ),$ and it follows from (17) and from proposition 8.2 that $z \in N(T,\lambda, \alpha, \beta).$ The fact that $z\in Dom(\phi_1*\phi_2)$ follows from proposition 11.6. A similar argument shows that the same property holds for $N_0(T,\lambda, \alpha, \beta).$ $\square$
\begin{thm} Let $\ma$  be a weakly cancellative commutative Banach algebra with dense principal ideals, let $a,b \in \R^k,$ let $T=(T_1,\dots,T_k)$ be a family of semigroups of  multipliers on $\ma$ satisfying the conditions of definition 8.1, let $\mathcal B$ be a normalization of $\ma$ with respect to $T,$ let $(\alpha, \beta) \in M_{a,b}$ and let $\lambda \in  \cup_{(\gamma, \delta) \in M_{a-\alpha, b -\beta}}\overline S_{\gamma, \delta}.$ 

If $N(T,\lambda, \alpha,\beta)\cap Dom(\FB(\phi))\neq \emptyset$ for some $\phi \in \fab,$ then the following properties hold

(i) $<T_{(\lambda)},\phi>\in \mathcal M(\mb)\subset \QM_r(\ma),$ and we have, for $z\in N(T,\lambda, \alpha,\beta)\cap Dom(\FB(\phi)),$
if $\nu$ is a $z$-representative measure for $\phi,$

$$<T_{(\lambda)},\phi>= \int_{\scab}e^{z\zeta}T(\lambda \zeta)d\nu(\zeta),$$ 

where the Bochner integral is computed with respect to the strong operator topology on ${\mathcal M}(\mb),$ and if $\chi$ is a character on $\ma,$ then we have

$$\tilde \chi\left (<T_{(\lambda)},\phi>\right )=\fb(\phi)(-\lambda_1\tilde \chi(\Delta_{T_1}),\dots,-\tilde \chi(\Delta_{T_k})),$$

where $\tilde \chi$ denotes the unique character on $\qm(\ma)$ such that $\tilde \chi_{|_{\ma}}=\chi.$

(ii) $$\lim_{\stackrel{\eta \to (0, \dots,0), \eta \in \scab}{_{\epsilon \to(0,\dots,0), \epsilon \in \sscab}}} \Vert < e_{-\epsilon}T_{(\lambda)},\phi*\delta_{\eta}>u \ - <T_{(\lambda)}, \phi>u\Vert_{\mb}=0 \ \mbox{for} \ u \in \mb.$$

(iii) If $\alpha_j<\beta_j <\alpha_j+\pi$ for $1 \le j \le k,$ then we have, for $\eta \in \sab, \epsilon \in \ssab,$ 
$$ <e_{-\epsilon}T_{(\lambda)},\phi*\delta_{\eta}>=e^{-z\eta}\int_{\tilde \partial  \overline S_{\alpha, \beta}}e^{(z-\epsilon )\sigma}{\mathcal C}_z(\phi)(\sigma-\eta)T(\lambda \sigma)d\sigma$$

where the Bochner integral is computed with respect to the strong operator topology on $\mathcal M(\mb).$ 

\smallskip

(iv) In the general case, set $W_n(\zeta)=\Pi_{1\le j \le k} {n^2\over\left ( n+\zeta_je^{i{\alpha_j+\beta_j\over 2}}\right )^2}$ for $n\ge 1, \zeta=(\zeta_1,\dots, \zeta_n)\in \overline S^*_{\alpha \beta}.$ Then we have

$$<T_{(\lambda)},\phi>$$ $$=\lim_{\stackrel{\epsilon \to 0}{_{\epsilon \in S^*_{\alpha,\beta}}}}\left (\lim_{n\to +\infty}{(-1)^k\over (2i\pi)^k}\int_{z+\tilde \partial S^*_{\alpha, \beta}}W_n(\sigma -z)\fb(\phi)(\sigma)((\sigma_1-\epsilon_1)I +\lambda_1\Delta_T)^{-1}\dots ((\sigma_k-\epsilon_k)I +\lambda_k \Delta_{T_k})^{-1}d\sigma\right ),$$

where the Bochner integral is computed with respect to the norm topology on $\mathcal M(\mb).$

(v) If, further, $\int_{z+\tilde \partial S^*_{\alpha, \beta}}\vert \fb(\phi)(\sigma)\vert \vert d\sigma)\vert <+\infty,$ then we have

$$<T_{(\lambda)},\phi> =\lim_{\stackrel{\epsilon \to 0}{_{\epsilon \in S^*_{\alpha, \beta}}}}<e_{-\epsilon}T_{(\lambda)},\phi>=$$ $$\lim_{\stackrel{\epsilon \to 0}{_{\epsilon \in S^*_{\alpha, \beta}}}}{(-1)^k\over (2i\pi)^k}\int_{z+\tilde \partial  S^*_{\alpha, \beta}} \FB(\phi)(\sigma)((\sigma_1-\epsilon_1)I +\lambda_1\Delta_T)^{-1}\dots ((\sigma_k-\epsilon_k)I +\lambda_k \Delta_{T_k})^{-1}d\sigma,$$ 
where the Bochner integral is computed with respect to the norm topology on $\mathcal M(\mb).$

(vi) If the condition of (v) is satisfied for some $z \in N_0(T, \lambda, \alpha,\beta)\cap Dom(\fb(\phi)),$ then we have

$$ <T_{(\lambda)},\phi>=$$ $${(-1)^k\over (2i\pi)^k}\int_{z+\tilde \partial  S^*_{\alpha, \beta}} \FB(\phi)(\sigma)((\sigma_1I +\lambda_1\Delta_T)^{-1}\dots (\sigma_kI +\lambda_k \Delta_{T_k})^{-1}d\sigma.$$ 

(vii) If $\phi_1\in \fab, \phi_2 \in \fab,$ and if $N(T, \lambda, \alpha, \beta)\cap Dom({\FB(\phi_1}))\neq \emptyset$ and $N(T, \lambda, \alpha, \beta)\cap Dom({\FB(\phi_2}))\neq \emptyset,$  then $N(T, \lambda, \alpha, \beta)\cap Dom(\FB(\phi_1*\phi_2))\neq \emptyset,$  and

$$<T_{(\lambda)},\phi_1*\phi_2>=<T_{(\lambda)},\phi_1><T_{(\lambda)},\phi_2>.$$

\end{thm}

Proof: (i)  Let $z \in N(T,\lambda, \alpha, \beta) \cap Dom(\FB(\phi)),$ and set $m=\sup_{\zeta \in \sab}\vert e^{z\zeta}\vert \Vert T_{\lambda}(\zeta)\Vert_{\mathcal M(\mb)}<+\infty.$  We have  

$$\Vert <T_{\lambda},\phi>u \Vert _{\mathcal B} \le m\Vert \phi e_z\Vert _{\uab'}\Vert u \Vert_{\mb},$$

and so $<T_{\lambda},\phi>\in {\mathcal M}(\mb).$ The integral formula in (i) follows then immediately from the definition given in proposition 10.2 and from definition 11.3.

Assume that $\ma$ is not radical, let $\chi$ be a character on $\ma,$ and let $\tilde \chi$ be the unique character on $\qm(\ma)$  such that $\tilde \chi(u)=\chi(u)$ for every $u \in \ma.$ Set $f_n(t)=0$ if $0\le t <{1\over n+1}$ or if $t>{1\over n},$ and $f_n(t)=n(n+1)$ if ${1\over n+1}\le t \le {1\over n},$ and let $\zeta$ be an element of the domain of definition of $T_j.$ Set $T_{j,\zeta}:=(T_j(t\zeta))_{t>0}.$ Then $(f_n)_{n\ge 1}\subset L^1_{\omega_{T_{j,\zeta}}}(\R^+)$ is a Dirac sequence, and since the map $t\to T_j(t\zeta)u$ is continuous on $(0,+\infty),$ a standard argument shows that we have, for $s>0,$ $u \in \ma,$

$$\lim \sup_{n\to +\infty}\left \Vert \Phi_{T_{j,\zeta}}(f_n)T_j(s\zeta)u-T_j(s\zeta)u\right \Vert $$ $$=\lim \sup_{n\to +\infty}\left \Vert \int_0^{+\infty}(f_n*\delta_s)(t)T_j(t\zeta)udt -T_j(s\zeta)u\right \Vert=0.$$

Since $\cup_{t>0}T_j(t\zeta)(\ma)$ is dense in $\ma,$ there exist $n\ge 1$ such that $\tilde \chi(\Phi_{T_{j,\zeta}}(f_n))\neq 0,$ and the restriction of $\tilde \chi$ to the Arveson ideal $\mathcal I_{T_{j,\zeta}}$ is a character on $\mathcal I_{T_{j,\zeta}}.$  It follows then from proposition 5.7 and proposition 7.1 that $\tilde \chi(T_j(t\zeta))=e^{t\tilde \chi(\Delta_{T_{j, \zeta}})}=e^{t\zeta\tilde \chi(\Delta_{T_j})}$ for $t>0,$ and so $\tilde \chi(T_j(\zeta))=e^{\zeta\tilde \chi(\Delta_{T_j})}$ for every $\zeta$ in the domain of definition of $T_j.$ Let $u\in \mathcal M(\mb).$ By continuity,
we see that $\tilde \chi(T_j(\zeta)u)=e^{\zeta\tilde \chi (\Delta_{T_j})}\tilde \chi(u)$ for every $\zeta \in \overline S_{\alpha_j,\beta_j}.$ Set $\lambda\zeta \tilde \chi(\Delta_T)=\lambda_1\zeta_1\tilde \chi(\Delta_{T_1})+\dots +\lambda_k\zeta_k\tilde \chi(\Delta_{T_k}).$   Consider again $z \in N(T,\lambda, \alpha, \beta) \cap Dom(\FB(\phi)).$ Since $\fb(\phi e_{-z})(\zeta)=<e_{-\zeta},\phi e_{-z}>=<e_{-z-\zeta},\phi>=\fb(\phi)(\zeta+z)$ for $\zeta \in \overline S^*_{\alpha, \beta},$ we obtain

$$\tilde \chi (<T_{(\lambda)},\phi>)\tilde \chi(u)=\tilde \chi\left (  \int_{\scab}e^{z\zeta}T(\lambda \zeta)ud\nu(\zeta)\right)=\int_{\scab}e^{z\zeta}\tilde \chi(T(\lambda \zeta))\tilde \chi(u)d\nu(\zeta)$$

$$=\left ( \int_{\scab}e^{z\zeta}e^{\lambda \zeta \tilde \chi(\Delta_T)}d\nu(t)\right ) \tilde \chi(u)= \fb(\phi e_{-z})(-\lambda\chi(\Delta_T)-z)\tilde \chi(u)=\fb(\phi)(-\lambda \chi(\Delta_T))\tilde \chi(u),$$

which concludes the proof of (i), since $\tilde \chi(u)\neq 0.$

(ii) Let $u \in \mb,$ and set $f(\zeta)=T(\lambda \zeta)u$ for $\zeta \in \scab.$ Using definition 11.3, we see that (ii) follows from (34) applied to $f.$

(iii) Define $f$ as above. We have, for $ \epsilon \in \ssab, \eta \in \sab,$ $u\in \mb,$

$$ <e_{-\epsilon}T_{(\lambda)},\phi*\delta_{\eta}>u=<e_{-\epsilon}f, \phi*\delta_{\eta}>=<(e_{-\epsilon})_{\eta}f_{\eta},\phi>=e^{-\epsilon \eta}<e_{-\epsilon}f_{\eta},\phi>,$$
so (iii) follows from (37).

\smallskip

(iv) The result follows from  proposition 11.9 (i) applied to $T_{(\lambda)} u_{|_{ \scab}},$ $z$ and $\phi$ for $u \in \mb.$

(v) The result follows from proposition 11.9 (ii) applied to $T_{(\lambda)} u_{|_{ \scab}},$ $z$ and $\phi$ for $u \in \mb.$

\smallskip

(vi) Now assume that the condition of (v) is satisfied for some $z \in N_0(T,\lambda, \alpha, \beta)\cap Dom(\fb(\phi)).$ Then there exists $\epsilon \in S^*_{\alpha, \beta}$ such that
$z+\epsilon\in N_0(T,\lambda, \alpha, \beta),$ and $z+\epsilon \in Dom(\fb(\phi e_{\epsilon})).$ We have

$$\int_{\partial \scab}e^{Re(z\sigma)}\Vert T_{\lambda}(\sigma)\Vert_{\mathcal M(\mb)} \vert d\sigma\vert$$ $$\le \left (\sup_{\zeta \in \scab}e^{Re((z+\epsilon)\zeta)}\Vert T_{(\lambda)}(\zeta)\Vert_{\mathcal M(\mb)}\right ) \int_{\partial \scab}e^{-Re(\epsilon \sigma)}\vert d \sigma \vert <+\infty,$$ 

and (vi) follows from proposition 11.9 (iii) applied to $T_{(\lambda)} u_{|_{ \scab}}$ $z$ and $\phi$ for $u \in \mb.$

\smallskip

(vii) Now assume that $\phi_1\in \fab, \phi_2 \in \fab$ satisfy the hypothesis of (vi) with respect to $T$ and $\lambda,$ and let $z^{(1)}\in N(T,\lambda, \alpha, \beta)\cap Dom({\FB(\phi_1}))$ and $z^{(2)}\in N(T,\lambda, \alpha, \beta)\cap Dom(\FB(\phi_2)).$ Set $z =\sup(z^{(1)},z^{(2)}).$ It follows from lemma 8.5 that  
$z \in N(T, \lambda,\alpha, \beta)\cap Dom({\FB(\phi_1}))\cap Dom(\FB(\phi_2))\subset N(T, \lambda,\alpha, \beta)\cap Dom({\FB(\phi_1}*\phi_2)).$ Let $\nu_1$ be a $z$-representative measure for $ \phi_1$ and let $\nu_2$ be a $z$-representative measure for $\phi_2.$ Then $\nu_1*\nu_2$ is a $z$-representative measure for $\phi_1*\phi_2,$ and we have, for $u \in \mb,$

$$<T_{(\lambda)},\phi_1*\phi_2>u=\int_{\scab}e^{z \zeta}T(\lambda \zeta)ud(\nu_1*\nu_2)(\zeta)$$

$$= \int_{\scab\times \scab}e^{z(\zeta_1+\zeta_2)}T(z(\zeta_1+\zeta_2))ud\nu_1(\zeta_1)d\nu_2(\zeta_2)$$

$$\left [ \int_{\scab}e^{z \zeta_1}T(\lambda \zeta_1)d\nu_1(\zeta_1)\right ]\left [ \int_{\scab}e^{z \zeta_2}T(\lambda \zeta_2)ud\nu_2(\zeta_2)\right ]$$
$$=<T_{(\lambda)},\phi_1>\left (<T_{(\lambda)},\phi_2>u\right),$$

which proves (vii).

$\square$

Let $\mathcal G_{a,b}=\cup_{(\alpha, \beta)\in M_{a,b}}{\mathcal F}_{\alpha, \beta}$ be the dual space introduced in definition 11.2, which is an algebra with respect to convolution according to proposition 11.13. If $\phi_1 \in \mathcal F_{\alpha^{(1)},\beta^{(1)}},$ and if $\phi_1 \in \mathcal F_{\alpha^{(2)},\beta^{(2)}},$ where $(\alpha^{(1)},\beta^{(1)})\in M_{a,b}, (\alpha^{(2)},\beta^{(2)})\in M_{a,b},$ then $\phi_1*\phi_2$ is well-defined but in general the fact that $N(T, \lambda, \alpha^{(1)}, \beta^{(1)})\cap Dom(\phi_1) \neq \emptyset$ and $N(T, \lambda, \alpha^{(2)}, \beta^{(2)})\cap Dom(\phi_2) \neq \emptyset$ does not seem to imply that $N(T, \lambda, \inf(\alpha^{(1)}, \alpha^{(2)}), \sup(\beta^{(1)}, \beta^{(2)}))\cap Dom(\phi_1*\phi_2)\neq \emptyset,$ which prevents from obtaining a direct extension of (vi) to the case where $\phi_1\in \mathcal G_{a,b}, \phi_2 \in \mathcal G_{a,b}.$ This difficulty will be circumvented in the next section by using Fourier-Borel transforms.

\section{Multivariable functional calculus for holomorphic semigroups associated to holomorphic functions of several complex variables}

In the following definition, the generator $\Delta_{T_j}$ of the strongly continuous semigroup $T_j$ and its Arveson spectrum $\sigma_{ar}(\Delta_{T_j})$ are defined according to section 5 if $a_j=b_j,$ and the generator $\Delta_{T_j}$ of the holomorphic semigroup $T_j$ and its Arveson spectrum $\sigma_{ar}(\Delta_{T_j})$ are defined according to section 7 if $a_j<b_j.$ 

\begin{defn} Let  $a=(a_1,\dots,a_p) \in \R^k, b=(b_1,\dots,b_p) \in \R^k$ such that $a_j \le b_j \le a_j+\pi$ for $j \le k,$ let $\ma$ be a weakly cancellative commutative Banach algebra having dense principal ideals, and let $T=({T_1,\dots,T_k})$ be a family of semigroups of multipliers on $\ma$ satisfying the conditions of definition 8.1. 
Let $(\alpha, \beta) \in M_{a,b}$ and let $\lambda \in \cup_{(\gamma, \delta) \in M_{a-\alpha, b-\beta}}\overline S_{\gamma, \delta}.$ 

An open set $U \subset \C^k$ is said to be admissible with respect to $(T, \lambda, \alpha, \beta)$ if $U= \Pi_{1\le j \le k}U_j$ where the open sets $U_j\subset \C$ satisfy the following conditions for some $z =(z_1,\dots, z_k)\in N_0(T,\alpha, \beta, \lambda)$ 

(i)  $U_j+\overline S^*_{\alpha_j,\beta_j}\subset U_j$ 

(ii) $U_j \subset z_j + S^*_{\alpha_j,\beta_j} ,$ and $\partial U_j =(z_j +e^{(-{\pi\over 2}-\alpha_j)i}.\infty, z_j +e^{(-\alpha_j-{\pi\over 2})i}s_{0,j})\cup \gamma([0,1])\cup(z_j+e^{({\pi\over 2}-\beta_j)i}s_{1,j},z_j + e^{({\pi\over 2}-\beta_j)i}.\infty),$ where $s_{0,j}\ge 0, s_{1,j}\ge 0,$ and where $\gamma: [0,1] \to  z_j+\overline {S}^*_{\alpha_j,\beta_j} \setminus \left (e^{(-{\pi\over 2}-\alpha_j)i}.\infty, e^{(-\alpha_j-{\pi\over 2})i}s_{0,j})\cup (e^{({\pi\over 2}-\beta_j)i}s_{1,j}, e^{({\pi\over 2}-\beta_j)i}.\infty)\right )$ is a one-to-one piecewise-$\mathcal C^1$ curve such that 
$\gamma(0)= e^{(-\alpha_j-{\pi\over 2})i}s_{0,j}$ and $\gamma_1=e^{({\pi\over 2}-\beta_j)i}s_{1,j}.$

(iii)  $\lambda_j\sigma_{ar}(-\Delta_{T_j})= \sigma_{ar}(-\Delta _{T_j(\lambda_j.)})\subset U_j.$

\end{defn}

Conditions (i) and (ii) mean that $U$ is admissible with respect to $(\alpha, \beta)$ in the sense of definition 12.1 and that some, hence all  elements $z\in \C^k$ with respect to which $U$ satisfies condition (ii) of definition 12.1 belong to $N_0(T,\alpha, \beta, \lambda).$ Hence $U_j$ is a open half-plane if $\alpha_j=\beta_j$, and the geometric considerations about $\partial U_j$ made in section 12 when $\alpha_j < \beta_j$ apply.

For $\alpha=(\alpha_1,\dots,\alpha_k)\in \R^k, \beta=(\beta_1,\dots, \beta_k)\in \R^k,$ we will use as in appendix 3 the obvious conventions $$\inf(\alpha, \beta)=(\inf(\alpha_1,\beta_1),\dots , \inf(\alpha_k, \beta_k)), \ \sup(\alpha, \beta)=(\sup(\alpha_1,\beta_1),\dots , \sup(\alpha_k, \beta_k)).$$

\begin{prop} If  $U^{(1)}$ is admissible with respect to $(T, \lambda, \alpha^{(1)},\beta^{(1)})$ and if $U^{(2)}$ is admissible with respect to $(T, \lambda, \alpha^{(2)},\beta^{(2)}),$ then
$U^{(1)}\cap U^{(2)}$ is admissible with respect to $(T, \lambda, \inf(\alpha^{(1)},\alpha^{(2)}), \sup(\beta^{(1)},\beta^{(2)})).$

\end{prop}

Proof: Set $\alpha^{(3)}= \inf(\alpha^{(1)},\alpha^{(2)}),\beta^{(3)}=\sup(\beta^{(1)},\beta^{(2)}),$ and set $U^{(3)}=U^{(1)}\cap U^{(2)}.$ Then $\left [\cup_{(\gamma, \delta) \in M_{a-\alpha^{(1)}, b-\beta^{(1)}}}\overline S_{\gamma, \delta}\right ]\cap \left [\cup_{(\gamma, \delta) \in M_{a-\alpha^{(2)}, b-\beta^{(2)}}}\overline S_{\gamma, \delta}\right ] \subset \left [\cup_{(\gamma, \delta) \in M_{a-\alpha^{(3)}, b-\beta^{(3)}}}\overline S_{\gamma, \delta}\right ],$ so it makes sense to check whether $U^{(1)}\cap U^{(2)}$ is admissible with respect to $(T, \lambda, \alpha^{(3)},\beta^{(3)}).$ The fact that $U^{(3)}$ satisfies (i) and (ii) follows from proposition 12.2, and the fact that $U^{(3)}$ satisfies (iii) is obvious. $\square$

If an open set $U \subset \C^k$ is admissible with respect to $(T, \lambda, \alpha, \beta)$, we denote as in section 12 by $H^{(1)}(U)$ the set of all holomorphic functions $F$ on $U$ 
satisfying the condition

$$\Vert F\Vert_{H^{(1)}(U)}:=\sup_{\epsilon \in S^*_{\alpha, \beta}}\int_{\sigma \in \tilde \partial U +\epsilon}\vert F(\sigma)\vert \vert d\sigma\vert <+\infty.$$

Notice that $\cup_{(\alpha, \beta)\in M_{a,b}}\left (\cup_{(\gamma, \delta)\in M_{a-\alpha, b -\beta}}\overline S_{\beta, \gamma}\right )=\cup_{(\alpha, \beta)\in M_{a,b}}\overline S_{a-\alpha, b-\beta}.$ The inclusion $\cup_{(\gamma, \delta)\in M_{a-\alpha, b -\beta}}\overline S_{\beta, \gamma}\subset \overline S_{a-\alpha, b-\beta}$ for $(\alpha, \beta)\in M_{a,b}$ is obvious. Conversely assume that $\lambda \in \overline S_{a-\alpha,b-\beta}$ for some $(\alpha, \beta) \in M_{a,b}.$ If $a_j=b_j$ then $a_j=\alpha_j=\beta_j=b_j,$ and so $\lambda_j$ is a nonnegative real number.  In this situation set $\alpha'_j=\beta'_j=a_j, \gamma_j=\delta_j=0.$ If $a_j<b_j,$ then $a_j<\alpha_j\le \beta_j <b_j,$ and $a_j-\alpha_j \le arg(\lambda_j)\le b_j-\beta_j$ if $\lambda_j\neq 0.$ In this situation set $\alpha'_j={a_j+\alpha_j\over 2}, \gamma_j=a_j -\alpha_j, \delta_j=b_j-\beta_j$ and $\beta'_j={b_j+\beta_j\over 2}.$ Then $(\alpha', \beta')\in M_{a,b},$  $(\gamma,\delta) \in M_{a-\alpha',b-\beta'},$ and $\lambda \in \overline S_{\gamma,\delta},$ which concludes the proof of the reverse inclusion.

\begin{cor} For $\lambda \in \cup_{(\alpha, \beta)\in M_{a,b}}\overline S_{a-\alpha, b-\beta} =\cup_{{\alpha, \beta}\in M_{a,b}}\left (\cup_{(\gamma, \delta) \in M_{a-\alpha, b-\beta}}\overline S_{\gamma, \delta}\right ),$ denote by $\mathcal N_{\lambda}$ the set of all $(\alpha, \beta) \in M_{a,b}$ such that $\lambda \in  \cup_{(\gamma, \delta) \in M_{a-\alpha, b-\beta}}\overline S_{\gamma, \delta}$, and denote by $\mathcal W_{T,\lambda}$ the set of all open sets $U \subset \C^k$ which are admissible with respect to $(T, \lambda, \alpha, \beta)$ for some $(\alpha, \beta) \in  \mathcal N_{\lambda}.$ 

Then $\mathcal W_{T,\lambda}$ is stable under finite intersections, and $\cup_{U\in \mathcal W_{T,\lambda}}H^{(1)}(U)$ is stable under products.

\end{cor}

Proof: The first assertion follows from the proposition and the second assertion follows from the fact that the restriction of $F \in H^{(1)}(U)$ is bounded on $U+\epsilon$ if $U$ is admissible with respect to $(\alpha,\beta)\in M_{a,b}$ and if $\epsilon \in S^*_{\alpha, \beta},$ see corollary 12.4. $\square$

A set ${\mathcal E} \subset \cup_{U\in \mathcal W_{T,\lambda}}H^{(1)}(U)$ will be said to be bounded if there exists $U\in \mathcal W_{T,\lambda}$ such that  $\mathcal E\subset H^{(1)}(U)$ and such that $\sup_{F\in {\mathcal E}}\Vert F \Vert_{H^{(1)}(U)}<+\infty,$ and bounded subsets of  $\cup_{U\in \mathcal W_{T,\lambda}}H^{\infty}(U)$ are defined in a similar way. A homomorphism $\phi: \cup_{U\in \mathcal W_{T,\lambda}}H^{(1)}(U)\to \mathcal M(\mb)$ will be said to be bounded if $\phi(\mathcal E)$ is bounded for every bounded subset $\mathcal E$ of $\cup_{U\in \mathcal W_{T,\lambda}}H^{(1)}(U),$ and a homomorphism $\phi : \cup_{U\in \mathcal W_{T,\lambda}}H^{\infty}(U)\to \qm_r(\mb)=\qm_r(\ma)$ will be said to be bounded if $\phi(\mathcal E)$ is pseudobounded for every bounded subset $\mathcal E$ of $\cup_{U\in \mathcal W_{T,\lambda}}H^{\infty}(U).$

Similarly let $\mathcal S(U)$ be the Smirnov class on $U \in \mathcal W_{T,\lambda}$ introduced in definition 12.6. A set ${\mathcal E} \subset \cup_{U\in \mathcal W_{T,\lambda}}\mathcal S(U)$ will be said to be bounded if there exists $U\in \mathcal W_{T,\lambda}$ such that  $\mathcal E\subset \mathcal S(U)$ and such that $\sup_{F\in {\mathcal E}}\Vert FG\Vert_{H^{\infty}(U)}<+\infty$ for some strongly outer function $G\in H^{\infty}(U),$ and a homomorphism $\phi : \cup_{U\in \mathcal W_{T,\lambda}}\mathcal S(U)\to \qm(\mb)=\qm(\ma)$ will be said to be bounded if $\phi(\mathcal E)$ is pseudobounded for every bounded subset $\mathcal E$ of $\cup_{U\in \mathcal W_{T,\lambda}}\mathcal S(U).$

Let $U=\Pi_{j\le k}U_j \in \mathcal W_{T,\lambda},$ and let $(\alpha, \beta) \in \mathcal N_{\lambda}$ such that $U$ is admissible with respect to $(T,\lambda, \alpha, \beta).$ Let $\partial U_j$ be oriented from $e^{-{\pi\over 2}-\alpha_j}.\infty$ to $e^{{\pi\over 2}-\beta_j}.\infty$. This gives an orientation on the distinguished boundary $\tilde \partial U=\Pi_{j\le k}\partial U_j$ of $U,$ to be used in the following theorem.

\begin{thm} Let $a=(a_1,\dots,a_p) \in \R^k, b=(b_1,\dots,b_p) \in \R^k$ such that $a_j \le b_j \le a_j+\pi$ for $j \le k,$ let $\ma$ be a weakly cancellative commutative Banach algebra with dense principal ideals, let $T=(T_1,\dots, T_k)$ be a family of semigroups of multipliers on $\ma$
satisfying the conditions of definition 8.1 with respect to $(a,b)$ and $\ma$ and let $\mb$ be a normalization of $\ma$ with respect to $T.$ 

(i) For $\lambda \in \cup_{(\alpha,\beta)\in M_{a,b}} \overline S_{\alpha, \beta},$ $U\in \mathcal W_{T,\lambda},$ $F\in H^{(1)}(U),$ set

$$F(-\lambda_1\Delta_{T_1},\dots,-\lambda_k\Delta_{T_k})$$ $$={1\over (2i\pi)^k}\int_{\tilde \partial U +\epsilon}F(\zeta_1,\dots,\zeta_k)(\lambda_1\Delta_{T_1} +\zeta_1 I)^{-1}\dots (\lambda_1\Delta_{T_k} +\zeta_k I)^{-1}d\zeta_1\dots d\zeta_k,$$

 where $U$ is admissible to respect to $(T, \lambda, \alpha, \beta),$ with $(\alpha, \beta)\in \mathcal N_{\lambda},$ and where $\epsilon \in S^*_{\alpha, \beta}$ is such that $U+\epsilon \in \mathcal W_{T,\lambda}.$ Then this definition does not depend on the choice of $U$ and $\epsilon,$ and the map $F \to F(-\lambda_1\Delta_{T_1},\dots,-\lambda_k\Delta_{T_k})$ is a bounded algebra homomorphism
 from $\cup_{U\in \mathcal W_{T, \lambda}}H^{(1)}(U)$ into $\mathcal M(\mb)\subset \qm_r(\ma).$

(ii) For every $U\in \mathcal W_{T,\lambda}$ there exists $G \in H^{(1)}(U)\cap H^{\infty}(U)$ such that $G(-\lambda_1\Delta_{T_1},\dots, -\lambda_k\Delta_{T_k})(\mb)$ is dense in $\mb,$ and for every $F \in H^{\infty}(U)$ there exists a unique $R_F \in \qm_r(\mb)=\qm_r(\ma)$  satisfying $$R_FG(-\lambda_1T_1, \dots, -\lambda_kT_k)=(FG)(-\lambda_1T_1, \dots, -\lambda_kT_k) \ \ (G\in H^{(1)}(U)).$$ The definition of $R_F$ does not depend on the choice of $U,$ and if we set $F(-\lambda_1T_1,\dots, -\lambda_kT_k)=R_F$ the definition of $F(-\lambda_1\Delta_{T_1}, \dots, -\lambda_k\Delta_{T_k})$ agrees with the definition given in (i) if $F \in \cup_{U\in \mathcal W_{T,\lambda}}H^{(1)}(U),$ and the map $F \to F(-\lambda_1T_1,\dots, -\lambda_kT_k)$ is a bounded homomorphism from $\cup_{U \in \mathcal W_{T,\lambda}}H^{\infty}(U)$ into $\qm_r(\mb)=\qm_r(\ma).$

(iii) If $(\alpha, \beta) \in \mathcal N_{\lambda},$ if $\phi \in \fab,$ and if $N_0(T,\lambda, \alpha, \beta)\cap Dom(\fb(\phi))\neq \emptyset,$ then 

$$\fb(\phi)(-\lambda_1\Delta_{T_1},\dots, -\lambda_k\Delta_{T_k})=<T_{(\lambda)},\phi>.$$

In particular if $F(\zeta) =e^{-\nu \zeta_j}$ for some $\nu \in \C$ such that  $\nu \lambda_j\in \cup_{(\gamma_j,\delta_j)\in M_{a_j,b_j}}\overline S_{\gamma_j,\delta_j}$ then $F(-\lambda_1\Delta_{T_1},\dots, -\lambda_k\Delta_{T_k}=T_j(\nu \lambda_j).$

(iv) If $(\alpha, \beta) \in \mathcal N_{\lambda},$ if $\phi \in \fab,$ and if $N(T,\lambda, \alpha, \beta)\cap Dom(\fb(\phi))\neq \emptyset,$ then 

$$<T_{(\lambda)},\phi>u=\lim_{\stackrel{\epsilon \to (0,\dots,0)}{_{\epsilon \in S^*_{\alpha, \beta}}}}\fb(\phi)(-\lambda_1\Delta_{T_1} +\epsilon_1I, \dots, -\lambda_k\Delta_{T_k}+\epsilon_kI)u \ \ (u\in \mb).$$

(v) If $U \in \mathcal W_{T,\lambda},$ and if $F\in H^{\infty}(U)$ is strongly outer on $U,$ then there exists $u\in \Omega(\mb)\cap Dom(F(-\lambda_1\Delta_{T_1}, \dots, -\lambda_k\Delta_{T_k}))$ such that $F(-\lambda_1\Delta_{T_1}, \dots, -\lambda_k\Delta_{T_k})u \in \Omega(\mb).$

(vi) For every $U\in \mathcal W_{T,\lambda}$ and every $F\in \mathcal S(U)$ there exists a unique $R_F \in \qm(\mb)=\qm(\ma)$  satisfying $R_FG(-\lambda_1T_1, \dots, -\lambda_kT_k)=(FG)(-\lambda_1T_1, \dots, -\lambda_kT_k) $ for every $G \in H^{\infty}(U)$ such that $FG\in H^{\infty}(U).$ The definition of $R_F$ does not depend on the choice of $U,$ and if we set $F(-\lambda_1T_1,\dots, -\lambda_kT_k)=R_F$ the definition of $F(-\lambda_1\Delta_{T_1}, \dots, -\lambda_k\Delta_{T_k})$ agrees with the definition given in (ii) if $F \in \cup_{U\in \mathcal W_{T,\lambda}}H^{\infty}(U),$ the map $F \to F(-\lambda_1T_1,\dots, -\lambda_kT_k)$ is a bounded homomorphism from $\cup_{U \in \mathcal W_{T,\lambda}}\mathcal S(U)$ into $\qm(\mb)=\qm(\ma),$ and we have, for $\chi \in \widehat \ma,$

$$\tilde \chi(F(-\lambda_1T_1, \dots, -\lambda_kT_k))=F(-\lambda_1\tilde \chi(\Delta_{T_1}),\dots, -\lambda_k\tilde \chi(\Delta_{T_k})) \ \ (F\in \cup_{U\in \mathcal W_{T,\lambda}}\mathcal S(U)),$$

where $\tilde \chi$ is the character on $\qm(\ma)$ such that $\tilde \chi_{|_{\ma}}=\chi.$

(vii) If $F(\zeta_1,\dots, \zeta_k)=-\zeta_j$ then $F(-\lambda_1\Delta_1,\dots,-\lambda_k\Delta_k)= \lambda_j\Delta_{T_j}.$
\end{thm}

Proof: In the following we will use the notations $d\zeta =d\zeta_1\dots d\zeta_k,$ $\lambda \Delta _T=(\lambda_1\Delta_{T_1},\dots, \lambda_k\Delta_{T_k }),$  $R(-\lambda \Delta_T, \zeta)=(-1)^k
(\lambda_1 \Delta_{T_1} +\zeta_1 I)^{-1}\dots (\lambda_k \Delta_{T_k} +\zeta_k I)^{-1})$  for $\zeta=(\zeta_1,\dots, \zeta_k) \in -Res_{ar}(\lambda\Delta _T):=-\Pi_{j=1}^kRes_{ar}(\Delta_{T_j(\lambda_j.)})$ With these notations, the formula given in (i) takes the form

$$F(-\lambda \Delta_T)={(-1)^k\over (2i\pi)^k}\int_{\tilde \partial U +\epsilon}F(\zeta)R(-\lambda \Delta_T,\zeta)d\zeta.$$ Clearly, $F(-\lambda \Delta_T)\in \mathcal M(\mb)\subset \qm_r(\ma).$ Let $U, U' \in \mathcal W_{T,\lambda},$ let $(\alpha, \beta)$ and $(\alpha', \beta')$ be the elements of $M_{a,b}$ associated to $U$ and $U'$ and let $\epsilon \in S^*_{\alpha, \beta}$ and $\epsilon' \in S^*_{\alpha', \beta'}$ such that $U+\epsilon \in \mathcal W_{T,\lambda}$ and $U'+\epsilon'\in \mathcal W_{T,\lambda}.$

Set $V=U+\epsilon, V'=U'+\epsilon', V''=V\cap V'.$ Then the function $G: \zeta \to F(\zeta)R(-\lambda\Delta_T,\zeta)$ is holomorphic on a neighborhood of $\overline V \setminus V'',$ 
and it follows from (43) that  there exists $M>0$ such that $\vert G(\zeta)\vert \le M$ for $\zeta \in \overline V\setminus V''.$
The open sets $V=\Pi_{j\le k}V_j$ and $V''=\Pi_{j\le k}V''_j$  have the form $(z +S^*_{\alpha,\beta})\setminus K)$ and $(z'' +S^*_{\alpha'',\beta''})\setminus K'')$ where $K$ and $K''$ are compact subsets of $\C^k,$ and where $\alpha''=\inf(\alpha, \alpha')$ and $\beta''=sup(\beta, \beta').$ Choose $\epsilon '' \in S^*_{\alpha'',\beta''},$ and denote by $V_{L,k}$ the intersection of $V_k\setminus \overline V''_k$ wtih the strip having for boundaries the lines $D^1_{L}=Le^{i(-{\pi\over 2} -\alpha_k)} +\R \epsilon''$ and  $D^2_L=Le^{i({\pi\over 2} +\beta_k)}+\R \epsilon".$ Set $W_{n,j}(\zeta_j)={n^2\over \left (n + 1  +(\zeta_j  - z_j)e^{i{\alpha_j +\beta_j\over 2}}\right )^2},$ and set $W_n(\zeta)=W_{n,1}(\zeta_1)\dots W_{n,k}(\zeta_k).$ Then $\vert W_n(\zeta)\vert < 1$ and $\lim_{n\to +\infty}W_n(\zeta)=1$ for $\zeta \in \overline V.$

It follows from Cauchy's theorem that we have, when $L$ is sufficiently large

$$0 =\int_{\Pi_{j\le k-1}\partial V_j} \left [ \int_{\partial V_{L,k}} W_n(\zeta)G(\zeta)d\zeta_k\right ]d\zeta_1\dots d\zeta_{k-1}.$$

We have, for $s=1,2,$

$$\left \Vert \int _{(\Pi_{j\le k-1}\partial V_j)\times (\partial V_k\cap D^s_L)}W_n(\zeta)G(\zeta)d\zeta\right \Vert \le M \left [\Pi_{j\le k-1}\int _{\partial V_j}\vert W_{n,j}(\zeta_j)\vert \vert d\zeta_j\vert \right ]\int_{\partial V_k\cap D^s_L}\vert W_{n,k}\vert \vert d\zeta_k\vert,$$

and so $\lim_{L\to +\infty}\int _{(\Pi_{j\le k-1}\partial V_j)\times (\partial V_k\cap D^s_L)}W_n(\zeta)G(\zeta)d\zeta=0.$ We obtain

$$\int_{\tilde \partial V}W_n(\zeta)G(\zeta)d\zeta=\int_{{\Pi_{j\le k-1}\partial V_j\times \partial V''_k}}W_n(\zeta)G(\zeta)d\zeta.$$

It follows then from the Lebesgue dominated convergence theorem that we have

$$\int_{\tilde \partial V}G(\zeta)d\zeta=\int_{{\Pi_{j\le k-1}\partial V_j\times \partial V''_k}}G(\zeta)d\zeta.$$

Using the same argument and a finite induction, we obtain

$$\int _{\tilde \partial V}G(\zeta)d\zeta=\int_{\tilde \partial V''}G(\zeta)d\zeta.$$

Similarly $\int _{\tilde \partial V'}G(\zeta)d\zeta=\int_{\tilde \partial V''}G(\zeta)d\zeta,$ which shows that the definition of $F(-\lambda_1\Delta_{T_1},\dots, -\lambda_k\Delta_{T_k})$ does not depend on the choice of $U$ and $\epsilon.$

Now let $F \in \cup_{U\in \mathcal W_{T,\lambda}}H^{(1)}(U),$ let $G \in \cup_{U\in \mathcal W_{T,\lambda}}H^{(1)}(U).$ There exists $U\in \mathcal W_{T,\lambda}$ such that $F_{|_U}\in H^{(1)}(U)$ and $G_{|_U}\in H^{(1)}(U).$ Choose $\epsilon \in S^*_{\alpha, \beta},$ where $(\alpha, \beta)$ is the element of  $M_{a,b}$ associated to $U,$ such that $U+\epsilon \in \mathcal W_{T,\lambda},$ and set $V=U+{\epsilon \over 2},$ $V'=U+\epsilon.$ For $M\subset\{1,\dots,k\},$ denote by $\vert M\vert$ the cardinal of $M.$ Then $\vert \{1,\dots,k\} \setminus M\vert =2^k-\vert M \vert.$ Since $(\lambda_j\Delta_{T_j}+\zeta_j I)^{-1}(\lambda_j\Delta_{T_j}+\sigma_j I)^{-1}={1\over \sigma_j-\zeta_j}\left ( (\lambda_j\Delta_{T_j}+\zeta_j I)^{-1}-(\lambda_j\Delta_{T_j}+\sigma_j I)^{-1}\right ),$ we have

$$F(-\lambda \Delta_T)G(-\lambda\Delta_T)={1\over (2i\pi)^{2k}}\int_{{\tilde \partial V}\times{\tilde \partial V'}}F(\zeta)G(\sigma)R(-\lambda\Delta_T, \zeta)R(-\lambda \Delta_T, \sigma)d\zeta d\sigma$$

$$={1\over (2i\pi)^{2k}}\sum_{M \subset \{1,\dots,k\}}L_M,$$ 

where

{\small $$L_M:=(-1)^{\vert M\vert}\int _{\tilde \partial V \times \tilde \partial V'}{1\over \sigma_1-\zeta_1}\dots{1\over \sigma_k-\zeta_k}F(\zeta)G(\sigma)\Pi_{j\in M}(\lambda_j\Delta_{T_j}+\zeta_jI)^{-1}\Pi_{j'\notin M}(\lambda_{j'}\Delta_{T_{j'}}+\sigma_{j'}I)^{-1}d\zeta d\sigma.$$} 

Assume that $M\neq \emptyset,$ and set $W_{n,M}((\sigma_j)_{j\in M})=\Pi_{j\in M} {n^2\over \left (n+1 +(\sigma_j-z_j)e^{i{\alpha_j+\beta_j\over 2}}\right )^2},$ where $z\in \C^k$ is choosen so that $z +S^*_{\alpha,\beta}\supset U.$ It follows from corollary 12.4 that $G$ is bounded on $V$, and so the function $(\sigma_j)_{j\in M} \to W_{n,M}((\sigma_j)_{j\in M})\Pi_{j\in M}{1\over \sigma_j-\zeta_j}G(\sigma)$ belongs to $H^{(1)}(\Pi_{j\in M}V_j +{1\over4}(\epsilon_j)_{j\in M} )$ for every $(\sigma_j)_{j\notin M}$ and every $\zeta \in \tilde \partial V.$

Since the open set $\Pi_{j\in M}V_j$ is admissible with respect to the family $\{ (\alpha_j, \beta_j)\}_{j\in M},$ it follows from theorem 12.5 that we have, for every $(\sigma_j)_{j\notin M}\partial V'_j$ and every $\zeta \in \tilde \partial V$

$$\int_{\Pi_{j\in M} \partial V_j}W_{n,M}((\sigma_j)_{j\in M})\Pi_{j\in M}{1\over \sigma_j-\zeta_j}G(\sigma)\Pi_{j\in M}d\sigma_j=0.$$

Set $P( \zeta, (\sigma_{j'})_{j'\notin M}):=\int_{\Pi_{j\in M} \partial V_j}\Pi_{j\in M}{1\over \sigma_j-\zeta_j}G(\sigma)\Pi_{j\in M}d\sigma_j.$

It follows then from the Lebesgue dominated convergence theorem  that we have, for $(\sigma_j)_{j\notin M}\in \Pi_{j'\notin M}\partial V'_j$ and $\zeta \in \tilde \partial V$

$$P( \zeta, (\sigma_{j'})_{j'\notin M})=0,$$

and so

$$(-1)^{\vert M\vert}L_M$$ {\footnotesize $$=\int_{\tilde  \partial V \times (\Pi_{j'\notin M}\partial V'_{j'})}\Pi_{j'\notin M}{1\over \sigma_{j'}-\zeta_{j'}}F(\zeta)P( \zeta, (\sigma_{j'})_{j'\notin M})\Pi_{j\in M}(\lambda_j\Delta_{T_j}+\zeta_jI)^{-1}\Pi_{j'\notin M}(\lambda_{j'}\Delta_{T_{j'}}+\sigma_{j'}I)^{-1}d\zeta \Pi_{j'\notin M}d\sigma_{j'}$$} $$=0.$$

We obtain

$$F(-\lambda \Delta_T)G(-\lambda \Delta_T)={1\over (2i\pi)^{2k}}L_{\emptyset}$$ $$={1\over (2i\pi)^{2k}}\int_{\tilde \partial V'}\left [ \int_{\tilde \partial V}{F(\zeta)\over (\sigma_1-\zeta_1)\dots (\sigma_k-\zeta_k)}d\zeta\right ]G(\sigma) (\lambda_1\Delta_1+\sigma_1I)^{-1}(\lambda_k\Delta_k+\sigma_kI)^{-1}d\sigma$$ $$={1\over (2i\pi)^k}\int_{\tilde \partial V'}F(\sigma)G(\sigma)(\lambda_1\Delta_1+\sigma_1I)^{-1}(\lambda_k\Delta_k+\sigma_kI)^{-1}d\sigma= (FG)(-\lambda \Delta_T),$$

and so the map $F\to F(-\lambda \Delta_{T})$ is an algebra homomorphism from $\cup_{U\in \mathcal W_{T,\lambda}}H^{(1)}(U)$ into $\mathcal M(\mb).$

Let $\mathcal E$ be a bounded subset of $\cup_{U\in \mathcal W_{T,\lambda}}H^{(1)}(U),$ let $U\in W_{T,\lambda}$ such that $\mathcal E$ is a bounded subset of $H^{(1)}(U),$ let $(\alpha,\beta)$ be the element of $M_{a,b}$ associated to $U,$  and let $\epsilon \in S^*_{\alpha, \beta}$ be such that $U+\epsilon \in \mathcal W_{T,\lambda}.$ Set $K =\sup_{\zeta \in \tilde \partial U +\epsilon}\Vert R(-\lambda \Delta_{T}, \zeta)\Vert_{\mathcal M(\mb)}.$ We have

$$\sup_{F \in \mathcal E}\Vert F(-\lambda\Delta_{T})\Vert_{\mathcal M(\mb)}\le {K\over (2\pi)^k}\sup_{F\in \mathcal E}\Vert F \Vert_{H^{(1)}(U)}<+\infty,$$

which shows that the map $F\to F(-\lambda\Delta_{T})$ is a bounded homomorphism from $\cup_{U\in \mathcal W_{T,\lambda}}H^{(1)}(U)$ into $\mathcal M(\mb)\subset \qm_r(\ma).$

(ii) Let $U \in \mathcal W_{T,\lambda},$ and let $(\alpha, \beta)\in M_{a,b}$ and $z \in \C^k$ be such that $U \subset  z+S^*_{\alpha, \beta}$ and $(z+S^*_{\alpha, \beta})\setminus U$ is bounded. For $j\le k,$ set

 $$s_j =1+\sup\left (\lim_{t\to +\infty}{log\left (\left \Vert T\left (t\lambda_je^{i{\alpha_j+\beta_j\over 2}}\right )\right \Vert \right )\over t}, -Re(z_je^{i{\alpha_j+\beta_j\over 2}}\right ).$$ Set $\tilde T_j(t)=T(t\lambda_je^{{\alpha_j+\beta_j\over 2}})$ for $t>0,$ with the convention $\tilde T_j(0)=I,$ and set, for $f \in \cap_{\zeta \in \C^k}e_{-\zeta}\mathcal U_{\alpha, \beta},$

$$<f,\phi>=\int_{(\R^+)^k}f(t_1e^{i{\alpha_1+\beta_1\over 2}},\dots,t_ke^{i{\alpha_k+\beta_k\over 2}})e^{-s_1t_1-\dots-s_kt_k}dt_1\dots dt_k.$$

Then $z \in Dom(\fb(\phi)),$ and we have, for $\zeta \in Dom(\fb(\phi)),$

$$\fb(\phi)(\zeta)=\int_{(R^+)^k}e^{-t_1\zeta_1e^{i{\alpha_1+\beta_1\over 2}}\dots-t_k\zeta_ke^{i{\alpha_k+\beta_k\over 2}}}e^{-s_1t_1+\dots-s_kt_k}dt_1\dots dt_k$$

$$={1\over (\zeta_1e^{i{\alpha_1+\beta_1\over 2}}+s_1)\dots(\zeta_ke^{i{\alpha_k+\beta_k\over 2}}+s_k)},$$

$$<T_{(\lambda)}, \phi>=\int_{(\R^+)^k}T_{1}(t_1\lambda_1e^{i{\alpha_1+\beta_1\over 2}})\dots T_{k}(t_1\lambda_ke^{i{\alpha_k+\beta_k\over 2}})e^{-s_1t_1-\dots-s_kt_k}dt_1\dots dt_k$$

$$=\left [ \int_0^{+\infty}\tilde T_1(t_1)e^{-s_1t_1}dt_1\right ] \dots \left [ \int_0^{+\infty}\tilde T_k(t_1)e^{-s_kt_k}dt_k\right ],$$

where the Bochner integrals are computed with respect to the strong operator topology on $\mathcal M(\mb).$

It follows from the observations in section 5 that $\left ( \int_0^{+\infty}\tilde T_j(t_1)e^{-s_jt_1}dt_j\right )(\mb)$ is dense in $\mb$ for $1\le j \le k,$ and so $<T_{(\lambda)},\phi>(\mb)$ is dense in $\mb.$ Now set $\phi_1=\phi*\phi.$ It follows from theorem 8.6 that we have

$$<T_{(\lambda)},\phi_1>=<T_{\lambda}, \phi>^2,$$

and so $<T_{(\lambda)},\phi_1>(\mb)$ is dense in $\mb.$

Set $F=\fb(\phi_1)=\fb(\phi)^2.$ Then $F\in H^{(1)}(U-\epsilon)\cap H^{\infty}(U-\epsilon)$ for some $\epsilon \in S^*_{\alpha,\beta},$ and we have, using assertion (vi) of theorem 8.6

$$F(-\lambda \Delta_T)={(-1)^k\over (2i\pi)^k} \int_{z+\tilde \partial S^*_{\alpha, \beta}}F(\zeta)(\lambda_1\Delta_{T_1}+\sigma_1)^{-1}\dots (\lambda_k\Delta_{T_k}+\sigma_k)^{-1}d\sigma_1\dots d\sigma_k$$ $$=<T_{(\lambda)}, \phi_1>,$$

which shows that $F(-\lambda \Delta_T)(\mb)$ is dense in $\mb.$

Now consider again $U\in \mathcal W_{T,\lambda},$ and let $F\in H^{\infty}(U).$ Let $G_0\in H^{(1)}(U)$ be such that $G_0(-\lambda \Delta_T)(\mb)$ is dense in $\mb,$ and let $u\in \Omega(\mb).$
Then $G_0(-\lambda \Delta_T)u\in \Omega(\mb),$ $FG_0\in H^{(1)}(U),$ and so there exists a unique quasimultiplier $R_F\in \qm_r(\mb)=\qm_r(\ma)$ such that $R_FG_0(-\lambda \Delta_T)u=(FG_0)(-\lambda\Delta_T)u,$ and $R_F=F(-\lambda \Delta_T)$ if $F\in H^{(1)}(U).$

Let $U'\in \mathcal W_{T,\lambda},$ and let $G\in H^{(1)}(U').$ We have

$$R_FG(-\lambda \Delta_T)G_0(-\lambda \Delta_T)=R_FG_0(-\lambda \Delta_T)G(-\lambda \Delta_T)=(FG_0)(-\lambda\Delta_T)G(-\lambda \Delta_T)$$
$$=(FG_0G)(-\lambda \Delta_T)=(FG)(-\lambda \Delta_T)G_0(-\lambda\Delta_T),$$
and so $R_FG(-\lambda \Delta_T)=(FG)(-\lambda \Delta_T),$ which shows that the definition of $R_F$ does not depend on the choice of $U.$ The map $F\to R_F$ is clearly linear.
Now let $F_1\in \cup_{U\in  \mathcal W_{T,\lambda}}H^{\infty}(U),$ let $F_2\in \cup_{U\in  \mathcal W_{T,\lambda}}H^{\infty}(U),$ and let $G\in \cup_{U\in  \mathcal W_{T,\lambda}}H^{(1)}(U)$ such that $G(-\lambda \Delta_T)(\mb)$ is dense in $\mb.$ We have

$$R_{F_1F_2}G^2(-\lambda\Delta_T)=(F_1F_2G^2)(-\lambda\Delta_T)=(F_1G)(\llt)(F_2G)(\llt)$$ $$=R_{F_1}R_{F_2}G^2(\llt),$$

and so $R_{F_1F_2}=R_{F_1}R_{F_2}$ since $G^2(\llt)\mb$ is dense in $\mb.$

Now let $\mathcal E$ be a bounded family of elements of $\cup_{U\in \mathcal W_{T,\lambda}}H^{\infty} (U).$ There exists $U \in \mathcal W_{T,\lambda}$ and $M>0$ such that
$F\in H^{(\infty)}(U)$ and $\Vert F\Vert_{H^{\infty}(U)}\le M$ for every $F \in \mathcal E.$ Let $G\in H^{(1)}(U)$ such that $G(-\lambda \Delta_T)(\mb)$ is dense in $\mb.$ Then the family $\{FG\}_{F\in \mathcal E}$ is bounded in $H^{(1)}(U),$ and it follows from (i) that there exists $u\in \Omega(\mb)$ such that $\sup_{F\in \mathcal E}\Vert (FG)(\llt)u\Vert_{\mb}<+\infty.$

We obtain

$$\sup_{F\in \mathcal E}\Vert R_FG(\llt)u\Vert_{\mb}=\sup_{F\in \mathcal E}\Vert (FG)(\llt)u\Vert_{\mb}<+\infty,$$

and so the family $\{R_F\}_{F\in \mathcal E}$ is pseudobounded in $\qm(\mb)=\qm(\ma)$ since $G(\llt)u\in \Omega(\mb).$ Since the family $\{\lambda^{-n}F^n\}_{n\ge 1} $ is bounded in 
$H^{\infty}(U)$ for $F\in H^{\infty}(U),$ $\lambda > (1+\Vert F\Vert_{H^{\infty}(U)})^{-1},$ this shows that $R_F\in \qm_r(\mb)=\qm_r(\ma)$ for $F\in \cup_{U\in \mathcal F}H^{\infty}(U),$ and that the map $F\to R_F$ is a bounded algebra homomorphism from $\cup_{U\in \mathcal F}H^{\infty}(U)$ into $\qm_r(\mb)=\qm_r(\ma),$ which concludes the proof of (ii).

(iii) Let $(\alpha, \beta) \in \mathcal N_{\lambda},$ let $\phi \in \fab,$ assume that $N_0(T,\lambda, \alpha, \beta)\cap Dom(\fb(\phi))\neq \emptyset,$ and let $z \in N_0(T,\lambda, \alpha, \beta)\cap Dom(\fb(\phi)).$ Then $z+S^*_{\alpha, \beta}$ is admissible with respect to $(T,\lambda, \alpha, \beta).$ As in the proof of (ii) we can construct $\phi_1\in \mathcal F_{\alpha,\beta}$ having the following properties
\begin{itemize}

 \item $z \in N_0(T,\lambda,\alpha, \beta)\cap Dom(\fb(\phi_1))$,
 
 \item $G:=\fb(\phi_1)\in H^{(1)}(z+S^*_{\alpha, \beta})\cap H^{\infty}(z+S^*_{\alpha, \beta}),$
 
 \item $<T_{(\lambda)}, \phi_1>=G(-\lambda \Delta_T),$ and $G(-\lambda \Delta_T)(\mb)$ is dense in $\mb.$
 
 \end{itemize}
 
 Let $\epsilon \in S^*_{\alpha, \beta}$ be such that $z+\epsilon +S^*_{\alpha, \beta}$ is admissible with respect to $(T,\lambda, \alpha, \beta).$ It follows from assertions (v) and (vi) of theorem 8.6 and from (i) and (ii) that we have
 
 $$<T_{(\lambda)},\phi>\fb(\phi_1)(-\lambda \Delta_T)=<T_{\lambda)},\phi><T_{(\lambda)},\phi_1>=<T_{(\lambda)},\phi*\phi_1>$$ $$=\int_{z+\epsilon+S^*_{\alpha, \beta}}\fb(\phi)(\sigma)\fb(\phi_1)(\sigma)(\lambda_1\Delta_{T_1}+ \sigma_1I)^{-1}\dots (\lambda_k\Delta_{T_k}+ \sigma_kI)^{-1}d\sigma_1\dots d\sigma_k$$ $$=(\fb(\phi)\fb(\phi_1))(-\lambda \Delta_T)= \fb(\phi)(-\lambda\Delta_T)\fb(\phi_1)(-\lambda \Delta_T),$$
 
 and so $<T_{(\lambda)},\phi>=\fb(\phi)(-\lambda\Delta_T)$ since $\fb(\phi_1)(-\lambda \Delta_T)(\mb)$ is dense in $\mb.$
 
 Now let $\nu \in \C$ such that $\nu \lambda_j\in  \cup_{(\gamma_j,\delta_j)\in M_{a_j,b_j}}\overline S_{\gamma_j,\delta_j},$ and let  $\nu_j=(\nu_{j,1},\dots, \nu_{j,k})$ be the $k$-tuple defined by the conditions $\nu_{j,s}=0$ if $s \neq j,$ $\nu_{j,j}=\nu.$ There exists $(\gamma_j,\delta_j)\in \mathcal N_{\lambda_j}$ such that $\nu \in \overline S_{\gamma_j,\delta_j},$ and there exists $(\alpha, \beta)\in {\mathcal N}_{\lambda}$ such that $\alpha_j=\gamma_j$ and $\beta_j=\delta_j.$
 Set $F(\zeta)=e^{-\nu \zeta_j}$ for $\zeta \in \C^k,$ and set $<f,\phi>=f(\nu_j)$ for $f\in \cap_{z\in \C^k}e_{-z}\mathcal U_{\alpha,\beta}.$ Then $Dom(\fb(\phi))=\C^k,$ and we have,  for $\zeta \in \C^k,$

 $$\fb(\phi)(\zeta)=<e_{-\zeta},\phi>=e^{-\nu_j\zeta}=e^{-\nu \zeta_j},$$
 
 and so $F=\fb(\phi).$ Let $z \in N_0(T,\lambda, \alpha, \beta)=N_0(T,\lambda, \alpha, \beta)\cap Dom(\fb(\phi)).$ Let $\delta_{\nu_j}$ be the Dirac measure at $\nu_j.$ Since $e_{-z}\delta_{\nu_j}$ is a representing measure for $\phi e_{-z}$ we have
 
 $$F(-\lambda \Delta_T)=<T_{(\lambda)},\phi>=T_j(\nu \lambda_j),$$
 
 which concludes the proof of (iii).

 (iv)  Let $(\alpha, \beta) \in \mathcal N_{\lambda},$ let $\phi \in \fab,$ and assume that $N(T,\lambda, \alpha, \beta)\cap Dom(\fb(\phi))\neq \emptyset.$ Set $e_{-\epsilon}T=(e_{-\epsilon_1}T_1,\dots,e_{-\epsilon_k}T_k).$ Then $N(T,\lambda,\alpha, \beta)\subset N_0(e_{-\epsilon}T, \lambda, \alpha, \beta)$ for $\epsilon \in S^*_{\alpha,\beta},$ and it follows from theorem 8.6 (ii) and from (iii) that we have, for $u\in \mb,$
 
 $$<T_{(\lambda)},\phi>u=\lim_{\stackrel{\epsilon \to (0,\dots,0)}{_{\epsilon \in S^*_{\alpha,\beta}}}}<e_{-\epsilon }T_{\lambda},\phi>u$$ $$=\lim_{\stackrel{\epsilon \to (0,\dots,0)}{_{\epsilon \in S^*_{\alpha,\beta}}}}\fb(\phi)(e_{-\epsilon}T_{(\lambda)})u=\lim_{\stackrel{\epsilon \to (0,\dots,0)}{_{\epsilon \in S^*_{\alpha,\beta}}}}\fb(\phi)(-\lambda_1T_1+\epsilon_1I,\dots,-\lambda_kT_k+\epsilon_kI)u,$$
 
 which concludes the proof of (iv).
 
 (v) Let $U\in \mathcal W_{t,\lambda},$ let $F\in H^{\infty}(U)$ be strongly outer, and let $(F_n)_{n\ge 1}$ be a sequence of invertible elements of $H^{\infty}(U)$ satisfying the conditions of definition 12.6 with respect to $F.$  It follows from (ii) that there exists $G\in H^{(1)}(U)\cap H^{\infty}(U)$ such that $G(-\lambda \Delta_T)(\mb)$ is dense in $\mb.$ Let $(\alpha,\beta) \in M_{a,b}$ and $z \in N_0(T,\lambda, \alpha, \beta)$ such that $U\subset z+S^*_{\alpha, \beta}$ and such that $(z+S^*_{\alpha, \beta})\setminus U$ is bounded. There exists $\epsilon \in \C^k$ such that $U+\epsilon \subset U$ is admissible with respect to $(T,\lambda, \alpha, \beta)$ and we have
 
 $$F(-\lambda \Delta_T)F_n^{-1}(-\lambda \Delta_T)G^2(-\lambda \Delta_T)$$ $$={1\over (2i\pi)^k}\int _{\epsilon +\tilde \partial U}F(\sigma)F^{-1}_n(\sigma)G^2(\sigma)(\lambda_1\Delta_{T_1} +\sigma_1I)^{-1}\dots(\lambda_k\Delta_{T_k}+\sigma_kI)^{-1}d\sigma_1\dots d\sigma_k,$$
 
 and it follows from  the Lebesgue dominated convergence theorem that $$\lim_{n\to +\infty}\Vert F(-\lambda \Delta_T)F_n^{-1}(-\lambda \Delta_T)G^2(-\lambda \Delta_T)-G^2(-\lambda \Delta_T)\Vert_{\mathcal M(\mb)}=0.$$
 
 Let $u \in \Omega(\mb).$ Then $G(-\lambda \Delta_T)u \in Dom(F(-\lambda \Delta_T))\cap\Omega(\mb).$ Set $u_n=F^{-1}_n(-\lambda \Delta_T)G(-\lambda \Delta_T)u\in \mb.$  We have $$G(-\lambda \Delta_T)^2u^2=\lim_{n\to +\infty}F(-\lambda \Delta_T)
 G(-\lambda \Delta_T)uu_n.$$  Since $G(-\lambda \Delta)^2u^2\in \Omega(\mb),$ this shows that $F(-\lambda \Delta_T)
 G(-\lambda \Delta_T)u \in \Omega(\mb),$ which proves (v).
 
 (vi) Let $U\in \mathcal W_{T,\lambda},$ let $F\in \mathcal S(U),$ let $G_0\in H^{\infty}(U)$ be a strongly outer function such that $FG_0\in H^{\infty}(U),$ and let $u\in Dom(G_0(-\lambda\Delta_T))$
 such that $G_0(-\lambda\Delta_{T})u\in \Omega(\mb).$ Let $v \in \Omega(\mb)\cap Dom(FG_0(-\lambda\Delta_T)).$ There exists a unique $R_F \in \qm(\mb)=\qm(\ma)$ satisfying the equation
 
 $$(FG_0)(-\lambda \Delta_T)uv=R_FG_0(-\lambda \Delta_T)uv,$$
 
 and we have $$(FG_0)(-\lambda \Delta_T)=R_FG_0(-\lambda \Delta_T),$$
 
 so that $R_F=F(-\lambda \Delta_T)$ if $F\in H^{\infty}(U).$
 
 Let $G\in \cup_{V\in \mathcal W_{T,\lambda}}H^{\infty}(V)$ such that $FG \in H^{\infty}(W)$ for some $W \in \mathcal W_{T,\lambda},$ and let $w\in \Omega(\mb)\cap Dom(G(-\lambda \Delta_T)).$ We have
 
 $$((FG)(-\lambda \Delta_T)vw)G_0(-\lambda \Delta_T)u=(FG_0)(-\lambda \Delta_T)G(-\lambda \Delta_T)uvw$$ $$=R_FG_0(-\lambda \Delta_T)G(-\lambda \Delta_T)uvw
 =(R_FG(-\lambda \Delta_T)vw)G_0(-\lambda\Delta_T)u.$$

 Since $vw\left (G_0(-\lambda\Delta_T)u\right)\in \Omega(\mb),$ this shows that $(FG)(-\lambda \Delta_T)=R_FG(-\lambda \Delta_T).$ So if we set $F(-\lambda \Delta_T)=R_F,$
 we obtain $F(-\lambda \Delta_T)G(-\lambda\Delta_T)=(FG)(-\lambda \Delta_T)$ for every $F\in \cup_{U\in \mathcal W_{T,\lambda}}\mathcal S(U)$ and for every $G\in \cup_{U\in \mathcal W_{T,\lambda}}H^{\infty}(U)$ such that $FG\in \cup_{U\in \mathcal W_{T,\lambda}}H^{\infty}(U).$ The map $F\to F(-\lambda \Delta_T)$ is clearly linear. Now let $F_1\in\cup
 _{U\in \mathcal W_{T,\lambda}}\mathcal S(U), F_2 \in \cup_{U\in \mathcal W_{T,\lambda}}\mathcal S(U),$ and let $G_1\in \cup_{U\in \mathcal W_{T,\lambda}}H^{\infty}(U)$ and $G_2\in \cup_{U\in \mathcal W_{T,\lambda}}H^{\infty}(U)$ be strongly outer functions such that $F_1G_1\in \cup_{U\in \mathcal W_{T,\lambda}}H^{\infty}(U)$ and $F_2G_2\in \cup_{U\in \mathcal W_{T,\lambda}}H^{\infty}(U).$ We have
 
 $$(F_1F_2)(-\lambda \Delta_T)G_1(-\lambda \Delta_T)G_2(-\lambda \Delta_T)=(F_1F_2G_1G_2)(-\lambda\Delta_T)$$ $$=(F_1G_1)(-\lambda\Delta_T)(F_2G_2)(-\lambda\Delta_T)
 =F_1(-\lambda\Delta_T)F_2(-\lambda\Delta_T)G_1(-\lambda \Delta_T)G_2(-\lambda \Delta_T),$$
 
 and so $(F_1F_2)(-\lambda \Delta_T)=F_1(-\lambda\Delta_T)F_2(-\lambda\Delta_T)$ since $ Dom(G_1(-\lambda \Delta_T))\cap \Omega(\mb)\neq \emptyset$ and $ Dom(G_2(-\lambda \Delta_T))\cap \Omega(\mb)\neq \emptyset,$ and the map $F \to F(-\lambda \Delta_T)$ is an algebra homomorphism from $\cup_{U\in \mathcal W_{T,\lambda}}\mathcal S(U)$
 into $\qm(\mb)=\qm(\ma).$

 Now let $\mathcal E$ be a bounded family of elements of $\cup_{U\in \mathcal W_{T,\lambda}}\mathcal S(U).$ There exists $U \in \mathcal W_{T,\lambda}$ and a strongly outer function $G\in H^{\infty}(U)$ such that
$FG\in H^{\infty}(U)$ for every $F \in \mathcal E$ and such that $\sup_{F\in \mathcal E}\Vert FG\Vert_{H^{\infty}(U)}<+\infty.$ So the family $\{ (FG)(-\lambda\Delta_T)\}_{F\in \mathcal E}$ is a pseudobounded family of elements of $\qm_r(\mb)=\qm_r(\ma),$ and there exists $u \in \Omega(\mb)\cap (\cap_{F\in \mathcal E}Dom((FG)(-\lambda \Delta_T)))$ such that $\sup_{F\in \mathcal E}\Vert (FG)(-\lambda \Delta_T)u\Vert_{\mathcal B}<+\infty.$ Let $v\in Dom(G(-\lambda \Delta_T))\cap\Omega(\mb)$, and set $w=G(-\lambda \Delta_T)uv.$ Then $w \in \Omega(\mb) \cap(\cap_{F\in \mathcal E}Dom(F(-\lambda \Delta_T))$ and

$$\sup_{F\in \mathcal E}\Vert F(-\lambda\Delta_T)w\Vert_{\mb}=\sup_{F\in \mathcal E}\Vert( F(-\lambda\Delta_T)G(-\lambda\Delta_T)uv\Vert_{\mb}$$ $$\le \sup_{F\in \mathcal E}\Vert (FG)(-\lambda\Delta_T)u\Vert_{\mb}\Vert v\Vert_{\mb}<+\infty,$$
and so the family $\{F(-\lambda\Delta_T\}_{F\in \mathcal E}$ is pseudobounded in $\qm(\mb)=\qm(\ma),$ and  the map $F\to F(-\lambda \Delta_T)$ is a bounded algebra homomorphism from $\cup_{U\in \mathcal W_{T,\lambda}}\mathcal S(U)$ into $\qm\mb)=\qm(\ma),$ which concludes the proof of (vi).

Now assume that $\ma$ is not radical, let $\chi \in \widehat \ma,$ and let $\tilde \chi$ be the unique character on $\qm(\ma)$ such that $\tilde \chi(u)=\chi(u)$ for every $u \in \ma.$

Let $F\in H^{(1)}(U),$ where $U\in \mathcal W_{T, \lambda},$ let $(\alpha, \beta)$ be the element of $M_{a,b}$ associated to $U,$ and let $\epsilon \in S^*_{\alpha, \beta}$ be such that $U+\epsilon$ is admissible with respect to $(T, \lambda, \alpha, \beta).$
Since Bochner integrals commute with linear functionals, we have

$$\tilde \chi \left (F(-\lambda_1\Delta_{T_1},\dots,-\lambda_k\Delta_{T_k})\right )$$ $$={1\over (2i\pi)^k}\int_{\tilde \partial U +\epsilon}F(\zeta_1,\dots,\zeta_k)(\lambda_1\tilde \chi(\Delta_{T_1} )+\zeta_1 I)^{-1}\dots (\lambda_1\tilde \chi(\Delta_{T_k}) +\zeta_k I)^{-1}d\zeta_1\dots d\zeta_k.$$

Since $U+\epsilon$ is admissible with respect to $(T,\lambda, \alpha, \beta),$ $(-\lambda_1\tilde \chi(\Delta_{T_1}),\dots,-\lambda_k\tilde \chi(\Delta_{T_k}))\in U+\epsilon,$ and it follows from theorem 12.5 that we have

$$\tilde \chi \left (F(-\lambda_1\Delta_{T_1},\dots,-\lambda_k\Delta_{T_k})\right )=F(-\lambda_1\tilde \chi (\Delta_{T_1}),\dots,-\lambda_k\tilde \chi (\Delta_{T_k})).$$

Now let $F\in H^{\infty}(U),$ where $U\in \mathcal W_{T,\lambda},$ and let $G\in H^{(1)}(U)$ such that $G(-\lambda \Delta_T)(\mb)$ is dense in $\mb$. Then $\tilde \chi(G(-\lambda \Delta_T))\neq 0, $ and we have

$$\tilde \chi (F(-\lambda \Delta_T))= {\tilde \chi ((FG)(-\lambda \Delta_T))\over \tilde \chi (G(-\lambda\Delta_T))}={(FG)(-\lambda \tilde \chi (\Delta_T))\over G(-\lambda \tilde \chi(\Delta_T))}=F(-\lambda \tilde \chi(\Delta_T)).$$

Finally let $F \in {\mathcal S}(U),$ where $U \in \mathcal W_{T,\lambda},$ and let $G \in H^{\infty}(U)$ be a strongly outer function such that $FG\in H^{\infty}(U).$ It follows from (v) that $G(-\lambda \Delta_{T})u \in \Omega(\mb)$ for some $u \in \mb,$ and so $\tilde \chi (G(-\lambda \Delta_T))\neq 0.$ The same argument as above shows then that $\tilde \chi (F(-\lambda \Delta_T))=F(-\lambda \tilde \chi(\Delta_T)),$ which concludes the proof of (vi).

(vii) Set $F(\zeta_1,\dots, \zeta_k)=-\zeta_j,$ choose $\nu_0 > \nu_1>\lim_{t\to +\infty}{log\Vert T_j(t\lambda_j)\Vert \over t},$ and set  again $v_{\nu_0}(t)=te^{-\nu_0 t}.$ It follows from proposition 12.8(ii) that $F\in \mathcal S(U)$ for every $U\in \mathcal W_{T,\lambda},$ and it follows from proposition 5.5(i) that we have

$$\lambda_j\Delta_{T_j}\int_{[0,\infty)^n}v_{\nu_0}(t)T_j(t\lambda_j)dt=-\int_0^{+\infty}v'_{\nu_0}(t)T_j(t\lambda_j)dt,$$

where the Bochner integrals are computed with respect to the strong operator topology on ${\mathcal M}(\mb).$

Now choose $(\alpha, \beta)\in \mathcal N_{ \lambda},$ and set, for $f \in \cap_{z\in \C^k}e_{-z}\mathcal U_{\alpha, \beta},$

$$<f, \phi_0>=\int_{[0,+\infty)^k}f(0,\dots,0,t_j,0, \dots,0)v_{\nu_0}(t_j)dt_j, $$ $$<f, \phi_1>=\int_{[0,+\infty)^k}f(0,\dots,0,t_j,0, \dots,0)v'_{\nu_0}(t_j)dt_j.$$

Then $\phi_0\in \fab, \phi_1 \in \fab,$ $-\nu_1\lambda_j +S^*_{\alpha, \beta}\in N_0(T, \lambda, \alpha, \beta)\cap Dom(\fb(\phi_0))\cap Dom(\fb(\phi_1)),$ and it follows from (iii) that we have

$$\int_{[0,\infty)^n}v_{\nu_0}(t)T_j(t\lambda_j)dt_j=<T_{(\lambda)},\phi_0>=\fb(\phi_0)(-\lambda\Delta_T), $$ $$\int_{[0,\infty)^n}v'_{\nu_0}(t)T_j(t\lambda_j)dt_j=<T_{(\lambda)},\phi_1>=\fb(\phi_1)(-\lambda\Delta_T).$$

But $\fb(\phi_0)(\zeta)={1\over \nu_0 +\zeta_j}, \fb(\phi_1)(\zeta)={\zeta_j\over \nu_0+\zeta_j}=-F(\zeta)\fb(\phi_0)(\zeta),$ which gives

$$\lambda_j\Delta_{T_j}\int_{[0,\infty)^n}v_{\nu_0}(t)T_j(t\lambda_j)dt=F(-\lambda\Delta_T)\int_{[0,\infty)^n}v_{\nu_0}(t)T_j(t\lambda_j)dt,$$

and so $F(-\lambda\Delta_T)= \lambda_j\Delta_{T_j},$ since $(\int_{[0,\infty)^n}v_{\nu_0}(t)T_j(t\lambda_j)dt)(\mb)$ is dense in $\mb,$ as observed in section 5. $\square$

\section{Appendix 1: Fourier-Borel and Cauchy transforms}

In this section we present some certainly well-known results about Fourier-Borel and Cauchy transforms of linear functionals on some spaces of holomorphic functions on sectors.

For $\alpha< \beta\le \alpha +\pi$ denote as usual by $\overline {S}_{\alpha,\beta}$ the closure of the open sector $S_{\alpha, \beta}$, and set by convention ${\overline S}_{\alpha, \alpha}:=\{te^{i\alpha}\}_{t\ge 0}.$ 

We set

\begin{equation} S^*_{\alpha,\beta}= S_{-\pi/2-\alpha,\pi/2-\beta}, {\overline S}^*_{\alpha,\beta}=\overline {S}_{-\pi/2-\alpha,\pi/2-\beta}.\end{equation}

Notice that $S^*_{\alpha,\alpha +\pi}=\emptyset,$ while ${\overline S}^*_{\alpha,\alpha +\pi}=\overline{S}_{-\pi/2-\alpha,-\pi/2-\alpha}=\{ -tie^{-i\alpha}\}_{t\ge 0}$

Now asssume that $\alpha \le \beta < \alpha +\pi.$ Let $\lambda=\vert \lambda \vert e^{i\omega} \in {\overline S}^*_{\alpha,\beta}$ and let $\zeta=\vert  \zeta \vert e^{i\theta} \in \overline S_{\alpha,\beta},$ with $-{\pi\over 2} -\alpha \le \omega \le {\pi \over 2} - \beta,$ $\alpha \le  \theta \le \beta.$ We have $-{\pi\over 2} \le \omega +\theta \le {\pi\over 2},$ $\vert e^{-\lambda\zeta} \vert =e^{-\vert \lambda\vert \vert \zeta \vert cos(\omega +\theta)},$ and we obtain

\begin{equation}\vert e^{-\lambda\zeta}\vert  < 1\ \ (\lambda \in S^*_{\alpha,\beta}, \zeta \in \overline{S}_{\alpha, \beta}\setminus \{0\}).
\end{equation}
\begin{equation}\vert e^{-\lambda\zeta}\vert  \le 1\ \ (\lambda \in \overline{S}^*_{\alpha, \beta}, \zeta \in \overline{S}_{\alpha, \beta}).
\end{equation}

\begin{defn} Let $\alpha=(\alpha_1,\dots,\alpha_k), \beta=(\beta_1,\dots \beta_k)\in \R^k$ such that $\alpha_j\le \beta_j < \alpha_j+\pi$ for $1\le j \le k.$ Set $\scab:= \Pi_{j=1}^k\overline S_{\alpha_j,\beta_j},$ $\sab^*:= \Pi_{j=1}^k S^*_{\alpha_j,\beta_j},$ $\scab^*:= \Pi_{j=1}^k\overline S^*_{\alpha_j,\beta_j}.$ If, further, $\alpha_j<\beta_j$ for $1\le j \le k,$ set
$\sab:= \Pi_{j=1}^k S_{\alpha_j,\beta_j}.$

Let $X$ be a Banach space. We denote by $\uab(X)$  the set of all continuous $X$-valued functions $f$ on ${\overline S}_{\alpha,\beta}$ satisfying $\lim_{\stackrel{\vert z \vert \to +\infty}{_{z\in \scab}}}\Vert f(z)\Vert_X=0$ such that the map $\zeta \to f(\zeta_1, \zeta_2, \dots, \zeta_{j-1}, \zeta, \zeta_{j+1}, \dots,\zeta_k)$ is holomorphic on $S_{\alpha_j,\beta_j}$ for every $(\zeta_1, \dots, \zeta_{j-1},\zeta_{j+1}, \dots,\zeta_k)\in \Pi_{\stackrel{1\le s \le k}{_{ s \neq j}}}\overline {S}_{\alpha_s, \beta_s}$ when $\alpha_j<\beta_j.$

 Similarly we denote by $\vab(X)$   the set of all continuous bounded $X$-valued  functions $f$ on $\overline \sab$  such that the map $\zeta \to f(\zeta_1, \zeta_2, \dots, \zeta_{j-1}, \zeta, \zeta_{j+1},\dots, \zeta_k)$ is holomorphic on $S_{\alpha_j,\beta_j}$ for every $(\zeta_1, \dots, \zeta_{j-1},\zeta_{j+1}, \dots,\zeta_k)\in \Pi_{\stackrel{1\le s \le k}{_{ s \neq j}}}\overline {S}_{\alpha_s, \beta_s}$ when $\alpha_j<\beta_j.$ The spaces $\uab(X)$ and $\vab(X)$ are equipped with the norm
$ \Vert f \Vert _{\infty}=\sup_{z \in \scab}\Vert f(z)\Vert_X,$ and we will write $\uab:=\uab(\C), \vab:=\vab(\C).$

A representing measure for $\phi  \in \uab'$ is a measure of  bounded variation $\nu$on ${\overline S}_{\alpha, \beta}$ satisfying

\begin{equation} <f,\phi>=\int_{{\overline S}_{\alpha, \beta}}f(\zeta)d\nu(\zeta) \ \ (f\in \uab).\end{equation}.

\end{defn}

Set $I:=\{j\le k \ | \ \alpha_j=\beta_j\}, J:=\{\{j\le k \ | \ \alpha_j<\beta_j\}.$ Since separate holomorphy with respect to each of the variables $z_j, j \in J$ implies holomorphy with respect to $z_J=(z_j)_{j\in J},$ the map $z_J\to f(z_I,z_J)$ is holomorphic on $\Pi_{j\in J}S_{\alpha_j,\beta_j}$ for every $z_I\in \Pi_{j\in I}\overline S_{\alpha_j,\alpha_j}.$

For $z=(z_1,\dots, z_k),\zeta=(\zeta_1,\dots,\zeta_k)\in \C^k,$ set again $ e_z(\zeta)=e^{z_1\zeta_1\dots +z_k\zeta_k}.$ Also set, if $X$ is a separable Banach space, and if $\alpha=(\alpha_1,\dots, \alpha_k)$ and $\beta=(\beta_1,\dots, \beta_k)$ satisfy the conditions above

\begin{equation}\usab(X)={\mathcal U}_{(-\pi/2 -\alpha_1, \dots, -\pi/2-\alpha_k), (\pi/2-\beta_1,\dots, \pi/2-\beta_k)}(X),\end{equation}
\begin{equation}\vsab(X)={\mathcal V}_{(-\pi/2 -\alpha_1, \dots, -\pi/2-\alpha_k), (\pi/2-\beta_1,\dots, \pi/2-\beta_k)}(X),\end{equation}

with the conventions $\usab=\usab(\C),$ $\vsab=\vsab(\C).$

\begin{prop} Let $\phi  \in \uab',$  and let $X$ be a separable Banach space. Set, for $f \in \vab(X),$

$$<f,\phi>=\int_{\scab}f(\zeta)d\nu(\zeta),$$

where $\nu$ is a representing measure for $\phi.$ Then this definition does not depend on the choice of $\nu,$ and we have

\begin{equation}<f,\phi>=\lim_{\stackrel{\epsilon \to 0}{\epsilon \in \overline S^*_{\alpha,\beta}}}<e_{-\epsilon f},\phi>.\end{equation}

\end{prop}

Proof: It follows from (16) and (17) that $ e_{-\epsilon}f\in \uab(X)$ for $f \in \vab(X), \epsilon \in S^*_{\alpha, \beta}.$ If $f \in \uab(X),$ then we have, for $l \in \uab(X)',$ 

$$<\int_{\scab}f(\zeta)d\nu(\zeta), l>=\int_{\scab}<f(\zeta),l>d\nu(\zeta)=<<f(\zeta),l>, \phi>,$$

which shows that the definition of $<f,\phi>$ does not depend on the choice of $\nu.$ Now if $f\in \vab(X),$ it follows from the Lebesgue dominated convergence theorem that we have

$$\int_{\scab}f(\zeta)d\nu(\zeta)=\lim_{\stackrel{\epsilon \to 0}{_{\epsilon \in \sscab}}}\int_{\scab}e_{-\epsilon}(\zeta)f(\zeta)d\nu(\zeta)=\lim_{\stackrel{\epsilon \to 0}{_{\epsilon \in \sscab}}}<e_{-\epsilon}f, \phi>,$$

and we see again that the definition of $<f, \phi>$ does not depend on the choice of the measure $\nu.$ $\square$

We now introduce the classical notions of Cauchy transforms and Fourier-Borel transforms. 




\begin{defn} Let $\phi \in \uab',$ and let $f\in \mathcal V_{\alpha,\beta}(X)$.

(i) The Fourier-Borel transform of $\phi$ is defined  on $\sscab$ by the formula

$$\FB(\phi)(\lambda)=<e_{-\lambda}, \phi> \ \ (\lambda \in \sscab).$$

(ii) The Cauchy transform of $\phi$ is defined on $\Pi_{1\le j \le k}(\C\setminus \overline S_{\alpha_j, \beta_j})$ by the formula

$${\mathcal C}(\phi)(\lambda)={1\over (2i\pi)^k}<{1\over (\zeta-\lambda)}, \phi_{\zeta}>$$ $$:={1\over (2i\pi)^k}<{1\over \zeta_1 -\lambda_1}\dots{1\over \zeta_k-\lambda_k}, \phi_{\zeta_1,\dots, \zeta_k}> $$ $$ (\lambda =(\lambda_1,\dots,\lambda_k)\in \Pi_{1\le j \le k}(\C\setminus \overline S_{\alpha_j, \beta_j})).$$

(iii) The Fourier-Borel transform of $f$ is defined on $\Pi_{1\le j \le k}(\C \setminus -\overline S^*_{\alpha_j,\beta_j})$ by the formula

$$\FB(f)(\lambda)=\int_0^{e^{i\omega}.\infty}e^{-\lambda \zeta}f(\zeta)d\zeta$$ $$:= \int_0^{e^{i\omega_1}.\infty}\dots \int_0^{e^{i\omega_k}.\infty}e^{-\lambda_1\zeta_1-\dots-\lambda_k \zeta_k}f(\zeta_1,\dots,\zeta_k)d\zeta_1\dots d\zeta_k $$ $$ (\lambda =(\lambda_1,\dots, \lambda_k) \in \Pi_{1\le j \le k}(\C\setminus -\overline S^*_{\alpha_j, \beta_j}) ),$$  

where $\alpha _j \le \omega_j \le \beta_j$ and where $Re(\lambda_je^{i\omega_j})>0$ for $1\le j \le k.$

\end{defn}

It follows from these definitions that ${\mathcal C}(\phi)$ is holomorphic on $\Pi_{1\le j \le k}(\C\setminus \overline S_{\alpha_j, \beta_j})$ for $\phi \in \uab',$ and that $\FB(f)$ is holomorphic on  $\Pi_{1\le j \le k}(\C \setminus -\overline S^*_{\alpha_j,\beta_j})$ for $f \in \mathcal V_{\alpha,\beta}(X)$. Also using proposition 10.2 we see that $\FB(\phi)\in \mathcal V^*_{\alpha, \beta}:=\mathcal V_{-{\pi\over 2}-\alpha, {\pi\over 2}-\beta}$ for $\phi \in \uab'.$

\begin{prop} Let $\phi \in \uab'.$ For $j \le k,$ set $I_{\eta,j}=({\pi\over 2}-\eta, {\pi\over 2}-\beta_j]$ for $\eta \in( \beta_j, \alpha_j +\pi],$ $I_{\eta,j}=(-{\pi\over 2} -\alpha_j, {\pi\over 2}-\beta_j)$ for $\eta \in (\alpha_j +\pi, \beta_j +\pi ],$ and set 
$I_{\eta}=  (-{\pi \over 2}- \alpha_j, {3\pi\over 2}-\eta)$ for $\eta \in (\beta_j +\pi, \alpha_j +2\pi).$ Then $I_{\eta,j} \subset [-{\pi\over 2}-\alpha_j,{\pi\over 2}-\beta_j],$ $cos(\eta +s)<0$ for $s \in I_{\eta,j},$ and if $\lambda =(\lambda_1,\dots, \lambda_k) \in \Pi_{1\le j \le k}(\C\setminus \overline S_{\alpha_j, \beta_j}),$ we have for $\omega=(\omega_1,\dots, \omega_k) \in \Pi_{1\le j \le k}I_{arg(\lambda_j),j},$

\begin{equation}{\mathcal C}(\phi)(\lambda)={1\over (2i\pi)^k}\int_{0}^{e^{i \omega} .\infty}e^{\lambda \sigma}\fb(\phi)(\sigma)d\sigma$$ $$:={1\over (2\pi i)^k}\int_0^{e^{i\omega_1}.\infty}\dots \int_0^{e^{i\omega_k}.\infty}e^{\lambda_1\sigma_1+\dots +\lambda_k\sigma_k}\FB(\phi)(\sigma_1,\dots, \sigma_k)d\sigma_1\dots d\sigma_k.\end{equation}

\end{prop}

Proof:  It follows from the definition of $I_{\eta,j}$ that $I_{\eta,j}\subset [-{\pi\over 2}-\alpha_j, {\pi\over 2}-\beta_j].$  In the second case we have obviously ${\pi\over 2} < \eta +s < {3\pi\over 2}$ for $s \in I_{\eta_j}.$ In the first case we have ${\pi\over 2} <\eta +s< {\pi\over 2} + \eta -\beta_j \le {3\pi\over 2} +\alpha_j -\beta_j<{3\pi\over 2}$ for $\omega \in I_{\eta,j}$ and in the third case we have ${3\pi\over 2} > \eta +s > \pi +\beta_j -{\pi\over 2} -\alpha_j >{\pi\over 2}$ for $\omega \in I_{\eta,j}.$ We thus see that $cos(\eta +s)<0$ for $\eta \in (\beta_j, 2\pi+\alpha_j),$ $\omega \in I_{\eta,j}.$

Now assume that $\lambda \in \Pi_{1\le j \le k}(\C\setminus \overline S_{\alpha_j, \beta_j}),$ let $\eta_j \in (\beta_j, 2\pi+\alpha_j)$ be a determination of $arg(\lambda_j),$ let $\nu$ be a representing measure for $\phi$ and let $\omega \in \Pi_{1\le j \le k}I_{\eta_j}.$ Then $\FB(\phi)$ is bounded on $\scabp,$ and since $cos(\eta_j +\omega_j)<0$ for $j\le k,$ we have

$${1\over (2i\pi)^k}\int_0^{e^{i\omega}.\infty}e^{\lambda \sigma}\fb(\phi)(\sigma)d\sigma={1\over (2i\pi)^k}\int_0^{e^{i\omega}.\infty}e^{\lambda \sigma}\left [ \int_{\overline S_{\alpha, \beta}}e^{-\sigma \zeta}d\nu(\zeta)\right ]d\sigma$$ $$= {1\over (2i\pi)^k}\int_{\overline S_{\alpha,\beta}}\left [\int_0^{e^{i\omega}.\infty}e^{\sigma(\lambda -\zeta)}d\sigma \right ]d\nu(\zeta)=\int_{\overline S_{\alpha,\beta}}{1\over \zeta -\lambda}d\nu(\zeta)=\mathcal C(\phi)(\lambda).$$ $\square$



 

Now identify the space ${\mathcal M}(\scab)$ of all measures of bounded variation on $\scab$ to the dual space of the space ${\mathcal{C}}_0(\scab)$ of continuous functions on $\scab$ vanishing at infinity via the Riesz representation theorem.The convolution product of two elements of ${\mathcal M}(\scab)$ is defined by the usual formula

$$\int_{\scab}f(\zeta)d(\nu_1*\nu_2)(\zeta):=\int_{\scab\times \scab}f(\zeta+\zeta')d\nu_1(\zeta)d\nu_2(\zeta') \ \  (f \in {\mathcal{C}}_0(\scab)).$$

\begin{prop} Let $X$ be a separable Banach space.

(i) For $f \in \vab(X), \lambda \in \scab,$ set $f_{\lambda}(\zeta)=f(\zeta +\lambda).$ Then $f_{\lambda}\in \vab(X)$ for $f \in \vab(X),$ $f_{\lambda}\in \uab(X)$ and the map $\lambda \to f_{\lambda}$ belongs to $\uab(\uab(X))$ for $f \in \uab(X).$ Moreover if we set, for $\phi \in \uab',$

$$f_\phi(\lambda)=<f_\lambda,\phi>,$$

then $f_{\phi}\in \vab(X)$ for $f \in \vab(X),$ and $f_{\phi}\in \uab(X)$ for $f \in \uab(X).$

(ii) For $\phi_1\in \uab',\phi_2\in \uab',$ set

$$<f, \phi_1*\phi_2>=<f_{\phi_1},\phi_2> \ \ \ (f\in \uab).$$

Then $\phi_1*\phi_2\in \uab',$ $\nu_1*\nu_2$ is a representing measure for $\phi_1*\phi_2$ if $\nu_1$ is a representing measure for $\phi_1$ and if $\nu_2$ is a representing measure for $\phi_2,$ and we have

$$<f , \phi_1*\phi_2>=<f_{\phi_1}, \phi_2> \ \ \ (f\in \vab(X)),$$

$$\FB(\phi_1*\phi_2)=\FB(\phi_1)\FB(\phi_2).$$

\end{prop}

Proof: These results follow from standard easy verifications which are left to the reader. We will just prove the last formula. Let $\phi_1\in \uab',$  $ \phi_2\in \uab'.$  We have, for $z=(z_1,\dots, z_k) \in \sscab, \lambda=(\lambda_1,\dots, \lambda_k)\in \scab, \zeta=(\zeta_1,\dots,\zeta_k)\in \scab,$

$$( e_{-z})_{\lambda}(\zeta)=e^{-z_1(\lambda_1+\zeta_1)\dots -z_k(\lambda_k+\zeta_k)}= e_{-z}(\lambda) e_{-z}(\zeta),$$

and so $( e_{-z})_{\lambda}= e_{-z}(\lambda) e_{-z},$ $( e_{-z})_{\phi_1}(\lambda)=<(e_{-z})_{\lambda},\phi_1>=e_{-z}(\lambda)\FB(\phi_1)(z),$ $( e_{-z})_{\phi_1}= \FB(\phi _1)(z) e_{-z},$ and

$$\FB(\phi_1*\phi_2)(z)= <( e_{-z})_{\phi_1}, \phi_2>=\FB(\phi _1)(z)< e_{-z},\phi_2>=\FB(\phi _1)(z)\FB(\phi_ 2)(z).$$ $\square$

For $\eta \in \scab,$ denote by $\delta_{\eta}$ the Dirac measure at $\eta.$ We identify $\delta_{\eta}$ to the linear functional $f \to f(\eta)$ on $\uab.$ With the above notations, we have, for $f \in \vab(X),$ $\phi \in \uab',$

$$f_{\delta _{\eta}}=f_{\eta}, <f,\phi*\delta_{\eta}>=<f_{\eta},\phi>.$$

If $f\in \uab(X),$ we have $\lim_{\stackrel{\epsilon \to 0}{_{\epsilon \in \sscab}}}\Vert  e_{-\epsilon}f -f\Vert_{\infty}=0=\lim_{\stackrel{\eta \to 0}{_{\eta \in \scab}}}\Vert f_{\eta} -f\Vert_{\infty}.$ We obtain, since  $\Vert  e_{-\epsilon}f\Vert_{\infty} \le \Vert f\Vert_{\infty}$ for $f \in \uab, \epsilon \in \sscab,$ 

\begin{equation}\lim_{\stackrel{\epsilon \to 0, \epsilon \in \sscab}{_{\eta \to 0, \eta \in \scab}}}\Vert e_{-\epsilon}f_\eta-f \Vert_{\infty}=0 \ \ (f \in \uab(X).\end{equation}

Now let $f \in \vab(X),$ let $\phi \in \uab',$ and let $\nu$ be a representative measure for $\phi.$ Since $<e_{\epsilon}f_{\eta}>=\int_{\scab}e^{-\epsilon \zeta}f(\zeta +\eta)d\nu(\zeta),$
and since $<e_{-\epsilon}f, \phi*\delta_{\eta}>=e^{-\epsilon \eta}<e_{-\epsilon}f_{\eta}, \phi>,$ it follows from the Lebesgue dominated convergence theorem that we have

\begin{equation}<f,\phi>=\lim_{\stackrel{\epsilon \to 0, \epsilon \in \sscab}{_{\eta \to 0, \eta \in \scab}}}<e_{-\epsilon} f_{\eta}, \phi> =\lim_{\stackrel{\epsilon \to 0, \epsilon \in \sscab}{_{\eta \to 0, \eta \in \scab}}}<e_{-\epsilon} f, \phi*\delta_{\eta}>\ \ (f \in \vab(X), \phi \in \uab').\end{equation}
In the following we will denote by $\tilde \partial \overline S_{\alpha, \beta}=\Pi_{1\le j \le k}\partial \overline S_{\alpha_j, \beta_j}$ the distinguished boundary of $ \overline S_{\alpha,\beta}$, where $\partial  \overline S_{\alpha_j, \beta_j}=(e^{i\alpha_j}.\infty,0]\cup[0, e^{i\beta_j}.\infty)$ is oriented from $e^{i\alpha_j}.\infty$ towards $e^{i\beta_j}.\infty.$

The following standard computations allow to compute in some cases $<f,\phi>$ by using the Cauchy transform when $\alpha_j <\beta_j$ for $j\le k.$ 


\begin{prop} Assume that $\alpha_j < \beta_j < \alpha_j+\pi$ for $1 \le j \le k,$ and let $\phi \in \uab'.$

If $f \in \vab(X),$ and if 

$$\int_{\tilde \partial \overline S_{\alpha,\beta}}\Vert f(\sigma)\Vert_X\vert d \sigma \vert <+\infty,$$


  then we have, for $\eta \in \sab,$

\begin{equation}<f_{\eta},\phi>=<f ,\phi*\delta_{\eta}>=\int_{\tilde \partial \overline S_{\alpha, \beta}}{\mathcal C}(\phi)(\sigma-\eta)f(\sigma)d\sigma .\end{equation}

 In particular we have, for  $f \in \vab(X),  \epsilon \in \ssab ,\eta \in \sab,$

\begin{equation} e^{-\epsilon \eta}<e_{-\epsilon}f_{\eta}, \phi>=<e_{-\epsilon}f ,\phi*\delta _{\eta}>=\int_{\tilde \partial \overline S_{\alpha, \beta}} e^{-\epsilon\sigma}{\mathcal C}(\phi)(\sigma-\eta)f(\sigma)d\sigma\end{equation}

\end{prop} 

Proof: Assume  that $f \in \vab(X)$ satisfies the condition $\sup_{\sigma \in \scab}(1+\vert \sigma\vert)^{2k}\Vert f(\sigma)\Vert<+\infty.$ 

 Let $\nu\in \mathcal M(\scab)$ be a representing measure for $\phi.$ For $R>0, j\le k,$  we denote by $\Gamma_{R,j}$ the Jordan curve $\{Re^{i\omega}\}_{\alpha_j \le \omega \le \beta_j}\cup [Re^{i\beta_j},0]\cup[0, Re^{i\alpha_j}],$ oriented counterclockwise. 

We have, for  $\eta \in \sab$, $\sigma \in \Pi_{1\le j \le k}\partial S_{\alpha_j, \beta_j},$ 

$$\left \vert {\mathcal C}(\phi)(\sigma-\eta)\right \vert \le {1\over (2\pi)^k} \Vert \phi\Vert_{\uab'}\Pi_{1\le j \le k}dist(\partial S_{\alpha_j,\beta_j} -\eta_j,\partial S_{\alpha_j,\beta_j})^{-1}.$$

It follows then from Fubini's theorem and Cauchy's formula that we have 
$$\int_{\tilde \partial \overline S_{\alpha, \beta}}{\mathcal C}(\phi)(\sigma-\eta)f(\sigma)d\sigma=\int_{\scab}\left [ {1\over 2i\pi)^k }\int_{\tilde \partial S_{\alpha, \beta}}{f(\sigma)\over \zeta -\sigma +\eta}d\sigma\right ]d\nu(\zeta)$$
$$=\int_{\scab}\lim_{R\to +\infty}{1\over (2i\pi)^k}\left [\int_{\Gamma_{R,1}}\dots \int_{\Gamma_{R,k}}{f(\sigma)\over  (\sigma_1 -\zeta_1 -\eta_1)\dots (\sigma_k-\zeta_k-\eta_k)} d\sigma\right ]
d\nu(\zeta)$$ $$= \int_{\scab}f(\zeta+\eta)d\nu(\zeta)= <f, \phi*\delta_{\eta}>.$$ Formula (26) follows from this equality applied to $ e_{-\epsilon}f.$ Taking the limit as $\epsilon \to 0, \epsilon \in \ssab$ 
in formula (26), we deduce formula (25) from the Lebesgue dominated convergence theorem.
$\square$

The following result is indeed standard, but we give a proof for the convenience of the reader.

\begin{prop} The linear span of the set $E_{\alpha,\beta}:=\{ f=e_{-\sigma} : \sigma \in \Pi_{j\le k}(0, e^{-i{\alpha_j +\beta_j\over 2}}.\infty)\}$ is dense in $\uab,$ and the Fourier-Borel transform is one-to-one on $\uab'.$
\end{prop}

Proof: Set $J_1=\{ j \in \{1,\dots,k\} \ | \ \alpha_j=\beta_j\},$ set $J_2:=\{ j \in \{1,\dots,k\} \ | \ \alpha_j<\beta_j\},$ denote by $\mathcal U_1$ the space of continuous functions on $\overline S_1=\Pi_{j \in J_1}\overline S_{\alpha_j, \beta_j}$ vanishing at infinity, set $S_2:= \Pi_{j \in J_2} S_{\alpha_j, \beta_j}$, and denote by $\mathcal U_2$ the space of continuous functions on $\overline S_2$ vanishing at infinity which satisfy the same analyticity condition as in definition 10.1 with respect to $\overline S_2.$ Also set $E_1 := \{ f=e_{-\sigma} : \sigma \in \Pi_{j\in J_1}(0, e^{-i{\alpha_j +\beta_j\over 2}}.\infty)\}$, and set $E_2:=\{ f=e_{-\sigma} : \sigma \in \Pi_{j\in J_2}(0, e^{-i{\alpha_j +\beta_j\over 2}}.\infty)\}.$ 

Assume that $J_1\neq \emptyset.$ Then the complex algebra $span(E_1)$ is self-adjoint and separates the point on $\mathcal U_1,$ and it follows from the Stone-Weierstrass theorem applied to the one-point compactification of $S_1$ that $span(E_1)\oplus \C.1$ is dense in $\mathcal U_1\oplus \C.1,$ which implies that $span(E_1)$ is dense in $\mathcal U_1$ since $\mathcal U_1$ is the kernel of a character on $\mathcal U_1\oplus \C.1.$

Now assume that $J_2\neq \emptyset,$ set $S_2^*=\Pi_{j \in J_2}S_{-{\pi\over 2}-\alpha_j, {\pi\over 2}-\beta_j},$ let $\phi \in U_2',$ and define the Cauchy transform and the Fourier-Borel transform of $\phi$ as in definition 10.3. Assume that $<f, \phi>=0$ for $\phi \in E_2.$ If $j \in J_2,$ then $g=0$ for every holomorphic function $g$ on $S^*_{\alpha_j,\beta_j}$  which vanishes on $(0, e^{-i{\alpha_j+\beta_j\over 2}}.\infty).$ An immediate finite induction shows then that $\fb(\phi)=0$ since $\fb(\phi)$ is holomorphic on $S_2^*.$ It follows then from  proposition 10.4 that $\mathcal C(\phi)=0,$ and it follows from (23) and (26) that $<f, \phi>=0$ for every $f \in \mathcal U_2.$ Hence $\phi=0,$ which shows that $span(E_2)$ is dense in $\mathcal U_2.$ This shows that $span(E_{\alpha,\beta})$  is dense in $\uab$ if $J_1=\emptyset$ or if $J_2=\emptyset.$


Now assume that $J_1\neq \emptyset$ and $J_2\neq \emptyset,$ and denote by $E\subset \uab$ the set of products $f=gh,$ where $g \in \mathcal U_1$ and $h\in \mathcal U_2.$
The space $\mathcal U_1=\mathcal C_0(\overline S_1)$ is a closed subsbace of codimension one of $\mathcal C(\overline S_1\cup\{\infty\}).$ Since the space $\mathcal C(K)$ has a Schauder basis for every compact space $K$, \cite{ba},\cite{se}, the space $\mathcal U_1$ has a Schauder basis. Identifying the dual space of $\mathcal U_1$ to the space of measures of bounded variation on $\overline S_1$,  this means that there exists a sequence $(g_n)_{n\ge 1}$ of elements of $\mathcal U_1$ and a sequence $(\nu_n)_{n\ge 1}$ of measures of bounded variation on $\overline S^1$ such that we have

$$g=\sum \limits_{n=1}^{+\infty}\left (\int_{\overline S_1}g(\eta)d\nu_n(\eta)\right )g_n \ \ (g \in \mathcal U_1),$$

where the series is convergent in $(\mathcal U_1, \Vert . \Vert_{\infty}).$

Set $P_m(g)=\sum \limits_{n=1}^m\left (\int_{\overline S_1}g(\eta)d\nu_n(\eta)\right )g_n$ for $g \in \mathcal U_1, m\ge 1.$ Then $P_m:\mathcal U_1 \to \mathcal U_1$ is a bounded linear operator, and $\lim \sup_{m\to +\infty }\Vert P_m(g)\Vert \le \Vert g \Vert <+\infty$ for every $g \in \mathcal U_1.$ It follows then from the Banach-Steinhaus theorem that there exists $M>0$ such that $\Vert P_m\Vert_{{\mathcal B}(\mathcal U_1)}\le M$ for $m\ge 1,$ a standard property of Schauder bases in Banach spaces.

Now let $\phi \in \uab'$ such that $<f, \phi>=0$ for $f \in E,$ let $\nu$ be a repesenting measure for $\phi,$ and let $f \in \uab.$ The function $f_{\zeta}= \eta \to f(\eta, \zeta)$ belongs to $U_1$ for $\zeta \in \overline S_2,$ and a routine verification shows that the function $h_n: \zeta \to \int_{\overline S_1}f_{\zeta}(\sigma)d\nu_n(\sigma)=\int_{\overline S_1}f(\zeta,\eta)d\nu_n(\eta)$ belongs to $U_2$ for $n\ge 1.$ Since the evaluation map $g\to g(\eta)$ is continuous on $\mathcal U_1$ for $\eta \in \overline S_1,$ we obtain, for $\eta \in \overline S_1, \zeta \in \overline S_2,$

$$f(\eta, \zeta)= \lim_{m \to +\infty} \sum \limits_{n=1}^{m}g_n(\eta)h_n(\zeta).$$

We have, for $m\ge 1,$  $\eta \in \overline S_1, \zeta \in \overline S_2,$

$$\left \vert \sum \limits_{n=1}^{m}g_n(\eta)h_n(\zeta)\right \vert \le \Vert P_m(f_{\zeta})\Vert_{\infty}\le M\Vert f_{\zeta}\Vert_{\infty}\le M\Vert f\Vert_{\infty}.$$

It follows then from the Lebesgue dominated convergence theorem that

$$\int_{\overline S_{\alpha, \beta}}f(\eta, \zeta)d\nu(\eta, \zeta)=\lim_{m\to +\infty}\sum \limits_{n=1}^m\int_{\overline S_{\alpha, \beta}}g_n(\eta)h_n(\zeta)d\nu(\eta, \zeta)=0.$$

This shows that $span(E)$ is dense in $\uab.$ Since $span(E_1)$ is dense in $\mathcal U_1$ and $span(E_2)$ is dense in $\mathcal U_2,$ $span(E_{\alpha, \beta})$ is dense in $span(E),$ and so $span(E_{\alpha, \beta})$ is dense in $\uab.$

Now let $\phi \in \uab'.$ If $\FB(\phi)=0,$ then $<f, \phi>=0$ for every $f \in E_{\alpha,\beta},$ and so $\phi=0$ since $span(E_{\alpha, \beta})$ is dense in $\uab,$ which shows that the Fourier-Borel transform is one-to-one on $\uab.$

$\square$

We will now give a way to compute $<f,\phi>$ for $f \in \vab(X), \phi \in \uab'$ by using Fourier-Borel transforms.
For $\sigma \in \Pi_{1\le j \le k}\left (\C \setminus \overline S^*_{\alpha_j, \beta_j}\right ),$ define $e^*_{\sigma}\in \uab'$ by using the formula

\begin{equation}<f,  e^*_{\sigma}>= \fb(f)(-\sigma). \end{equation}

Also for $\phi \in \uab', $ $g \in \uab,$ define $\phi g \in \uab'$ by using the formula

$$<f, \phi g>= <fg, \phi> \ \ \ (f\in \uab).$$

It follows from definition 10.3 that if $\sigma=(\sigma_1,\dots, \sigma_k) \in \Pi_{1\le j \le k}\left (\C \setminus \overline S^*_{\alpha_j, \beta_j}\right ),$ we have, for $f\in \uab,$

$$<f,  e^*_{\sigma}> =\int_{0}^{e^{i\omega}.\infty}e^{\sigma \zeta}f(\zeta)d\zeta,$$

where $\omega=(\omega_1,\dots,\omega_k)$ satisfies $\alpha_j \le \omega_j\le \beta_j, Re(\sigma_j \omega_j)<0$ for $1\le j \le k,$ which gives

$$\Vert e^*_{\sigma}\Vert_{\infty}\le \Pi_{1\le j \le k}\int_0^{\infty}e^{tRe(\sigma_j \omega_j)}dt ={1\over \Pi_{1\le j \le k}(-Re(\sigma_j\omega_j))}.$$

The same formula as above holds with the same $\omega$ to compute $<f, e^*_{\sigma'}>$ for $\sigma'\in  \Pi_{1\le j \le k}\left (\C \setminus \overline S^*_{\alpha_j, \beta_j}\right )$ when $\vert \sigma -\sigma'\vert$ is sufficiently small, and so the map $\sigma \to  e^*_\sigma\in \uab'$ is holomorphic on  $\Pi_{1\le j \le k}\left (\C \setminus \overline S^*_{\alpha_j, \beta_j}\right )$ since the map $\lambda \to e_{-\lambda}\in L^1(\R^+)$ is holomorphic on the open half-plane $P^+:=\{ \lambda \in \C \ | \ Re(\lambda)>0\}.$

Now let $\epsilon \in S^*_{\alpha, \beta}$ and let $\omega \in \Pi_{1\le j \le k}[\alpha_j, \beta_j]$ such that $Re(\epsilon_j e^{i\omega_{j}})>0$ for $j\le k.$ Then $\sigma - \epsilon \in 
\Pi_{1\le j \le k}\left (\C \setminus \overline S^*_{\alpha_j, \beta_j}\right )$ for $\sigma \in \tilde \partial S^*_{\alpha, \beta},$ and $Re((\sigma_j -\epsilon_j) \omega_j)\le -Re(\epsilon_j \omega_j)<0$ for $1 \le j \le k.$ We obtain

$$\Vert e^*_{\sigma -\epsilon}\Vert \le {1\over \Pi_{1\le j \le k}Re(\epsilon_j \omega_j)},$$

and so $\sup_{\sigma \in \tilde \partial S^*_{\alpha, \beta}}\Vert e^*_{\sigma - \epsilon}\Vert_{\infty}< +\infty \ \ (\epsilon \in S^*_{\alpha, \beta}).$

We now give the following certainly well-known natural result.

\begin{prop} Let $\phi \in \uab'.$ Assume that

$$\int_{\tilde \partial S^*_{\alpha, \beta}}\vert \fb(\phi)(\sigma)\vert \vert d\sigma\vert <+\infty.$$

Then we have, for $\epsilon \in S^*_{\alpha, \beta},$

$$\phi e_{-\epsilon}={1\over (2i\pi)^k}\int_{\tilde \partial S^*_{\alpha, \beta}}\fb(\phi)(\sigma) e^*_{\sigma -\epsilon}d\sigma,$$

where the Bochner integral is computed in $(\uab', \Vert .\Vert_{\infty}),$ which gives, for $f \in \mathcal V_{\alpha, \beta}(X),$

\begin{equation} <f e_{-\epsilon}, \phi>={1\over (2i\pi)^k}\int_{\tilde \partial S^*_{\alpha, \beta}} \fb(\phi)(\sigma)\fb(f)(-\sigma +\epsilon)d\sigma.\end{equation}

\end{prop}

Proof:  Since the map $\sigma \to  e^*_{\sigma -\epsilon}\in \uab'$ is continuous on $\tilde \partial S^*_{\alpha, \beta},$ and since $\sup_{\sigma \in \tilde \partial S^*_{\alpha, \beta}}\Vert  e^*_{\sigma - \epsilon}\Vert_{\infty}< +\infty, $ the Bochner integral $\int_{\tilde \partial S^*_{\alpha, \beta}}\fb(\phi)(\sigma) e^*_{\sigma -\epsilon}d\sigma$ is well-defined in 
$(\uab', \Vert .\Vert_{\infty}).$ Set $\phi_{\epsilon}:={1\over (2i\pi)^k}\int_{\tilde \partial S^*_{\alpha, \beta}}\fb(\phi)(\sigma) e^*_{\sigma -\epsilon}d\sigma\in \uab'.$ Since the map $\phi \to \FB(\phi)(\zeta)$ is continuous on $\uab',$  we have, for $\zeta \in \overline S_{\alpha, \beta},$

$$\fb(\phi_{\epsilon})(\zeta)={1\over (2i\pi)^k}\int_{\tilde \partial S^*_{\alpha, \beta}}\fb(\phi)(\sigma) \fb( e^*_{\sigma -\epsilon})(\zeta)d\sigma$$ $$={1\over (2i\pi)^k}\int_{\tilde \partial S^*_{\alpha, \beta}}\fb(\phi)(\sigma)<e_{-\zeta}, e^*_{\sigma -\epsilon}>d\sigma.$$

It follows from definition (27) that $<e_{-\zeta}, e^*_{\sigma -\epsilon}>= \fb(e_{-\zeta})(\epsilon -\sigma).$ Let $\omega \in \Pi_{1\le j \le k}[\alpha_j, \beta_j]$ such that $Re(\epsilon_j e^{i\omega_{j}})>0$ for $j\le k.$ Since  $Re((\sigma_j -\epsilon_j) \omega_j)\le -Re(\epsilon_j \omega_j)<0$ for $1 \le j \le k,$ we have, for $\sigma \in \tilde \partial S^*_{\alpha, \beta},$

$$\fb(e_{-\zeta})(\epsilon -\sigma)=\int_0^{e^{i\omega}.\infty}e^{(\sigma-\epsilon)\eta}e_{-\zeta}(\eta)d\eta=\int_0^{e^{i\omega}.\infty}e^{(\sigma-\epsilon-\zeta)\eta}d\eta$$ $$={1\over \Pi_{1\le j \le k}(\zeta_j+\epsilon_j-\sigma_j)}.$$

Using the notation ${1\over \zeta +\epsilon -\sigma}:={1\over \Pi_{1\le j \le k}(\zeta_j+\epsilon_j-\sigma_j)},$ this gives

$$\fb(\phi_{\epsilon})(\zeta)={1\over (2i\pi)^k}\int_{\tilde \partial S^*_{\alpha, \beta}}{\fb(\phi)(\sigma)\over \zeta +\epsilon -\sigma}d\sigma.$$

As in appendix 3 , set $W_{j,n}(\zeta_j)={n^2 \over \left (n +e^{{\alpha_j+\beta_j\over 2}i}\zeta_j\right )^2}$ for $n\ge 1, \zeta_j \in \overline S^*_{\alpha_j, \beta_j},$ and set $W_n(\zeta)=\Pi_{j\le k}W_{n,j}(\zeta_j)$
for $\zeta \in \overline S^*_{\alpha, \beta}.$
Then $\vert W_{n,j}(\zeta_j)\vert \le1$ for
$\zeta_j\in \overline S^*_{\alpha_j, \beta_j}$, $W_n(\zeta) \to 1$ as $n\to \infty$ uniformly on compact sets of $\overline S^*_{\alpha, \beta},$ and $\lim_{\stackrel{\vert \zeta \vert \to \infty}{_{\zeta \in S^*_{\alpha, \beta}}}}W_n(\zeta)
=0.$ The open set $S^*_{\alpha,\beta}$ is admissible with respect to $(\alpha, \beta)$ in the sense of definition 12.1 and, since $\fb(\phi)$ is bounded on  $S^*_{\alpha,\beta},$   $\fb(\phi)W_n\in H^{(1)}(S^*_{\alpha, \beta})$ for $n\ge 1.$ It follows then from theorem 12.5 that we have, for $t \in (0,1), \zeta \in \overline S^*_{\alpha, \beta},$

$${1\over (2i\pi)^k}\int_{\tilde \partial S^*_{\alpha, \beta}}{\fb(\phi)(\sigma +t\epsilon)W_n(\sigma+t\epsilon)\over \zeta +(1-t)\epsilon - \sigma}d\sigma$$ $$={1\over (2i\pi)^k}\int_{\tilde \partial S^*_{\alpha, \beta}+t\epsilon}{\fb(\phi)(\sigma)W_n(\sigma)\over \zeta +\epsilon - \sigma}d\sigma=\fb(\phi)(\zeta +\epsilon)W_n(\zeta+\epsilon),$$

and it follows from the Lebesgue dominated convergence theorem that we have

$${1\over (2i\pi)^k}\int_{\tilde \partial S^*_{\alpha, \beta}}{\fb(\phi)(\sigma)W_n(\sigma)\over \zeta +\epsilon - \sigma}d\sigma=\fb(\phi)(\zeta +\epsilon)W_n(\zeta+\epsilon).$$

Taking the limit as $n\to +\infty,$ and using again the Lebesgue dominated convergence theorem, we obtain, for $\zeta \in \overline S^*_{\alpha, \beta},$ 

$$\fb(\phi_{\epsilon})(\zeta)={1\over (2i\pi)^k}\int_{\tilde \partial S^*_{\alpha, \beta}}{\fb(\phi)(\sigma)\over \zeta +\epsilon - \sigma}d\sigma=\fb(\phi)(\zeta +\epsilon)$$ $$=<e_{-\zeta}e_{-\epsilon},\phi>=<e_{-\zeta},\phi e_{-\epsilon}>=\fb(\phi e_{-\epsilon})(\zeta),$$

and it follows from the injectivity of the Fourier-Borel transform on $\uab'$ that $\phi_{\epsilon}=\phi e_{-\epsilon}.$ 

This gives, for $f \in \mathcal V_{\alpha,\beta}(X),$ since $<f, e^*_{\sigma -\epsilon}>= \fb(f)( -\sigma+\epsilon)$ for $\sigma \in \tilde \partial S^*_{\alpha, \beta},$

$$<fe_{-\epsilon}, \phi>=<f, \phi e_{-\epsilon}>={1\over (2i\pi)^k}\int _{\tilde \partial S^*_{\alpha, \beta}}\fb(\phi)(\sigma)<f,  e^*_{\sigma -\epsilon}>$$ $$={1\over (2i\pi)^k}\int_{\tilde \partial S^*_{\alpha, \beta}} \fb(\phi)(\sigma)\fb(f)(-\sigma +\epsilon)d\sigma.$$

$\square$

For $J\subset \{1,\dots, k\},$ set $P_{J,j}=\C\setminus -S_{-{\pi\over 2}-\alpha_j, 
{\pi\over 2}-\alpha_j}$,  $\omega_{J,i}=\alpha_j$ for $j \in J,$ $P_{J,j}=\C \setminus -S_{-{\pi\over 2}-\beta_j, 
{\pi\over 2}-\beta_j}$, $\omega_{J,i}=\beta_j$ for $j\in \{1,\dots, k\} \setminus J,$ and set $P_J=\Pi_{1\le j \le k}P_{J,i},$ $\omega_J=(\omega_{J,1},\dots,\omega_{J,k}).$ If $f \in \vab(X),$ and if $ \int_{\tilde \partial \overline S_{\alpha, \beta}}\Vert f(\zeta)\Vert _X \vert d\zeta\vert <+\infty,$ then the formula $\fb(f)(\sigma)=\int_0^{e^{i\omega_J}.\infty}e^{-\zeta \sigma}f(\zeta)d\sigma$ defines a continuous bounded extension of $\fb(f)$ to $P_J.$ So in this situation $\fb(f)$ has a continuous bounded extension to $\cup_{J\subset \{1,\dots, k\}}P_J=\Pi_{1\le j \le k}\left (\C \setminus -S^*_{\alpha_j, \beta_j}\right ).$ Applying formula (28) to the sequence $(\epsilon_n)=({\epsilon \over n})$ for some $\epsilon \in S^*_{\alpha, \beta},$ we deduce from the Lebesgue dominated convergence theorem and from formula (23) the following result.

\begin{cor}

Let $f \in \vab(X),$ and let $\phi \in \uab'.$ Assume that the following conditions are satisfied

 $$(i) \int_{\tilde \partial \overline S_{\alpha, \beta}}\Vert f(\zeta)\Vert _X \vert d\zeta\vert <+\infty.$$ 

 $$(ii) \int_{\tilde \partial S^*_{\alpha, \beta}}\vert \fb(\phi)(\sigma)\vert \vert d\sigma\vert <+\infty.$$ 
 
 Then
 
 \begin{equation}<f, \phi>={1\over (2i\pi)^k}\int_{\tilde \partial S^*_{\alpha, \beta}} \fb(\phi)(\sigma)\fb(f)(-\sigma)d\sigma.\end{equation}
 
 \end{cor}

In the following we will denote by $\tilde \nu$ the functional $f \to \int_{\sab}f(\zeta)d\nu(\zeta)$ for $\nu \in \mathcal M(\scab).$ In order to give a way to compute $<f,\phi>$ for $\phi \in \uab', f \in \vab(X),$ we will use the following easy observation.

\begin{prop} Let $\nu$ be a probability measure on $\scab,$ let $R>0,$ and let $X$be a separable Banach space. Set $\nu_R(A)=\nu(RA)$ for every Borel set $A\subset \scab.$ Then $\lim_{R\to +\infty}\Vert f_{\tilde \nu_R
}-f\Vert_{\infty}=0$ for every $f \in \uab(X).$
\end{prop}

Proof: Let $f\in \uab(X).$ Then $f$ is uniformly continuous on $\scab,$ and so for every $\delta >0$ there exists $r >0$ such that $\Vert f(\zeta+\eta)-f(\zeta)\Vert_X< \delta$ for every $\zeta \in \scab$ and for every $\eta \in \scab \cap \overline B(0,r).$ It follows from the Lebesgue dominated convergence theorem that $\lim_{R\to +\infty}\nu_R(B(0,r))=\lim_{R\to +\infty}\nu(B(0,rR))=\nu(\scab)=1.$ This gives

$$\lim \sup_{R\to +\infty}\Vert f_{\tilde \nu_R}-f\Vert_{\infty}=\lim \sup_{R\to +\infty}\left (\sup_{\zeta \in \sab}\left \Vert \int_{\scab}(f(\zeta +\eta)-f(\zeta))d\nu_R(\eta)\right \Vert _X\right )$$ $$\le \lim \sup_{R\to +\infty}\left (\sup_{\zeta \in \scab}\int_{\scab \cap B(0,r)}\Vert f(\zeta +\eta) -f(\zeta)\Vert _X d\nu_n(\eta)\right )$$ $$+2\Vert f \Vert_{\infty}\lim \sup_{R \to +\infty}\int_{\scab\setminus (\scab \cap B(0,r))}d\nu_n(\eta)\le \delta.$$

Hence $\lim_{R\to +\infty}\Vert f_{\tilde \nu_R}-f\Vert_{\infty}=0.$ $\square$

It follows from the definition of $\nu_R$ that $<f,\tilde \nu_R>=<f_{1\over R}, \tilde \nu>$ for $f \in \vab(X),$ where $f_{1\over R}(\zeta)=f(R^{-1}\zeta) \ (\zeta \in \sab).$ In particular if $\fb(\tilde \nu_R)=\fb(\nu)_{1\over R},$ and $(\tilde \nu_1)_R*(\tilde \nu_2)_R=(\tilde \nu_1*\tilde \nu_2)_R=(\widetilde {\nu_1*\nu_2})_R$ for $R>0$ if $\nu_1$ and $\nu_2$ are two probability measures on $\sab.$

We deduce from proposition 10.8 and proposition 10.9 the following corollary, in which the sequence $(W_n)_{n\ge 1}$ of functions on $\overline S^*_{\alpha,\beta}$ introduced in appendix 3 and used in the proof of proposition 10.8 allows to compute 
$<f,\phi>$ for $\phi \in \uab',$ $f \in \vab(X)$ in the general case.

\begin{cor}  Set $W_n(\zeta)=\Pi_{1\le j \le k}{n^2\over \left (n+\zeta_j e^{i{\alpha_j+\beta_j\over 2}}\right )^2}$ for $n\ge 1,$ $\zeta=(\zeta_1,\dots, \zeta_k)\in \overline S^*_{\alpha, \beta}.$ Then
we have, for $\phi \in \uab',$ $f \in \vab(X),$
\begin{equation} <f,\phi>=\lim_{\stackrel{\epsilon \to 0}{_{\epsilon \in \sab}}}\left ( \lim_{n\to +\infty}{1\over (2i\pi)^k}\int_{\tilde \partial S^*_{\alpha, \beta}}W_n(\sigma)\fb(\phi)(\sigma)\fb(f)(\epsilon - \sigma)d\sigma \right ).\end{equation}
\end{cor}

Proof: Define a measure $\nu_0$ on $\scab$ by using the formula

$$<f, \nu_0>=\int_{[0,+\infty)^k}e^{-t_1\dots -t_k}f(t_1e^{i{\alpha_1+\beta_1\over 2}},\dots,t_ke^{i{\alpha_k+\beta_k\over 2}})dt_1\dots dt_k \ \ (f \in \mathcal C_{0}(\sab)).$$

Then $\nu_0$ and $\nu=\nu_0*\nu_0$ are probability measures on $\scab,$ and we have, for $\zeta \in \overline S^*_{\alpha, \beta},$

$$\fb(\tilde \nu_0)(\zeta)=\int_{[0,+\infty)^k}e^{-t_1\dots -t_k}e^{-t_1\zeta_1e^{i{\alpha_1+\beta_1\over 2}}-\dots-t_k\zeta_ke^{i{\alpha_k+\beta_k\over 2}}}dt_1\dots dt_k$$ 
$$ = \Pi_{1\le j \le k}{1\over 1+\zeta_je^{i{\alpha_j+\beta_j\over 2}}}.$$

Hence $\fb(\tilde \nu)=\fb(\tilde \nu_0^2)=W_1,$ and $\fb(\tilde \nu_n)=(W_1)_{1\over n}=W_n.$ It follows from (29) that we have, for $\epsilon \in \overline S^*_{\alpha, \beta},$

$$<fe_{-\epsilon}, \phi>=\lim_{n\to +\infty}<(fe_{-\epsilon})_{\tilde \nu_n},\phi>=\lim_{n\to +\infty}<fe_{-\epsilon}, \phi*\tilde \nu_n>$$ $$=\lim_{n\to+\infty}{1\over (2i\pi)^k}\int_{\tilde \partial S^*_{\alpha,\beta}}\fb(\phi*\tilde \nu_n)(\sigma)\fb(f)(\epsilon -\sigma)d\sigma$$ $$=\lim_{n\to +\infty}{1\over (2i\pi)^k}\int_{\tilde \partial S^*_{\alpha, \beta}}W_n(\sigma)\fb(\phi)(\sigma)\fb(f)(\epsilon - \sigma)d\sigma,$$

and the result follows from the fact that $<f, \phi>=\lim_{\stackrel{\epsilon \to 0}{_{\epsilon \in \sab}}}<fe_{-\epsilon}, \phi>.$ $\square$

\section{Appendix 2: An algebra of fast-decreasing holomorphic functions on products of sectors and half-lines and its dual}

In this section we will use the notations introduced in definition 4.1 for $\alpha=(\alpha_1,\dots,\alpha_k)\in \R^k$ and $\beta =(\beta_1, \dots, \beta_k)\in \R^k)$ satisfying $\alpha_j\le \beta_j<\alpha_j+\pi$ for $1\le j\le k.$ Notice that is $x\in \C, y \in \C,$ there exists $z \in \C$ such that $\left(x+\overline S^*_{\alpha_j, \beta_j}\right )\cap \left (y+\overline S^*_{\alpha_j, \beta_j}\right )=z+\overline S^*_{\alpha_j, \beta_j}.$ Such a complex number $z$ is unique if $\alpha_j < \beta_j.$ If $\alpha_j=\beta_j,$ then $\overline S^*_{\alpha_j, \beta_j}=\overline S_{\alpha_j-\pi/2, \alpha_j +\pi/2}$ is a closed half-plane, the family $\left \{x +\overline S^*_{\alpha_j, \beta_j} \right \}_{x\in \C}$ is linearly ordered with respect to inclusion and the condition $\left (x+\overline S^*_{\alpha_j, \beta_j}\right )\cap \left (y+\overline S^*_{\alpha_j, \beta_j}\right )=z+\overline S^*_{\alpha_j, \beta_j}$ defines a real line of the form $z_0+e^{i\alpha_j}\R,$ where $z_0 \in \{x,y\}.$

The following partial preorder  on $\C^k$ is the partial order associated to the cone $\sscab$ if $\alpha_j<\beta_j$ for $1\le j \le k.$

\begin{defn} (i) For $z=(z_{1},\dots, z_{k})\in \C^k$ and $z'=(z'_{1},\dots, z'_{k})\in \C^k,$  set $z\preceq z'$ if $z'\in z+\sscab.$

(ii) if $(z^{(j)})_{1\le j \le m}$ is a finite family of elements of $\C^k$ denote by $sup_{1\le j \le m}z_j$ the set of all $z \in \C^k$ such that $\cap_{1\le j \le k}\left (z^{(j)}+\overline S^*_{\alpha, \beta}\right )=z + \overline S^*_{\alpha, \beta}.$
\end{defn}

For $z=(z_1,\dots, z_k) \in \C^k,$  set $e^{z}=(e^{z_1},\dots, e^{z_k}),$ and denote again by $e_{z}:\C^k\to \C$ the map $(\zeta_1,\dots, \zeta_k) \to e^{z \zeta}=e^{z_1\zeta_1+\dots +z_k\zeta_k}.$ 

 It follows from (17) that $ e_{-z'}{\mathcal U}_{\alpha, \beta} \subseteq  e_{-z}{\mathcal U}_{\alpha, \beta}$ if $z\preceq z'.$

For $f\in  e_{-z}\vab,$ set $\Vert f \Vert _{ e_{-z}\vab}=\Vert  e_{z}f\Vert_{\infty},$ which defines a Banach space norm on $ e_{-z}\uab$ and $e_{-z}\vab.$

\begin{prop} 

\smallskip

(i)Set $\gamma_n= ne^{-i{\alpha +\beta\over 2}}$ for $n\ge 1.$ Then the sequence $(\gamma_n)_{n\ge 1}$ is cofinal in $(\C^k, \preceq).$

\smallskip

(ii) If $z\preceq z',$ then $e_{-z'}\uab$ is a dense subset of $(e_{-z}\uab, \Vert .\Vert_{ e_{-z}\uab}).$

\smallskip

(iii)  The set  $\cap_{z \in \C^k} e_{-z}\uab$ is a dense ideal of $\uab,$ which if a Fr\'echet algebra with respect to the family $(\Vert .\Vert_{e_{-\gamma_n}\uab})_{n\ge 1}.$

(iv) If $X$ is a separable Banach space, and if $z \in sup_{1\le j \le m}z^{(j)},$ then $e_{-z}\uab(X) =\cap_{1\le j \le m}e_{-z^{(j)}}\uab(X),$ $e_{-z}\vab(X) =\cap_{1\le j \le m}e_{-z^{(j)}}\vab(X),$ and $\Vert f \Vert_{e_{-z}\vab}=\max_{1\le j \le m}\Vert f \Vert_{e_{-z^{(j)}}\vab}$ for $f \in e_{-z}\vab.$

\end{prop}

Proof: (i) Let  $z=(z_1,\dots, z_k)\in \C,$ and let $j\le k.$ Since $\left ({\pi\over 2} -\beta_j \right)+\left (-{\pi\over 2}-\alpha_j \right )= -(\alpha_j +\beta_j),$ $t_{0,j}e^{-i{(\alpha_j +\beta_j)\over 2}}\in \partial (z_j+\overline{S}^*_{\alpha_j, \beta_j})$ for some $t_{0,j} \in \R,$ so $t e^{-i{(\alpha_j +\beta_j)\over 2}}\in z_j+\overline{S}^*_{\alpha_j, \beta_j}$ for every $t\ge t_{0,j}$, and  (i) follows.

\smallskip

(ii) Assume that $z\preceq z'.$ The fact that $e_{-z'}\uab \subset e_{-z}\uab$ follows from (16). Let $z'' \in z'+\ssab \subset z+\ssab.$ We have $z''=z+re^{i\eta}$ where $r>0,$ and where $\eta=(\eta_1,\dots,\eta_k)$ satisfies  $-{\pi\over 2}-\alpha_j < \eta_j < {\pi\over 2}-\beta_j$ for $j\le k.$

The semigroup $(e_{-te^{i\eta}})_{t >0}$ is analytic and bounded in the Banach algebra $\uab,$ and lim$_{t\to 0^+}\Vert f -fe_{-te^{i\eta}}\Vert_{\infty}=0$ for every $f \in \uab.$

It follows then from the analyticity of this semigroup that $\left [e_{-re^{i\eta}}\uab\right ]^-=\left [\cup_{t>0}e_{-te^{i\eta}}\uab\right ]^-=\uab.$ Hence $e_{-z''}\uab$ is dense in $e_{-z_1}\uab,$ which proves (ii) since $e_{-z''}\uab \subset  e_{-z'}\uab.$

\smallskip

(iii) Denote by $i_{z,z'}: f \to f$ the inclusion map from $ e_{-z'}\uab$ into $ e_{-z}\uab$ for $z\preceq z'.$ Equipped with these maps, the family $( e_{-z}\uab)_{z\in \C^k}$ is a projective system of Banach spaces, and we can identify $\left (\cap_{z\in \C^k} e_{-z}\uab, (\Vert .\Vert_{ e_z\uab})_{z\in \C^k}\right )$ to the inverse limit of this system, which defines a structure of complete locally convex topological space on $\left (\cap_{z\in \C} e_{-z}\uab, (\Vert .\Vert_{ e_{-z}\uab})_{z\in \C^k}\right )$. 
It follows from (i) that the sequence $(\Vert .\Vert_{e_{-\gamma_n}\mathcal U_{\alpha, \beta}})_{n\ge 1}$ of norms defines the same topology as the family $\left(\Vert .\Vert_{e_{-z}\uab})_{z\in \C}\right )$ on $\cap_{z \in \C} e_{-z}\uab=\cap_{n\ge 1}e_{-\gamma_n}\uab,$ which defines a Fr\'echet algebra structure on $\cap_{z \in \C} e_{-z}\uab$. 

It follows from (ii) that $e_{-\gamma_{n+1}}\uab$ is dense in $e_{-\gamma_{n}}\uab$ for $n\ge 0,$ and a standard application of the Mittag-Leffler theorem of projective limits of complete metric spaces, see for example theorem 2.14 of \cite{e1}, shows that $\cap_{z \in \C} e_{-z}\uab=\cap_{n\ge 1}e_{-\gamma_n}\uab$ is dense in $e_{-\gamma_0}\uab=\uab.$

\smallskip
(iv) Let $z=(z_1,\dots, z_k)\in \C^k,$ let  $z'=(z'_1,\dots z'_k) \in \C^k,$ and let $z"=(z"_1,\dots,z"_k)\in sup(z,z').$ Then $e_{-z"+z}\in \vab, e_{-z"+z'}\in \vab,$ $\Vert e_{-z"+z}\Vert_{\infty} \le 1,  \Vert e_{-z"+z}\Vert _{\infty}\le 1,$ and so $e_{-z"}\uab(X)\subset e_{-z}\uab(X)\cap e_{-z'}\uab(X),$ $e_{-z"}\vab(X)\subset e_{-z}\vab(X)\cap e_{-z'}\vab(X),$
and $\max(\Vert f \Vert_{e_{-z}\vab(X)}, \Vert f \Vert_{e_{-z'}\vab(X)})\le \Vert f \Vert_{e_{-z"}\vab(X)}$ for $f \in e_{-z"}\vab(X).$
We claim that $\vert e^{z_j"\zeta}\vert \le \min ( \vert e^{z_j\zeta}\vert , \vert e^{z'_j\zeta}\vert)$ for $\zeta \in \partial \overline S_{\alpha_j, \beta_j},$  $1 \le j \le k.$ If $z_j \in z'_j+\overline S^*_{\alpha_j,\beta_j},$ or if $z'_j \in z_j+\overline S^*_{\alpha_j,\beta_j},$ this is obviously true. Otherwise we have $\alpha_j<\beta_j$ and, say, $z_j"=z_j+re^{i(-\alpha_j-{\pi\over 2})}=z_j'+r'e^{i(-\beta_j+{\pi\over 2})},$ with $r>0,r'>0.$ Let $\zeta =\rho e^{i\theta} \in \overline S_{\alpha_j,\beta_j},$ where $\rho\ge 0, \theta \in [\alpha_j,\beta_j].$ We have $Re((z"_j-z_j)\zeta)=r\rho cos( \theta-\alpha_j -{\pi\over 2})\ge 0,$ and $Re((z"_j-z_j)\zeta)= r'\rho cos (\theta -\beta_j +{\pi\over 2})\ge 0.$ So  $\vert e^{z_j"\zeta} \vert =\vert e^{z\zeta}\vert \le \vert e^{z'\zeta}\vert $ if $\theta =\alpha_j,$ and $\vert e^{z_j"\zeta} \vert =\vert e^{z'\zeta}\vert \le \vert e^{z'\zeta}\vert $ if $\theta =\beta_j,$ which proves the claim. 

We now use the Phragm\'en-Lindel$\ddot{\mbox{o}}$f principle. Let $s \in \cup_{1\le j \le k}(1,{\pi\over \beta_j-\alpha_j})$ and for $1 \le j \le k$ let $\zeta_j^s$ be a continuous determination of the $s$-power of $\zeta$ on $\overline S_{\alpha_j,\beta_j}$ which is holomorphic on $S_{\alpha_j,\beta_j}$ when $\alpha_j < \beta_j.$ Set $\zeta^s=\zeta_1^s\dots\zeta_k^s$ for $\zeta \in \overline \sab.$ Let $f \in e_{-z}\vab(X)\cap e_{-z'}\vab(X),$ and let $\epsilon >0.$ Set $g_{\epsilon}(s)=e^{-\epsilon \zeta^s}e^{z"\zeta }f(\zeta)$ for $\zeta \in \overline \sab.$ It follows from the claim that  $g_\epsilon \in \uab(X),$ and it follows from the maximum modulus principle that there exists $\zeta_0\in \tilde \partial S_{\alpha, \beta}$ such that $\Vert g_\epsilon\Vert_{\uab(X)}= \Vert g_{\epsilon}(\zeta_0)\Vert
\le \vert e^{-\epsilon \zeta_0^s}\vert \max(\Vert f\Vert _{e^{-z}\vab(X)}, \Vert  f\Vert _{e^{-z'}\vab(X)})\le  \max(\Vert f\Vert _{e^{-z}\vab(X)}, \Vert  f\Vert _{e^{-z'}\vab(X)}).$ Since $\lim_{\epsilon \to 0}e^{-\epsilon \zeta^s}=1$ for every $\zeta \in \overline \sab,$ this shows that $f \in e^{-z"}\vab(X),$ and $\Vert f \Vert_{e_{-z"}\vab(X)}=\max(\Vert f \Vert_{e_{-z}\vab(X)},\Vert f \Vert_{e_{-z'}\vab(X)}).$

Now let $f \in e_{-z}\uab(X)\cap e_{-z'}\uab(X).$ Then $f \in e_{-z"}\vab,$ and $\lim_{\stackrel{\vert \zeta \vert\to 0}{_{\zeta \in \partial \overline S_{\alpha,\beta}}}}\Vert e^{z"\zeta}f(\zeta)\Vert=0.$ The Banach algebra
$\uab$ possesses a bounded approximate identity $(g_n)_{n\ge 1},$ one can take for example $g_n(\zeta_1,\dots , \zeta_n)=\Pi_{1 \le j \le k}{n\zeta_j\over n\zeta_j +e^{i{\alpha_j+\beta_j\over 2}}}.$ We have $$\lim_{n\to +\infty}\Vert e_{z"}f g_n-e_{z"}f\Vert_{\infty}= \lim_{n\to +\infty} \max_{\zeta \in \partial \scab}\Vert e_{z"\zeta}f(\zeta)g_n(\zeta)- e_{z"\zeta}f(\zeta)\Vert =0,$$ and so $ e_{z"}f\in \uab$ since $\uab$ is a closed subalgebra of $\vab.$ This concludes the proof of (iv) when $m=2.$ The general case follows by an immediate induction,
since $\sup(\zeta, z^{(l)})=\sup_{1\le j \le l} z^{(j)})$ for every $\zeta \in \sup_{1\le j \le l-1} z^{(j)}$ if $(z^{1)},\dots, z^{(l)})$ is a finite family of elements of $\C^k.$

 $\square$ 

\smallskip

Notice that assertions (ii) and (iii) of the proposition do not extend to the case where $\beta_j =\alpha_j +\pi$ for some $j\le k.$ It suffices to consider the case where $\alpha_j=-{\pi\over 2}, \beta_j={\pi\over 2}.$ Set $\lambda_j(t)= (\lambda_{s,t})_{1\le s \le k},$ where $\lambda_{s,t}=0$ for $s\neq j$ and $\lambda_{j,t}=t.$ Then the map $f \to e_{-\lambda_j(t)}f$ is an isometry on $\uab$ for every $t \ge 0$ and $\cap_{t>0} u_{\lambda_j(t)}\uab=\{0\}$ since the zero function is the only bounded holomorphic function $f$ on the right-hand open half-plane  satisfying $\lim_{r\to +\infty}\vert e^{tr}f(r)\vert =0$ for every  $t>0.$

Let $i_{\zeta}: f \to f$ be the inclusion map from $\cap_{z \in \C}e_{-z}\uab$ into $e_{-\zeta}\uab.$ Since $i_{\zeta}$ has  dense range, the map $i_{\zeta}^*: \phi \to \phi_{|_{\cap_{z\in \C} e_{-z}\uab}}$ is a one-to-one map from $(e_{-\zeta} \uab)'$ into $\cap_{z \in \C} e_{-z}\uab)',$ which allows to identify $(e_{-\zeta}\uab)'$ to a subset of $(\cap_{z \in \C} e_{-z}\uab)',$ so that we have

\begin{equation}(\cap_{z \in \C} e_{-z}\uab)'=\cup_{z\in \C}( e_{-z} \uab)'=\cup_{n\ge 1}(e_{-ne^{-i{\alpha +\beta\over 2}}}\uab)'.\end{equation}


\begin{defn}  Set ${\mathcal F}_{\alpha, \beta}:=(\cap_{z\in \C^k} e_{-z}\uab)'.$ Let $\phi \in \fab,$ and let $X$ be a separable Banach space.

\smallskip

(i) The domain of the Fourier-Borel transform of $\phi$ is defined by the formula

 $$Dom(\FB(\phi)):=\{z \in \C^k \ | \ \phi 
\in (e_{-z}\uab)'\}.$$ 

\smallskip

(ii) For $z \in Dom(\FB(\phi))$ the functional $\phi e_{-z}\in \uab'$ is defined by the formula

$$<f, \phi e_{-z}>=< e_{-z}f, \phi> \ \ (f\in \uab),$$

\smallskip

and $<g,\phi>$ is defined for $g\in  e_{-z}\vab(X)$ by the formula

$$<g,\phi>=< e_{z}g,\phi e_{-z}>.$$

(iii) The Fourier-Borel transform of $\phi$ is defined for $z \in \dfi$ by the formula

$$\FB(\phi)(z) = < e_{-z},\phi>.$$

(iv) The $z$-Cauchy transform of $\phi$ is defined on $\C^k\setminus -\sscab$ for $z \in \dfi$ by the formula

$${\mathcal C}_z(\phi)  =\mathcal C(\phi  e_{-z}).$$

(v) If $z \in \dfi$ a measure $\nu$ of bounded variation on $\scab$ is said to be a $z$-representing measure for $\phi$ if $\nu$ is a representing measure for $\phi  e_{-z}.$

\end{defn}

Since the map $\zeta \to e_{-\zeta}$ is holomorphic on $S^*_{\alpha, \beta},$ the map $z \to  e_{-z}$ is a holomorphic  map from $\lambda +S^*_{\alpha, \beta}$ into $e_{-\lambda}\uab$ for every $\lambda \in \dfi,$ and so $\FB(\phi)$ is holomorphic on the interior of $\dfi$ for $\phi \in \fab.$ Also the $z$-Cauchy transform ${\mathcal C}_z(\phi)$ is holomorphic on $\C \setminus \scab$ for every $z \in \dfi.$
Notice also that if $\phi \in \uab',$ then $\sscab\subset Dom(\FB(\phi))$ and so the function $\FB(\phi)$ defined above is an extension to $Dom(\FB(\phi))$ of the Fourier-Borel transform already introduced in definition 10.3 on $\sscab.$

Now let $z\in \dfi, z'\in \dfi$ and assume that $g \in e_{-z}\vab(X)\cap e_{ -z'}\vab(X).$ Let $z''\in\sup(z,z')\subset Dom(\FB(\phi)).$  Then $g \in e_{-z"}\vab(X).$ Let $\nu$ be a $z$-representative measure for $\phi.$ We have, for $h\in \cap_{\lambda \in \C^k}e_{-\lambda} \uab,$ since $ e_{-z}=e_{-z''}e_{z''-z}.$

$$<h, \phi>=\int_{\scab} e_{z}(\zeta)h(\zeta)d\nu({\zeta})= \int_{\scab}e_{z''}(\zeta)h(\zeta)e_{z-z''}(\zeta)d\nu(\zeta).$$

Since $e_{z-z''}\nu$ is a measure of bounded variation on $\sab,$ $e_{z-z"}\nu$ is a $z''$-representative measure for $\phi.$ Similarly if $\nu'$ is a $z'$-representative measure for $\phi$ then $e_{z'-z''}\nu'$ is a $z''$-representative measure for $\phi,$ and we have

$$\int_{\scab} e_{z}(\zeta)g(\zeta)d\nu(\zeta)=\int_{\scab}e_{z''}(\zeta) g(\zeta)e_{z-z''}(\zeta)d\nu(\zeta)$$ $$=\int_{\scab}e_{z''}(\zeta) g(\zeta)e_{z'-z''}(\zeta)d\nu'(\zeta)=\int_{\scab}e_{z'}(\zeta)g(\zeta)d\nu'(\zeta),$$

which shows that the definition of $<g, \phi>$ does not depend on the choice of $z \in Dom(\FB(\phi))$ such that $g\in e_{-z}\vab(X).$

\begin{prop}  Let $\phi \in \fab.$ 

(i) The set $\dfi$ is connected.

(ii)    $z +\sscab \subset Dom(\FB(\phi))$, and $\FB(\phi)$ is continuous on $z+\sscab$ and holomorphic on $z+\ssab$ for every $z \in \dfi.$
\end{prop}

Proof: (i) The fact that $\dfi$ is connected follows from the fact that the arcwise connected set $(z_1+\sscab)\cup(z_2+\sscab)$ is contained in $\dfi$ for $z_1\in \dfi, z_2\in \dfi.$

(ii) Let $z \in \dfi.$ It follows from (16) that $z+\sscab \subset Dom(\FB(\phi))$ and so $\FB(\phi)$ is holomorphic on the open set $z+\ssab \subset Dom(\FB(\phi)).$ Let $\nu$ be a  measure of bounded variation on $\scab$  which is $z$-representing measure  for $\phi.$ We have, for $\eta \in \scab,$

$$\FB(\phi)(z+\eta)=< e_{-z-\eta},\phi>=<e_{-\eta},\phi e_{-z}>\int_{\scab}e^{-\eta \zeta}d\nu(\zeta),$$

and the continuity of $\FB(\phi)$ on $z+\sscab$ follows from the Lebesgue dominated convergence theorem.

$\square$

Notice that $\dfi$ is not closed in general: for example if we set $<f,\phi> =\int_{\overline{S}_{-{\pi\over 4},{\pi\over 4}}}\zeta f(\zeta)dm(\zeta)$ for $f\in \cap_{z\in \C} e_{-z}{\mathcal U}_{-{\pi\over 4},{\pi\over 4}},$ where $m$ denotes the Lebesgue measure on $\C,$ then $t\in \dfi$ for every $t>0,$ but $0\notin \dfi.$ Notice also that if $\nu$ is a measure supported by a compact subset of $\sab,$ and if we set $<f,\phi>:=\int_{\sab}f(\zeta)d\nu(\zeta)$ for $f \in \cap_{z \in \C^k} e_{-z}\uab,$ then $\phi \in \cap_{z\in \C^k}( e_{-z}\uab)',$  so that $\dfi=\C^k,$ and ${\mathcal {FB}}(\phi)$ is the entire function defined on $\C^k$ by the formula

$${\mathcal {FB}}(\phi)(z)=\int_{\sab}e^{-z\zeta}d\nu(\zeta).$$


\smallskip

We now introduce the convolution product of elements of $\fab.$ If $\phi \in \fab, f \in \cap_{z \in \C^k} e_{-z}\uab,$ $\lambda \in \scab,$ set again $f_\lambda(\zeta)=f( 
\zeta +\lambda)$ for $\zeta \in \scab.$ Then $f_\lambda \in  \cap_{z^k \in \C} e_{-z}\uab,$ and we can compute $<f_\lambda, \phi>.$ The map $\lambda \to f_\lambda$ is a continuous map from $\scab$ into the Fr\'echet algebra $\cap_{z \in \C} e_{-z}\uab$ which is holomorphic  on $\sab.$ We obtain  

\begin{lem} Let $\phi \in \fab.$ Then the function $f_\phi:\lambda \to <f_\lambda,\phi>$ belongs to $\cap_{z \in \C^k} e_{-z}\uab$ for every $f$ in $\cap_{z \in \C^k} e_{-z}\uab,$ and the linear map $f\to f_{\phi}$ is continuous on $\cap_{z \in \C^k} e_{-z}\uab.$

\end{lem}

Proof: Let $f \in \cap_{z \in \C^k} e_{-z}\uab,$ let $z_0\in \dfi,$ let $\nu$ be a $z_0$-representing measure for $\phi$ on $\scab,$ and  let $z \in \C^k.$ Let $z_1\in \sup(z_0,z),$ so that $(z_0 +\sscab)\cap(z+\sscab)=z_1+\sscab,$ and set $\eta_0=z_1-z_0,\eta=z_1-z.$ We have, for $\lambda \in \scab,$

$$e^{z\lambda}<f_{\lambda},\phi>=\int_{\scab}e^{z\lambda+z_0\zeta}f(\zeta +\lambda)d\nu(\zeta)= \int_{\scab}e^{-\eta \lambda -\eta_0\zeta}e^{z_1(\zeta +\lambda)}f(\zeta +\lambda)d\nu(\zeta).$$

Since $\vert e^{-\eta \lambda -\eta_0\zeta}e^{z_1(\zeta +\lambda)}f(\zeta +\lambda)\vert\le \Vert e_{z_1}f\Vert_{\infty},$ it follows from Lebesgue's dominated convergence theorem that
$\lim_{\stackrel{\vert \lambda \vert \to +\infty}{_{\lambda \in \scab}}}\vert e^{z \lambda}<f_{\lambda},\phi>\vert=0,$ and so $f_{\phi}\in \cap_{z \in \C^k} e_{-z}\uab.$ Also $\Vert  e_{z}f_{\phi}\Vert _{\infty}\le \Vert e_{z_1}f\Vert_{\infty}\int_{\scab}d\vert \nu\vert(\zeta),$ which shows that the map $f \to f_{\phi}$ is continuous on $\cap_{z^k \in \C} e_{-z}\uab.$ $\square$

\smallskip

Notice that it follows from the Hahn-Banach theorem that given $\phi \in (e_{-z_0}\uab)'$ there exists a $z_0$-representing measure $\nu$ for $\phi$ such that $\int_{\scab}d\vert \nu\vert(\zeta)= \Vert \phi\Vert_{(e_{-z_0}\uab)'}.$ The calculation above shows then that we have, for $z \in \C,$ $f \in \cap_{z \in \C} e_{-z}\uab,$ $\phi \in \fab,$ $z_0\in Dom(\fb(\phi)),$
$z_1\in \sup(z_0,z),$
\begin{equation} \Vert  e_{z}f_\phi\Vert_{\infty} \le \Vert e_{z_1}f\Vert_{\infty}\Vert \phi\Vert_{(e_{-z_0}\uab)'}.\end{equation}

\begin{prop} For $\phi_1\in \fab, \phi_2\in \fab,$ define the convolution product $\phi_1*\phi_2\in \fab$ by the formula

$$<f, \phi_1*\phi_2>=<f_{\phi_1},\phi_2> \ \ (f\in \cap_{z \in \C^k} e_{-z}\uab).$$

Then $sup(z_1,z_2)\subset Dom(\FB(\phi_1*\phi_2))$ for $z_1 \in Dom(\FB(\phi_1)), z_2\in Dom(\FB(\phi_2)),$ and we have, for $z \in \sup(z_1,z_2),$
$$\Vert \phi_1*\phi_2\Vert_{(e_{-z}\uab)'}\le \Vert \phi_1\Vert_{( e_{-z_1}\uab)'}\Vert \phi_2\Vert_{( e_{-z_2}\uab)'}.$$

More generally $Dom(\FB(\phi_1))\cap Dom(\FB(\phi_2))\subset Dom(\FB(\phi_1*\phi_2)),$ and if $z \in Dom(\FB(\phi_1))\cap Dom(\FB(\phi_2))$ then $(\phi_1*\phi_2) e_{-z}=(\phi_1 e_{-z})*(\phi_2 e_{-z}),$ so that $\nu_1*\nu_2$ is a $z$-representative measure for $\phi_1*\phi_2$ if $\nu_1$ is a $z$-representing measure for $\phi_1$ and  if $\nu_2$ is a $z$-representing measure for $\nu_2,$  and we have 

$$\FB(\phi_1*\phi_2)(z)=\FB(\phi_1)(z)\FB(\phi_2)(z) \ \ \ (z \in Dom(\FB(\phi_1))\cap Dom(\FB(\phi_2)).$$

\end{prop}

Proof: Let $z_1 \in Dom(\FB(\phi_1))$, let $ z_2\in Dom(\FB(\phi_2)),$ and let $z \in \sup(z_1,z_2).$ It follows from (32) that we have, for $f \in  \cap_{z \in \C^k} e_z\uab,$

$$\vert < f, \phi_1*\phi_2>\vert=\vert <f_{\phi_1},\phi_2>\vert \le \Vert e_{z_2}f_{\phi_1}\Vert_{\infty}\Vert \phi_2\Vert_{( e_{-z_2}\uab)'}$$ $$\le \Vert e_{z}f\Vert_{\infty} \Vert \phi_1\Vert_{( e_{-z_1}\uab)'}\Vert \phi_2\Vert_{( e_{-z_2}\uab)'}.$$

Hence $\phi_1*\phi_2 \in \fab,$ $sup(z_1,z_2) \subset Dom(\FB(\phi_1*\phi_2)),$ and $\Vert \phi_1*\phi_2\Vert_{(e_{-z}\uab)'}\\ \le \Vert \phi_1\Vert_{( e_{-z_1}\uab)'}\Vert \phi_2\Vert_{( e_{-z_2}\uab)'}$ for $z \in \sup(z_1,z_2).$

Let $z \in Dom(\FB(\phi_1))\cap \FB(\phi_2)).$ Then $z \in \sup(z,z) \subset Dom(\FB(\phi_1*\phi_2)).$

Let $\nu_1$ be a  $z$-representing measure for $\phi_1$ and  let $\nu_2$ be a $z$-representing measure for $\phi_2.$ We have, for $f \in  \cap_{s \in \C^k}e_{-s}\uab,$

$$<f,\phi_1*\phi_2>=<f_{\phi_1},\phi_2>=\int_{\scab}e^{z\lambda}f_{\phi_1}(\lambda)d\nu_2(\lambda)$$ $$=\int_{\scab}\left [ \int_{\scab}e^{z\zeta}f(\zeta +\lambda)d\nu_1(\lambda)\right ]e^{z\lambda}d\nu_2(\lambda)$$ $$=\int \int_{\scab\times \scab}e^{z(\zeta +\lambda)}f(\zeta +\lambda)d\nu_1(\zeta)d\nu_2(\lambda)=\int_{\scab}e^{zs}d(\nu_1*\nu_2)(s),$$

and so $\nu_1*\nu_2$ is a representing measure for $(\phi_1*\phi_2) e_{-z},$ which means that $\nu_1*\nu_2$ is a $z$-representative measure for $\phi_1*\phi_2.$ Since $\nu_1$ is a representative measure for $\phi_1 e_{-z},$ and since $\phi_2$ is a representative measure for $\phi_2 e_{-z},$ it follows from proposition 10.5 (ii) that $(\phi_1*\phi_2) e_{-z}=(\phi_1 e_{-z})*(\phi_2 e_{-z}).$ 
 
 It follows also from proposition 10.5(ii) that
 
 $$\FB(\phi_1*\phi_2)(z)=\FB( (\phi_1*\phi_2)e_{-z})(1)=\FB(( \phi_1e_{-z})*( \phi_2e_{-z}))(1)$$ $$=\FB( \phi_1e_{-z})(1)\FB(\phi_2e_{-z})(1)=\FB(\phi_1)(z)\FB(\phi_2)(z).$$
 
$ \square$

Using proposition 10.4, we obtain the following link between $z$-Cauchy transforms and Fourier-Borel transforms of elements of $\fab.$

\begin{prop} Let $\phi \in \fab.$ For $j \le k,$ set $I_{\eta,j}=({\pi\over 2}-\eta, {\pi\over 2}-\beta_j]$ for $\eta \in( \beta_j, \alpha_j +\pi],$ $I_{\eta,j}=(-{\pi\over 2} -\alpha_j, {\pi\over 2}-\beta_j)$ for $\eta \in (\alpha_j +\pi, \beta_j +\pi ],$ and set 
$I_{\eta}=  (-{\pi \over 2}- \alpha_j, {3\pi\over 2}-\eta)$ for $\eta \in (\beta_j +\pi, \alpha_j +2\pi).$ Then $I_{\eta,j} \subset [-{\pi\over 2}-\alpha_j,{\pi\over 2}-\beta_j],$ $cos(\eta +s)<0$ for $s \in I_{\eta,j},$ and if $\lambda =(\lambda_1,\dots, \lambda_k) \in \C^k\setminus \scab,$ we have for $\omega=(\omega_1,\dots, \omega_k) \in \Pi_{1\le j \le k}I_{arg(\lambda_j),j},$ $z\in Dom(\FB(\phi)),$

$${\mathcal C}_z(\phi)(\lambda)={1\over (2 i\pi)^k}\int_0^{e^{i\omega}.\infty}e^{\lambda \sigma}\FB(\phi)(\sigma+z)d\sigma$$ \begin{equation}:={1\over (2\pi i)^k}\int_0^{e^{i\omega_1}.\infty}\dots \int_0^{e^{i\omega_k}.\infty}e^{\lambda \sigma}\FB(\phi)(\sigma+z)d\sigma.\end{equation}

\end{prop}

Proof: We have ${\mathcal C}_z(\phi)={\mathcal C}(\phi  e_{-z}),$ and, for $\sigma \in \sscab,$ $$\FB(\phi)(\sigma+z)= <e_{-\sigma-z},\phi>=\\<e_{-\sigma} e_{-z},\phi>=<e_{-\sigma},\phi  e_{-z}>=\FB(\phi  e_{-z})(\sigma).$$ Applying formula (22) to $\phi  e_{-z},$ we obtain (33). $\square$

Let $X$ be a separable Banach space. For $\eta \in \scab,$ $z\in \C^k,$ $f \in  e_{-z}\vab(X),$ set $f_{\eta}(\zeta)=f(\zeta +\eta) \ (\zeta \in \scab).$  If $\phi \in \fab,$ and if $z \in Dom(\fb(\phi)),$ we have

$$<f, \phi*\delta_{\eta}>=<e_zf, (\phi*\delta_{\eta})e_{-z}>=<e_zf, (\phi e_{-z})*(\delta_{\eta}e_{-z})>=e^{-z\eta}<e_zf,(\phi e_{-z})*\delta_{\eta}>$$
$$=e^{-z\eta}<(e_zf)_{\eta},(\phi e_{-z})>=<e_zf_{\eta},\phi e_-z>=<f_\eta,\phi>.$$

We also have, for $f \in e_{-z}\uab(X),$

$$\lim_{\stackrel{\eta \to 0}{_{\eta \in \scab}}}\Vert f_{\eta}-f\Vert_{e_{-z}\uab(X)}=\lim_{\stackrel{\eta \to 0}{_{\eta \in \scab}}} \sup_{\zeta \in \scab}\Vert e^{z\zeta}f(\zeta +\eta)-e^{z\zeta}f(\zeta)\Vert_{\infty}$$ $$\le \lim_{\stackrel{\eta \to 0}{_{\eta \in \scab}}}\left (\Vert (e_zf)_{\eta}-e_{z}f\Vert_{\infty} +\vert 1-e^{-z\eta}\vert \Vert (e_{z}f)_{\eta})\Vert_{\infty}\right )=0,$$

and so, since $(e_{-\epsilon}f)_{\eta}=e^{-\epsilon \eta}e_{-\epsilon}f_{\eta},$

\begin{equation} \lim_{\stackrel{\eta \to 0, \eta \in \scab}{_{\epsilon \to 0, \epsilon \in  \sscab}}} \Vert (e_{-\epsilon}f)_{\eta}-f\Vert_{ e_{-z}\uab}$$ $$=\lim _{\stackrel{\eta \to 0, \eta \in \scab}{_{\epsilon \to 0, \epsilon \in  \sscab}}} \Vert  e_{-\epsilon}f_\eta -f\Vert_{ e_{-z}\uab}=0  \ \ (f \in  e_{-z}\uab(X), z \in \C^k)\end{equation}

Now let $f \in e_{-z}\vab(X),$ and let $\phi \in (e_{-z}\uab)'.$ If $\nu$ is a $z$-representative measure for $\phi,$ we have, for $\eta \in \scab, \epsilon \in \overline S^*_{\alpha,\beta},$

$$<(e_{-\epsilon}f)_{\eta},\phi>=e^{-\epsilon \eta}<e_{-\epsilon}f_{\eta},\phi>=e^{-(\epsilon+z) \eta}\int_{\scab} e^{-\epsilon \zeta}e^{z(\zeta+\eta)}f(\zeta +\eta)d\nu(\zeta),$$

and it follows from the Lebesgue dominated convergence theorem that we have

\begin{equation} \lim_{\stackrel{\eta \to 0, \eta \in \scab}{_{\epsilon \to 0, \epsilon \in  \sscab}}} \Vert <(e_{-\epsilon}f)_{\eta},\phi>-<f,\phi>\Vert_{X}$$ $$=\lim _{\stackrel{\eta \to 0, \eta \in \scab}{_{\epsilon \to 0, \epsilon \in  \sscab}}} \Vert <e_{-\epsilon}f_\eta  ,\phi>-<f,\phi>\Vert _{ X}=0  \ \ (f \in  e_{-z}\vab(X), \phi \in (e_{-z}\uab)', z \in \C^k)\end{equation}

The following consequence of proposition 10.6 allow to compute in some cases $<f,\phi>$ for $\phi \in \uab', f \in \vab(X), z \in Dom(\fb(\phi))$ by using the $z$-Cauchy transform. 


\begin{prop} Assume that $\alpha_j < \beta_j < \alpha_j+\pi$ for $1 \le j \le k,$ let $\phi \in \fab,$ let $z \in Dom(\fb(\phi)),$ and let $X$ be a separable Banach space.

If $f \in  e_{-z}\vab(X),$ and if

$$\int_{\tilde \partial S_{\alpha,\beta}}e^{Re(z\sigma)}\Vert f(\sigma)\Vert_X\vert d \sigma \vert <+\infty,$$

  then we have, for $\eta \in \sab,$

\begin{equation}<f_{\eta},\phi>=<f ,\phi*\delta_{\eta}>=\int_{\tilde \partial \overline S_{\alpha, \beta}}{e^{z(\sigma - \eta)}C}_z(\phi)(\sigma-\eta)f(\sigma)d\sigma.\end{equation}

 In particular we have, for  $f \in \vab(X),  \epsilon \in \ssab ,\eta \in \sab,$

\begin{equation}e^{-\epsilon \eta}<e_{-\epsilon}f_{\eta}, \phi>=<e_{-\epsilon}f ,\phi*\delta _{\eta}>=\int_{\tilde \partial \overline S_{\alpha, \beta}}e^{(z-\epsilon )(\sigma-\eta)} {\mathcal C}_z(\phi)(\sigma-\eta)f(\sigma)d\sigma.\end{equation}

\end{prop} 

Proof: Assume that $f\in e_{-z}\vab(X)$ satisfies the condition $\int_{\tilde \partial S_{\alpha,\beta}}\Vert f(\sigma)\Vert_X\vert d \sigma \vert <+\infty.$ We have, for $\eta \in \sab,$ $\epsilon \in S^*_{\alpha, \beta},$

$$<f_{\eta},\phi>=<e_zf_{\eta},\phi_{e_{-z}}>=e^{-z\eta}<(e_{z}f)_{\eta},\phi e_{-z}>, e_{-\epsilon}f_{\eta}=e^{\epsilon \eta}(e_{-\epsilon}f)_{\eta}$$

so (36) follows from (25) applied to $e_zf$ and $\phi e_{-z},$ and (37) follows from (36) applied to $e_{-\epsilon}f.$ $\square$

\smallskip

For $z \in \C^k, f \in e_{-z}H^{\infty}(\sab,X),$ define the Fourier-Borel transform of $f$ for $\zeta =(\zeta_1\dots,\zeta_k) \in \Pi_{1\le j \le k}\left (\C \setminus (-z_j -\overline S^*_{\alpha_j, \beta_j})\right )$ by the formula

$$ \FB(f)(\zeta)=\FB( e_{z}f)(z+\zeta)=\int_0^{e^{i\omega}.\infty}e^{-\zeta \sigma}f(\sigma)d\sigma$$ \begin{equation}:=\int_{0}^{e^{i\omega_1}.\infty}\dots \int_{0}^{e^{i\omega_k}.\infty}e^{-\zeta_1\sigma_1 \dots  -\zeta_k\sigma_k}f(\sigma_1,\dots,\sigma_k)d\sigma_1\dots d\sigma_k ,\end{equation}

where $\alpha _j \le \omega_j \le \beta_j$ and where $Re((z_j +\zeta_j)e^{i\omega_j})>0$ for $1\le j \le k.$

\smallskip

The following consequences of proposition 10.8, corollary 10.9 and corollary 10.11 allow to interpret the action of $\phi \in \fab$ on $ e_{-z}\uab$ for $z \in Dom(\FB(\phi))$ in terms of Fourier-Borel transforms.

\begin{prop} Let $\phi \in \fab,$ let $z =(z_1,\dots,z_k)\in \dfi,$ and let $f \in e_{-z}\vab(X).$ Set again $W_n(\zeta)=\Pi_{1\le j \le k}{n^2\over \left (n+\zeta_je^{i{\alpha_j+\beta_j\over 2}}\right )^2}$ for $\zeta =(\zeta_1, \dots, \zeta_k)\in S^*_{\alpha, \beta}, n \ge 1.$ Then

$$ (i) \ \ <f , \phi>=\lim_{\stackrel{\epsilon \to 0}{_{\epsilon \in \sab}}}\left ( \lim_{n\to +\infty}{1\over (2i\pi)^k}\int_{z+\tilde \partial \overline S^*_{\alpha, \beta}}  W_n( \sigma -z)\FB(\phi)(\sigma)\FB(f)(-\sigma+\epsilon)d\sigma\right ).$$

\smallskip

(ii)  If, further,  $\int _{\tilde \partial S^*_{\alpha, \beta}}\vert \fb(\phi)(\sigma))\vert \vert d\sigma \vert<+\infty,$ then we have, for $\epsilon \in S^*_{\alpha, \beta},$

$$<e_{-\epsilon}f, \phi>={1\over (2i\pi)^k} \int_{z+\tilde \partial \overline S^*_{\alpha, \beta}} \FB(\phi)(\sigma)\FB(f)(-\sigma +\epsilon)d\sigma,$$

and so

$$< f, \phi>=\lim_{\stackrel{\epsilon \to 0}{_{\epsilon \in \sab}}}{1\over (2i\pi)^k} \int_{z+\tilde \partial \overline S^*_{\alpha, \beta}} \FB(\phi)(\sigma)\FB(f)(-\sigma +\epsilon)d\sigma.$$ 

\smallskip

If, further, $\int_{\tilde \partial \overline \sab} e^{Re(z\sigma)}\Vert f(\sigma)\Vert \vert d\sigma \vert <+\infty,$ then

(iii) 

$$< f, \phi>={1\over (2i\pi)^k} \int_{z+\tilde \partial \overline S^*_{\alpha, \beta}} \FB(\phi)(\sigma)\FB(f)(-\sigma)d\sigma.$$ 
\end{prop}

Proof: We have $<f,\phi> = <e_zf, \phi e_{-z}>.$ Since $<\fb(f)(-\zeta-z)= \fb(e_zf)(-\zeta)$ for $\zeta \in \Pi_{1\le j \le k}\left (\C \setminus \overline S^*_{\alpha, \beta}\right ),$ and since
$\fb(\phi e_{-z})(\zeta)=<e_{-\zeta -z},\phi>=\fb(\phi)(\zeta +z)$ for $\zeta \in \overline S^*_{\alpha, \beta},$ it follows from corollary 10.11 that we have

$$<f,\phi>=<e_zf, \phi e_{-z}>=\lim_{\stackrel{\epsilon \to 0}{_{\epsilon \in S^*_{\alpha, \beta}}}}\left ( {1\over (2i\pi)^k}\int_{\tilde \partial S^{*}_{\alpha, \beta}}W_n(\zeta)\fb(\phi e_{-z})(\zeta)\fb(e_zf)(\epsilon -\zeta)d\zeta \right )$$ $$=\lim_{\stackrel{\epsilon \to 0}{_{\epsilon \in S^*_{\alpha, \beta}}}}\left ( {1\over (2i\pi)^k}\int_{\tilde \partial S^{*}_{\alpha, \beta}}W_n(\zeta)\fb(\phi)(z+\zeta)\fb(f)(-z +\epsilon -\zeta )d\zeta \right ),$$ and we obtain (i) by using the change of variables $\sigma = z +\zeta$ for $\zeta \in \tilde \partial S^*_{\alpha, \beta}.$ Using the same change of variables we deduce (ii) from proposition 10.8 and (iii) from corollary 10.9.
 $\square$

\begin{lem} Let $\alpha=(\alpha_1, \dots, \alpha_k)\in \R^k, \alpha'=(\alpha'_1, \dots, \alpha'_k)\in \R^k,$ and assume that $\alpha'_j\le \alpha_j\le \beta_j\le \beta'_j<\alpha'_j+\pi$ for $j\le k.$ Then $\cap_{z\in \C^k} e_{-z}{\mathcal U}_{\alpha',\beta'}$ is dense in  $\cap_{z\in \C^k} e_{-z}{\mathcal U}_{\alpha,\beta}.$

\end{lem}

Proof: Let $\phi \in {\mathcal F}_{\alpha,\beta},$ and assume that $<f,\phi>=0$ for every $f \in \cap_{z\in \C^k} e_{-z}{\mathcal U}_{\alpha',\beta'}.$ Let $z\in \dfi.$ Then $\FB(\phi)(z+\zeta)=0$ for every $\zeta \in \overline S^*_{\alpha',\beta'}.$ Since $\dfi$ is connected, we have $\FB(\phi)=0.$ Hence $\phi =0,$ since the Fourier-Borel transform
is one-to-one on ${\mathcal F}_{\alpha,\beta}.$ $\square$

We  can thus identify ${\mathcal F}_{\alpha,\beta}$ to a subset of ${\mathcal F}_{\alpha',\beta'}$ for $\alpha'_j\le \alpha_j\le\beta_j\le \beta'_j<\alpha'_j+\pi$ for $j\le k.$

A standard application of the Mittag-Leffler theorem of projective limits of complete metric spaces, see for example \cite{e1}, theorem 2.14, shows that we have the following result,
where as before $M_{a,b}= \{(\alpha,\beta)\in \R^k\times \R^k \ | \ a_j<\alpha_j\le \beta_j <b_j \ \mbox{if} \ a_j < b_j,  \alpha_j=\beta_j=a_j \ \mbox{if} \ a_j=b_j\}.$
\begin{prop} Let $a=(a_1,\dots,a_k) \in \R^k, b=(b_1,\dots,b_k) \in \R^k$ such that $a_j \le b_j \le a_j+\pi$ for $j \le k.$ Then $\cap_{(\alpha' ,\beta')\in M_{a,b}, \lambda \in \C^k}e_{-\lambda}{\mathcal U}_{\alpha',\beta'}$ is dense in
$ e_{-z}{\mathcal U}_{\alpha,\beta}$ for every $z \in \C$ and every $(\alpha, \beta)\in M_{a,b}.$
\end{prop}
Let $(a,b)\in \R^k\times \R^k$ be as above, and denote by $\Delta_{a, b}$ the set of all triples $(\alpha ,\beta, z)$ where $(\alpha,\beta)\in M_{a,b}$ and $z\in \C^k.$ Denote by $\preceq$ the product partial order on $\R^k$ associated to the usual order on $\R.$  If $(\alpha,\beta, z)\in \Delta_{a,b},  (\alpha', \beta', z')\in \Delta_{a, b},$  set $(\alpha, \beta, z)\preceq (\alpha' , \beta', z')$ if $\alpha' \preceq \alpha, \beta \preceq \beta'$ and $z' \in z+\overline S^*_{\alpha ', \beta'}.$
For every finite family $F= \{ (\alpha^{(l)}, \beta^{(l)},z^{(l)})\}_{1\le l \le m}$ of elements of $\Delta_{a, b},$  set

$$\sup(F)=\{ \inf_{1\le l\le m}\alpha^{(l)}\}\times\{\sup_{1\le l\le m}\beta^{(l)}\}\times \sup_{1\le l \le m}z^{(j)},$$

where $\sup_{1\le l \le m}z^{(j)}$ denotes the set of all $z\in \C^k$ satifsfying the condition

 $$z+{\overline S}^*_{\inf_{1\le l \le m}\alpha^{(l)},\sup_{1\le l\le m}\beta^{(l)}}=\cap_{1\le l\le m}\left (z^{(l)}+{\overline S}^*_{\inf_{1\le l \le k}\alpha^{(l)}, \sup_{1\le l\le k}\beta^{(l)}}\right ),$$
 
 so that $\sup_{1\le l \le m}z^{(j)},$ is the set introduced in definition 9.1(ii) when $\alpha= \inf_{1\le l \le m}\alpha^{(l)}$ and $\beta = \sup_{1\le l\le m}\beta^{(l)}.$ Notice that $\sup_{1\le l \le m}z^{(j)},$ is a singleton if $(\inf_{1\le l\le m}\alpha^{(l)})_j < (\sup_{1\le l\le m}\beta^{(l)})_j$ for $1 \le j \le k.$

It follows from the proposition that we can identify the dual of the projective limit $\cap_{(\alpha, \beta,z)\in \Delta_{a,b}}e_{-z}{\mathcal U}_{\alpha, \beta}$ to the inductive limit $\cup_{(\alpha, \beta,z)\in \Delta_{a, b}}( e_{-z}{\mathcal U}_{\alpha, \beta})'.$ This suggests the following definition.

\begin{defn} Let $a=(a_1,\dots,a_k) \in \R^k, b=(b_1,\dots,b_k) \in \R^k$ such that $a_j \le b_j \le a_j+\pi$ for $j \le k.$ Set

$${\mathcal G}_{a,b}=(\cap_{(\alpha, \beta,z)\in \Delta_{a, b}} e_{-z}{\mathcal U}_{\alpha, \beta})'=\cup_{(\alpha, \beta,z)\in \Delta_{a, b}}( e_{-z}{\mathcal U}_{\alpha, \beta})'.$$

For $\phi \in \gab,$ set $dom(\phi)=\{ (\alpha, \beta, z)\in \dab \ | \ \phi \in ( e_{-z}{\mathcal U}_{\alpha, \beta})'\}.$

\end{defn}

We thus see that the inductive limit $\gab =\cup_{(\alpha,\beta)\in M_{a,b}}{\mathcal F}_{\alpha, \beta}$ is an associative unital pseudo-Banach algebra with respect to the convolution product introduced above on the spaces ${\mathcal F}_{\alpha, \beta}.$ A subset $V$ of $\gab$ is bounded if and only if there exists $(\alpha, \beta)\in M_{a,b}$ and $z \in \C^k$ such that $V$ is a bounded subset of $( e_{-z}\uab)'.$

The proof of the following proposition is left to the reader.

\begin{prop} Let $\phi \in \gab,$ and let $(\alpha, \beta,z)\in dom(\phi).$ Then $(\alpha', \beta', z') \in dom(\phi)$ if $(\alpha, \beta, z )\preceq(\alpha', \beta', z').$ In particular if $(\phi_j)_{1\le j \le m}$ is a finite family of elements of $\gab,$ and if $(\alpha^{(j)},\beta^{(j)},z^{(j)})\in dom(\phi_j)$ for $1\le j \le m,$ then $sup_{1\le j \le m}(\alpha^{(j)},\beta^{(j)},z^{(j)}) \subset  \cap_{1\le j \le m}
dom(\phi_j)\subset dom(\phi_1*\dots*\phi_m).$

\end{prop}

\section{Appendix 3: Holomorphic functions on admissible open sets}

\begin{defn} Let  $a=(a_1,\dots,a_p) \in \R^k, b=(b_1,\dots,b_p) \in \R^k$ such that $a_j \le b_j \le a_j+\pi$ for $j \le k.$

An open set $U \subset \C^k$ is said to be admissible with respect to $(\alpha, \beta) \in M_{a,b}$ if  $U=\Pi_{1\le j \le k}U_j,$ where the open sets $U_j \subset \C$ satisfy the following conditions  for some $z =(z_1,\dots, z_k)\in \C^k,$

(i)  $U_j+\overline S^*_{\alpha_j,\beta_j}\subset U_j$ 

(ii) $U_j \subset z_j + S^*_{\alpha_j,\beta_j} ,$ and $\partial U_j -z_j=(e^{(-{\pi\over 2}-\alpha_j)i}.\infty, e^{(-\alpha_j-{\pi\over 2})i}s_{0,j})\cup \theta_j([0,1])\cup(e^{({\pi\over 2}-\beta_j)i}s_{1,j}, e^{({\pi\over 2}-\beta_j)i}.\infty),$ where $s_{0,j}\ge 0, s_{1,j}\ge 0,$ and where $\theta_j: [0,1] \to  \overline {S}^*_{\alpha_j,\beta_j} \setminus \left (e^{(-{\pi\over 2}-\alpha_j)i}.\infty, e^{(-\alpha_j-{\pi\over 2})i}s_{j,0})\cup (e^{({\pi\over 2}-\beta_j)i}s_{j,1}, e^{({\pi\over 2}-\beta_j)i}.\infty)\right )$ is a one-to-one piecewise-$\mathcal C^1$ curve such that 
$\theta_j(0)= e^{(-\alpha_j-{\pi\over 2})i}s_{j,0},$ and $\theta_j(1)=e^{({\pi\over 2}-\beta_j)i}s_{j,1}.$

If $U$ is an admissible open set with respect to some $(\alpha, \beta)\in M_{a,b}$,  $H^{(1)}(U)$ denotes the space of all functions $F$ holomorphic on $U$ such that $\Vert F \Vert_{H^{(1)}(U)}:=\sup_{\epsilon \in S^*_{\alpha, \beta}}\int_{ \tilde \partial U+\epsilon}\vert F(\sigma)\vert \vert d\sigma\vert <+\infty.$

\end{defn}

For example if $\alpha_j=\beta_j$ conditions then  conditions $(i)$ and $(ii)$ are satisfied if an only if $U_j$ is a half-plane of the form $\{ z_j\in \C \ | \ Re(z_je^{i\alpha j})> \lambda\}$
for some $\lambda \in \R.$ 

If $\alpha_j<\beta_j,$ define $\tilde x_j=\tilde x_j(\zeta_j)$ and $\tilde y_j=\tilde y_j(\zeta_j)$ for $\zeta_j \in \C$ by the formula

\begin{equation} \zeta_j= z_j +\tilde x_je^{(-{\pi\over 2}-\alpha_j)i}+\tilde y_je^{({\pi\over 2}-\beta_j)i}.\end{equation}



Notice that if $\zeta_j \in U_j,$ and if $\tilde x_j(\zeta'_j)\ge \tilde x_j(\zeta_j)$ and $\tilde y_j(\zeta'_j)\ge \tilde y_j(\zeta_j)$ then $\zeta'_j \in z_j +\overline S^*_{\alpha_j, \beta_j}\subset U_j.$ This shows that there exists
$t_{j,0}\in [0,s_{j,0}]$ and $t_{j,1}\in [0, s_{j,1}]$ and continuous piecewise ${\mathcal C}^1$-functions $f_{j}$ and $g_{j}$ defined respectively on $[0, t_{j,0}]$ and $[0,t_{j,1}]$ such that $$U_j-z_j=\{ \zeta_j \in S^*_{\alpha_j, \beta_j} \ | \ \tilde x(\zeta_j)\in (0,t_{j,0}], \tilde y(\zeta_j)>f_{j}(\tilde x(\zeta_j))\} \cup \{  \zeta_j \in S^*_{\alpha_j, \beta_j}  \ | \ \tilde x(\zeta_j)>t_{j,0}\}$$ $$
=\{ \zeta_j \in S^*_{\alpha_j, \beta_j} \ | \ \tilde y(\zeta_j)\in (0,t_{j,1}], \tilde x(\zeta_j)>g_{j}(\tilde y(\zeta_j))\} \cup \{  \zeta_j \in S^*_{\alpha_j, \beta_j}  \ | \ \tilde y(\zeta_j)>t_{j,1}\}.$$

We have $f_{j}(0)=t_{j,1}, f_{j}(t_{j,0})=0,$ $g_{j}(0)=t_{j,1}, g_{j}(t_{j,1})=0,$ $f_{j}$ and $g_{j}$ are strictly decreasing and $f_{j}=g_{j}^{-1}$ if $t_{j,s}>0$ for some, hence for all $s\in\{1,2\}.$

For $\alpha=(\alpha_1,\dots,\alpha_k)\in \R^k, \beta=(\beta_1,\dots, \beta_k)\in \R^k,$ we will use the obvious conventions $$\inf(\alpha, \beta)=(\inf(\alpha_1,\beta_1),\dots , \inf(\alpha_k, \beta_k)), \ \sup(\alpha, \beta)=(\sup(\alpha_1,\beta_1),\dots , \sup(\alpha_k, \beta_k)).$$ Clearly, $(\inf(\alpha^{(1)},\alpha^{(2)}), \sup(\beta^{(1)},\beta^{(2)})\in M_{a,b}$ if $(\alpha^{(1)},\beta^{(1)})\in M_{a,b}$ and $(\alpha^{(2)},\beta^{(2)})\in M_{a,b}.$

\begin{prop} If  $U^{(1)}$ is admissible with respect to $(\alpha^{(1)},\beta^{(1)})\in M_{a,b}$ and if $U^{(2)}$ is admissible with respect to $(\alpha^{(2)},\beta^{(2)})\in M_{a,b},$ then
$U^{(1)}\cap U^{(2)}$ is admissible with respect to $(\inf(\alpha^{(1)},\alpha^{(2)}), \sup(\beta^{(1)},\beta^{(2)})).$

\end{prop}

Set $\alpha^{(3)}= \inf(\alpha^{(1)},\alpha^{(2)}),\beta^{(3)}=\sup(\beta^{(1)},\beta^{(2)}),$ and set $U^{(3)}=U^{(1)}\cap U^{(2)}.$  The fact that $U^{(3)}$ satisfies  (i) follows from the fact that $\overline S^*_{\alpha^{(3)}_j,\beta^{(3)}_j}=\overline S^*_{\alpha^{(1)}_j,\beta^{(1)}_j}\cap \overline S^*_{\alpha^{(2)}_j,\beta^{(2)}_j}.$ The fact that $U^{(3)}$ satisfies (ii) follows easily from the fact that $\left [z^{(1)}+S^*_{\alpha^{(1)},\beta^{(1)}}\right ]\cap\left [z^{(2)}+S^*_{\alpha^{(2)},\beta^{(2)}}\right ]$ is itself admissible with respect to $( \alpha^{(3)}, \beta^{(3)})$ if $U^{(1)}$ satisfies definition 12.1 with respect to $z^{(1)}$ and if $U^{(2)}$ satisfies definition 12.1 with respect to $z^{(2)}.$  $\square$

\begin{lem} Let $U$ be an admissible open set with respect to some $(\alpha, \beta)\in M_{a,b}$,  and let $F \in H^{(1)}(U).$

(i)   We have, for  $\epsilon=(\epsilon_1, \dots, \epsilon_k) \in S^*_{\alpha, \beta},$     

$$\int_{\Pi_{j\le k}\left (U_j\setminus (\overline U_j +\epsilon_j)\right )}\vert F(\zeta)\vert dm(\zeta)\le  \vert \epsilon_1\vert \dots \vert \epsilon_k\vert\Vert F\Vert_{H^{(1)}(U)}.$$

where $m$ denotes the Lebesgue measure on $\C^k\approx \R^{2k}.$

(ii) We have, for $\zeta \in U,$

$$\vert F(\zeta)\vert  \le {2^k\over \pi^k cos \left ({\beta_1-\alpha_1\over 2}\right )\dots cos \left ({\beta_k-\alpha_k\over 2}\right )}{dist(\zeta_1,\partial S^*_{\alpha_1,\beta_1})\dots dist(\zeta_k,\partial S^*_{\alpha_k,\beta_k})\over \left [dist(\zeta_1,\partial U_1)\dots dist(\zeta_k, \partial U_k)\right ]^2}\Vert F \Vert_{H^{(1)}(U)}.$$

\end{lem}

Proof: (i) Let $F\in H^{(1)}(U),$ let $\epsilon=(\epsilon_1,\dots, \epsilon_k)\in S^*_{\alpha, \beta},$ for $j\le k$ let  $\gamma_j \in (-{\pi \over 2} -\alpha_j,{\pi \over 2}-\beta_j)$ be a determination of $arg(\epsilon_j),$ and set $r_j=\vert \epsilon _j\vert >0.$  

Set $U_{j,1}=z_j+t_{j,0}e^{(-{\pi\over 2}-\alpha_j)i} +S_{-{\pi\over 2}-\alpha_j, \gamma_j},$ $U_{j,2}=z_j+t_{j,1}e^{({\pi\over 2} -\beta_j)i} +S_{\gamma_j, {\pi\over 2}-\beta_j},$ and $U_{j,3}=z_j +\cup_{\rho >0}\left ( \rho e^{\gamma_j i}+\left ( \partial U_j\cap S^*_{\alpha_j,\beta_j}\right )\right ),$ with the convention $U_{j,3}=\emptyset$ if $t_{j,0}=t_{j,1}=0.$ Also for $\zeta_j \in \C$  set $x_i=Re(\zeta_j), y_j=Im(\zeta_j).$

 For $t_j<0,$ $0 < \rho_j < r_j,$ set $\zeta_j=\zeta_j(\rho_j,t_j)=  \rho_j e^{i\gamma_j}+(t_{j,0} -t_j)e^{-i({\pi \over 2}-\alpha_j)}.$ This gives a parametrization of $U_{j,1} \setminus (U_{j,1}+\epsilon_j),$ and we have

$$dx_jdy_j= \left | \begin{array}{lc} cos(\gamma_j)& sin(\alpha_j) \cr sin(\gamma_j) & cos(\alpha_j) \end{array}\right | d\rho_j dt_j=cos(\alpha_j+\gamma_j)d\rho_j dt_j.$$

Similarly for $t_j >t_{j,1},$ $0<\rho_j<r_j,$ set $\zeta_j=\zeta_j(\rho_j,t_j)= \rho_j e^{i\gamma_j} + t_je^{i({\pi \over 2}-\beta_j)}.$  This gives a parametrization of $U_{j,2}\setminus (U_{j,2} +\epsilon_j),$ and we have

$$dx_jdy_j= \left | \begin{array}{lc} cos(\gamma_j) & sin(\beta_j)\cr  sin(\gamma_j)& cos(\beta_j)\end{array}\right | dt_jd\rho_j=cos(\beta_j+\gamma_j)d\rho_j dt_j.$$

Now assume that  $U_{j,3} \neq \emptyset,$ so that  $t_{j,0}>0$ and $t_{j,1}>0.$ For $0<t_j<t_{j,1},$ $0<\rho_j <r_j$ set $\zeta_j=\zeta_j(\rho_j,t_j)=\rho_je^{i\gamma_j} +g_j(t_j)e^{(-{\pi\over 2}-\alpha_j)i} +t_je^{({\pi\over 2}-\beta_j)i}.$ This gives a parametrization of $U_{j,3}\setminus(U_{j,3}+\epsilon_j),$ and we have

$$dx_jdy_j=\left | \begin{array}{lc} cos(\gamma_j) & -g_j'(t)sin(\alpha_j)+sin(\beta_j)\cr sin(\gamma_j)& -g_j'(t)cos(\alpha_j)+cos(\beta_j)\end{array} \right |d\rho_jdt_j$$ $$=(cos(\beta_j+\gamma_j) -g'_j(t)cos(\alpha_j+\gamma_j))d\rho_jdt_j.$$

We have $0 < cos(\alpha_j+\gamma_j)<1,$ $0 < cos(\alpha_j+\gamma_j)<1,$ $g'_j(t_j)<0,$ and using the Cauchy-Schwartz inequality, we obtain

$$0< cos(\beta_j+\gamma_j) -g_j'(t)cos(\alpha_j+\gamma_j)$$ $$=cos(\gamma_j) (cos(\beta_j) -g'_j(t)cos(\alpha_j)) -sin(\gamma_j)(sin(\beta_j) -g'_j(t)sin(\alpha_j))$$
$$ \le \sqrt {(cos(\beta_j)-g'_j(t)cos(\alpha_j))^2 +(sin(\beta_j)-g'_j(t)sin(\alpha_j))^2}$$ $$=\sqrt{ 1-2g'_j(t)cos(\beta_j-\alpha_j)+g'_j(t)^2}.$$

On the other hand we have

$$\left \vert {\partial \zeta_j\over \partial t_j}(\rho_j,t_j)\right \vert ^2= \left (g'_j(t)e^{(-{\pi \over 2}-\alpha_j)i} +e^{({\pi\over 2}-\beta_j)i}\right ) \left (g'_j(t)e^{({\pi \over 2}+\alpha_j)i} +e^{(-{\pi\over 2}+\beta_j)i}\right )$$ $$= 1-2g'_j(t)cos(\beta_j-\alpha_j)+g'_j(t)^2.$$

The boundary $\partial U_j +\rho_je^{i\gamma_j}$ being oriented from $e^{(-{\pi\over 2}-\alpha_j)i}.\infty$ to $e^{({\pi\over 2}-\beta_j)}.\infty,$ we obtain

$$\int_{\Pi_{j\le k}\left (U_j\setminus (\overline U_j +\epsilon_j)\right )}\vert F(\zeta)\vert dm(\zeta)$$ $$\le \int _{(0,r_1)\times \dots \times (0,r_k)}\left [\int_{\Pi_{j\le k}(\partial U_j +\rho_j e^{i\gamma_j})}\vert F(\sigma_j)\vert \vert d \sigma_1\vert \dots \vert d\sigma_k\vert \right ] d\rho_1\dots d\rho_k$$ $$\le r_1\dots r_k\Vert F\Vert_{H^{(1)}(U)},$$

which proves (i).

(ii) Let $F\in H^{(1)}(U),$ let $\zeta \in U$, set $r_j=dist(\zeta_j, \partial U_j),$ set $r =(r_1,\dots,r_k),$ and set $B(\zeta,r)=\Pi_{j\le k}B(\zeta_j,r_j).$ Using Cauchy's formula and polar coordinates, we obtain the standard formula

\begin{equation}F(\zeta)={1\over \vert B(\zeta, r)\vert}\int _{B(\zeta, r)}F(\eta)dm(\eta).\end{equation}

 where $\vert B(\zeta, r)\vert=\pi^kr_1^2\dots r_k^2$ denotes the Lebesgue measure of $B(\zeta, r).$

Denote by $u_j$ the orthogonal projection of $\zeta_j$ on the real line $z_j +\R e^{i{(-{\pi\over 2}-\alpha_j)}}$, denote by $v_j$ the orthogonal projection of $\zeta_j$ on the real line $z_j +\R e^{({\pi\over 2}-\beta_j)i},$ and let $w_j\in \{u_j,v_j\}$ be such that $\vert \zeta_j-w_j\vert =\min(\vert \zeta_j-u_j\vert, \vert \zeta_j-v_j\vert).$ An easy topological argument shows that $w_j\in \partial S^*_{\alpha_j,\beta_j},$ so that $\vert \zeta_j -w_j\vert =dist (\zeta_j, \partial S^*_{\alpha_j, \beta_j}) \ge dist(\zeta_j, U_j)=r_j.$ 
For $\lambda \in \R,$ we have $\zeta_j \notin \overline S^*_{\alpha_j, \beta_j} +z_j+2(\zeta_j-w_j)+\lambda i (\zeta_j -w_j)\supset \overline U_j+2(\zeta_j-w_j)+\lambda i(\zeta_j-w_j).$ If $\pi -\beta_j+\alpha_j >{\pi \over 2},$ then $\zeta_j -w_j \in S^*_{\alpha_j, \beta_j}.$ If $\pi - \beta_j +\alpha_j \le {\pi \over 2},$ then we can choose $\lambda \in \R$ such that $\zeta_j -w_j +\lambda i (\zeta_j-w_j) \in \overline S^*_{\alpha_j, \beta_j}$ and such that $\vert \zeta_j -w_j +\lambda i (\zeta_j-w_j)\vert ={\vert \zeta_j -w_j\vert \over cos \left ({\beta_j -\alpha_j\over 2}\right )}.$ So there exists in both cases $\epsilon_j \in S^*_{\alpha_j, \beta_j}$ such that $\zeta_j \notin \overline U_j +\epsilon_j$ and $\vert \epsilon_j \vert ={2dist(\zeta_j ,\partial S^*_{\alpha_j, \beta_j})\over cos \left ({\beta_j-\alpha_j\over 2}\right )}.$

Using(40) and (i), we obtain

$$\vert F(\zeta)\vert \le {1\over \pi^kr_1^2\dots r_k^2}\int_{\Pi_{j\le k}\left (U_j \setminus(\overline U_j +\epsilon_j)\right )}\vert F(\eta)\vert dm(\eta)$$ $$\le {2^k\over \pi^k cos \left ({\beta_1-\alpha_1\over 2}\right )\dots cos \left ({\beta_k-\alpha_k\over 2}\right )}{dist(\zeta_1,\partial S^*_{\alpha_1,\beta_1})\dots dist(\zeta_k,\partial S^*_{\alpha_k,\beta_k})\over \left [dist(\zeta_1,\partial U_1)\dots dist(\zeta_k, \partial U_k)\right ]^2}\Vert F \Vert_{H^{(1)}(U)},$$

which proves (ii). $\square$

\begin{cor} $(H^{(1)}(U), \Vert . \Vert_{H^{(1)}(U)})$ is a Banach space, $F_{\overline U +\epsilon}$ is bounded on $\overline U+\epsilon,$ and $\lim_{\stackrel{dist(\zeta, \partial U )\to +\infty}{_{\zeta \in \overline U +\epsilon}}}F(\zeta)=0$ for every $F \in H^{(1)}(U)$ and every $\epsilon \in S^*_{\alpha, \beta}.$

\end{cor}

Proof: It follows from (ii) that for every compact set $K \subset U$ there exists $m_K>0$ such that $\max_{\zeta \in K}\vert F(\zeta)\vert \le m_K\Vert F \Vert_{H^{(1)}(U)}.$ So every Cauchy sequence $(F_n)_{n\ge 1}$ in $(H^{(1)}(U), \Vert . \Vert_{H^{(1)}(U)})$ is a normal family which converges uniformly on every compact subset of $U$ to a holomorphic function $F:U\to \C.$ Since $\int_{\tilde \partial U +\epsilon}\vert F(\sigma)-F_n(\sigma)\vert \vert d\sigma\vert =\lim_{R\to +\infty}\int_{B(0,R)\cap\left ( \tilde \partial U +\epsilon\right ) }\vert F(\sigma)-F_n(\sigma)\vert d\sigma\vert,$ an easy argument shows that $F\in H^{(1)}(U)$ and that $\lim_{n\to +\infty} \Vert F-F_n\Vert_{H^{(1)}(U)}=0.$ 

Now let $\epsilon >0.$ There exists $m_j>0$ such that $dist(\zeta_j, \partial U)\ge m_jdist( \zeta_j, \partial S^*_{\alpha_j, \beta_j})$ for every $\zeta_j \in \overline U_j +\epsilon_j,$ which gives, for $\zeta \in \overline U +\epsilon,$

$$\vert F(\zeta)\vert \le {2^k\over \pi^k m_1\dots m_kcos \left ({\beta_1-\alpha_1\over 2}\right )\dots cos \left ({\beta_k-\alpha_k\over 2}\right )dist (\zeta_1,\partial U_1)\dots dist(\zeta_k, \partial U_k)}\Vert F \Vert_{H^{(1)}(U)}.$$

Since $\inf_{\zeta_j \in \overline  U_j +\epsilon_j}dist(\zeta_j, \partial U_j) >0$ for $j\le k,$ this shows that $F$ is bounded on $\overline U +\epsilon,$ and that $\lim_{\stackrel{dist(\zeta, \partial U )\to +\infty}{_{\zeta \in \overline U +\epsilon}}}F(\zeta)=0.$
$\square$

\begin{thm}  Let $U$ be an admissible open set with respect to some $(\alpha, \beta)\in M_{a,b}$,  and let $F \in H^{(1)}(U).$ Then 

$$\int_{\tilde \partial U+\epsilon}F(\sigma)d\sigma=0 \ \mbox{for every} \ \epsilon \in S^*_{\alpha, \beta}, \ \mbox{and}$$

$$F(\zeta)={1\over (2i\pi)^k}\int_{\tilde \partial U+\epsilon} {F(\sigma)d\sigma\over (\zeta_1-\sigma_1)\dots (\zeta_k-\sigma_k)} \ \mbox{for every} \ \epsilon \in S^*_{\alpha, \beta} \ \mbox{and for every} \ \zeta \in U +\epsilon,$$

where $\partial U_j$ is oriented from $e^{i(-{\pi\over 2}-\alpha _j)}.\infty$ to $e^{i({\pi\over 2}-\beta j)}.\infty$ for $j\le k.$

\end{thm}

Proof:  Let $z \in \C^k$ satisfying the conditions of definition 12.1 with respect to $U,$ let $\epsilon \in S^*_{\alpha, \beta},$ let $L>1$ such that  $(z_j + e^{i\alpha_j}.\infty, z_j+Le^{i\alpha_j}] \subset \partial U_j$ and  $[z_j + Le^{i\beta_j}, z_j+e^{i\beta_j}.\infty)\subset \partial U_j$ for $ j\le k,$  and let $M>1.$ Set $$\Gamma_{L, j,1}= (\epsilon_j +\partial U_j)\setminus \left ((z_j +\epsilon_j + Le^{i\alpha_j}, z_j+\epsilon_j+ e^{i\alpha_j}.\infty) \cup (z_j +\epsilon_j + Le^{i\beta_j}, z_j+\epsilon_j +e^{i\beta_j}.\infty)\right ),$$ $$\Gamma_{L, j,2}=[z_j+\epsilon_j +Le^{i\beta_j}, z_j+L\epsilon_j+Le^{i\beta j}],$$ $$\Gamma_{L, j,3}= (L \epsilon_j +\partial U_j)\setminus \left ((z_j +L\epsilon_j + Le^{i\alpha_j}.\infty, z_j+L \epsilon_j +Le^{i\alpha_j})\\ \cup (z_j +L \epsilon _j+ Le^{i\beta}, z_j+\epsilon_j +e^{i\beta_j}.\infty)\right ),$$
$$\Gamma_{L,j,4}=[z_j +L \epsilon_j +Le^{i\alpha_j}, z_j+\epsilon_j +Le^{i\alpha_j}], \Gamma_{L,j}=\cup_{1\le s \le 4}\Gamma_{L,j,s}, $$ where the Jordan curve $\Gamma_{L,j}$ is oriented clockwise.

For $n\ge 1, \zeta_j \in \overline S^*_{\alpha_j, \beta_j},$ set $W_{j,n}(\zeta_j)={n^2 \over \left (n +e^{{\alpha_j+\beta_j\over 2}i}\zeta_j\right )^2},$ and set $W_n(\zeta)=\Pi_{j\le k}W_{n,j}(\zeta_j)$
for $\zeta \in \overline S^*_{\alpha, \beta}.$
Then $\vert W_{n,j}(\zeta_j)\vert \le 1$ for
$\zeta_j\in \overline S^*_{\alpha_j, \beta_j}$, $W_n(\zeta) \to 1$ as $n\to \infty$ uniformly on compact sets of $\overline S^*_{\alpha, \beta},$ and $\lim_{\stackrel{\vert \zeta \vert \to \infty}{_{\zeta \in S^*_{\alpha, \beta}}}}W_n(\zeta)
=0.$

Denote by $V_{L,j}$ the interior of $\Gamma_{L,j}$ and set $V_{L}=\Pi_{1\le j \le k}V_{L,j}.$ If $\zeta \in V_{L},$ it follows from Cauchy's theorem that we have
$$\sum_{l \in \{1, 2, 3, 4\}^k}\int_{\Pi_{j\le k}\Gamma_{L,j,l(j)}}W_n(\sigma-z-\epsilon)F(\sigma)d\sigma= \int_{\tilde \partial V_{L}}W_n(\sigma-z-\epsilon)F(\sigma)d\sigma=0.$$ 

Set $l_0(j)=1$ for $j\le k.$ It follows from the corollary that there exists $M>0$ such that $\vert F(\zeta)\vert \le M$ for $\zeta \in \overline U +\epsilon,$ and there exists $R_n>0$ such that $\int _{\Gamma_{L,j}}\vert W_n(\sigma_j-z_j-\epsilon_j)\vert \vert d\sigma \vert \le R_n$ for every $L.$ Also $\lim \sup_{L\to +\infty}\int_{\Gamma _{L,j,s}}\vert W_{n,j}(\sigma_j-z_j-\epsilon_j)\vert \vert d(\sigma_j)\vert=0$ for $s\ge 2, j\le k.$

Let $l\neq l_0,$ and let $j_l\le k$ such that $j_l \ge 2.$ We have

$$\lim \sup_{L\to +\infty}\left | \int_{\Pi_{j\le k}\Gamma_{L,j,l(j)}}\vert W_n(\sigma-z-\epsilon)F(\sigma)d\sigma \right \vert $$ $$\le MR_n^{k-1}\int_{\Gamma_{L,j_l,l(j_l)}}\vert W_{n,j_l}(\zeta_{j_l}-z_{j_l}-\epsilon_{j_l})\vert \vert d\sigma_{j_l}\vert=0.$$

This gives

$$\int_{\partial U +\epsilon}W_n(\sigma-z-\epsilon)F(\sigma)d\sigma=  \lim_{L\to +\infty}\int_{\Pi_{j\le k}\Gamma_{L,j,l_0(j)}}W_n(\sigma-z-\epsilon)F(\sigma)d\sigma$$ $$=\lim_{L\to +\infty}\sum_{l \in \{1, 2, 3, 4\}^k}\int_{\Pi_{j\le k}\Gamma_{L,j,l(j)}}W_n(\sigma-z-\epsilon)F(\sigma)d\sigma=0.$$

It follows then from the Lebesgue dominated convergence theorem  that $\int_{\partial U +\epsilon}F(\sigma)d\sigma=0.$

Similarly, applying Cauchy's formula when $\zeta \in U+\epsilon$ is contained in $V_L,$ we obtain

$${1\over (2i\pi)^k}\int_{\partial U +\epsilon}{W_n(\sigma-z-\epsilon)F(\sigma)\over (\zeta_1-\sigma_1)\dots (\zeta_k-\sigma_k)}d\sigma=  \lim_{L\to +\infty}{1\over (2i\pi)^k}\int_{\Gamma_{L,j,l_0(j)}}{W_n(\sigma-z-\epsilon)F(\sigma)\over (\zeta_1-\sigma_1)\dots (\zeta_k-\sigma_k)}d\sigma$$ $$=\lim_{L\to +\infty}{1\over (2i\pi)^k}\sum_{l \in \{1, 2, 3, 4\}^k}\int_{\Pi_{j\le k}\Gamma_{L,j,l(j)}}{W_n(\sigma-z-\epsilon)F(\sigma)\over (\zeta_1-\sigma_1)\dots (\zeta_k-\sigma_k)}d\sigma$$ $$=\lim_{L\to +\infty}{1\over (2i\pi)^k}\int_{\tilde \partial V_L}{W_n(\sigma-z-\epsilon)F(\sigma)\over (\zeta_1-\sigma_1)\dots (\zeta_k-\sigma_k)}d\sigma=W_n(\zeta-z-\epsilon)F(\zeta).$$

It follows then again from the Lebesgue dominated convergence theorem that we have

$$F(\zeta)=\lim_{n\to +\infty}W_n(\zeta-z-\epsilon)F(\zeta)=\lim_{n\to +\infty}{1\over (2i\pi)^k}\int_{\partial U +\epsilon}{W_n(\sigma-z-\epsilon)F(\sigma)\over (\zeta_1-\sigma_1)\dots (\zeta_k-\sigma_k)}d\sigma$$ $$={1\over (2i\pi)^k}\int_{\partial U +\epsilon}{F(\sigma)\over (\zeta_1-\sigma_1)\dots (\zeta_k-\sigma_k)}d\sigma.$$ $\square$

Let $\zeta \in U,$ and let $\epsilon \in S^*_{\alpha, \beta}.$ It follows from the theorem that there exists $\rho >0$ such that $\vert F(\zeta)\vert \le {1\over (2\pi)^k}{\Vert F\Vert_{H^{(1)}(U)}\over \Pi_{j\le k}dist(\zeta_j,\partial U_j+t\epsilon_j)}$ for $t\in (0,\rho].$ 
Since $\lim_{t\to 0^+}dist(\zeta_j, \partial U_j +t\epsilon_j )=\lim_{t\to 0^+}dist(\zeta_j-t\epsilon_j, \partial U_j)=dist(\zeta_j, \partial U_j),$ we obtain, for $F\in H^{(1)}(U),$ $\zeta \in U,$

\begin{equation} \vert F(\zeta)\vert \le {1\over (2\pi)^k}{\Vert F\Vert_{H^{(1)}(U)}\over \Pi_{j\le k}dist(\zeta_j,\partial U_j)}  \end{equation}

which improves inequality (ii) of lemma 12.3.

If $\alpha_j=\beta_j$ for $j\le k$, then every $(\alpha, \beta)$ admissible open set $U$ is a product of open half-planes and the space $H^{(1)}(U)$ is the usual Hardy space $H^{1}(U).$ The standard conformal mappings of the open unit disc $D$ onto half planes induce an isometry from the Hardy space $H^1(D^k)$ onto $H^1(U).$ It follows then from standard results about $H^1(D^k)$, see theorems 3.3.3 and 3.3.4 of \cite{ru}, that $F$ admits a.e. a nontangential limit $F^*$ on $\tilde \partial U,$ and that $\lim_{\epsilon \to 0}\vert \int_{\tilde \partial U}\vert F^*(\sigma) -F(\sigma +\epsilon)\vert \vert d\sigma\vert =0.$ This gives the formula

\begin{equation} F(\zeta)={1\over (2i\pi)^k}\int_{\tilde \partial U} {F^*(\sigma)d\sigma\over (\zeta_1-\sigma_1)\dots (\zeta_k-\sigma_k)} \ \mbox{for every}  \ \zeta \in U.\end{equation}

We did not investigate whether such nontangential limits of $F$ on $\tilde \partial U$ exist in the general case.

Recall that the Smirnov class $\mathcal N^+(P^+)$ on the right-hand open half-plane $P^+$ consists in those functions $F$ holomorphic on $P^+$ which can be written under the form $F=G/H$ where $G \in H^{\infty}(P^+)$ and where $H\in H^{\infty}(P^+)$ is outer, which means that we have, for $Re(\zeta)>0,$

$$F(\zeta)=exp\left ({1\over \pi}\int_{-\infty}^{+\infty}{1-iy\zeta\over( \zeta-iy)(1+y^2)}log\vert F^*(iy)\vert dy\right ),$$

where $F^*(iy)=\lim_{x\to 0^+}F(x+iy)$ if defined a.e. on the vertical axis and satisfies $\int_{-\infty}^{+\infty}{\vert log\vert F^*(iy)\vert\over1+y^2}dy<+\infty.$

Set , for $Re(\zeta)>0,$

$$F_n(\zeta)=exp\left ({1\over \pi}\int_{-\infty}^{+\infty}{1-iy\zeta\over(\zeta- iy)(1+y^2)}\sup(log\vert F^*(iy)\vert ,-n)dy\right ).$$

It follows from the positivity of the Poisson kernel on the real line that $\vert F(\zeta)\vert \le \vert F_n(\zeta)\vert$ and that $\lim_{n\to +\infty}F_n\zeta)=F(\zeta)$ for $Re(\zeta)>0.$
Also the nontangential limit $F_n^*(iy)$ of $F$ at $iy$ exists a.e. on the imaginary axis and $\vert F_n^*(iy)\vert=\sup(e^{-n},\vert F^*(iy)\vert)$ a.e., which shows that $\sup_{\zeta \in P^+}\vert F_n(\zeta)\vert =\sup_{\zeta \in P^+}\vert F(\zeta)\vert$ when $n$ is sufficiently large. Hence $\lim_{n\to +\infty}F(\zeta)F_n^{-1}(\zeta)=1$ for $\zeta \in P^+.$

This suggests the following notion;

\begin{defn} Let $U\subset \C^k$ be a connected open set. A holomorphic function $F\in H^{\infty}(U)$ is said to be strongly outer on $U$ if there exists a sequence $(F_n)_{n\ge 1}$ of invertible elements of $H^{\infty}(U)$ satisfying the following conditions

(i) $\vert F(\zeta)\vert\le \vert F_n(\zeta)\vert \ \ (\zeta \in U, n\ge 1),$

(ii) $\lim_{n\to +\infty}F(\zeta)F_n^{-1}(\zeta)=1 \ \ (\zeta \in U).$

The Smirnov class $\mathcal \mathcal S(U)$ consists of those holomorphic functions  $F$ on $U$ such that $FG\in H^{\infty}(U)$ for some strongly outer function $G \in H^{\infty}(U).$
\end{defn}

It follows from (ii) that $F(\zeta)\neq 0$ for every $\zeta \in U$ if $F$ is strongly outer on $U,$ and $F_{|_V}$ is strongly outer on $V$ if $V\subset U.$ Similarly if $F \in \mathcal S(U)$ then $F_{|_V}\in  \mathcal S(V).$ Also it follows immediately from the definition that the set of bounded strongly outer functions on $U$ is stable under products, and that if there is a conformal mapping $\theta$ from an open set $V \subset \C^k$ onto $U$ then $F \in H^{\infty}(U)$ is strongly outer on $U$ if and only if $F \circ \theta$ is strongly outer on $V$, and if $G$ is holomorphic on $U$ then $G\in {\mathcal S}(U)$ if and only $F\circ \theta \in \mathcal S(V).$ 

Now let $(\alpha,\beta)\in M_{a,b}$ and let $U=\Pi_{j\le k}U_j$ be an admissible open set with respect to $(\alpha, \beta).$ Then each set $U_j$ is conformally equivalent to the open unit disc $\D,$ and so there exists a conformal mapping $\theta$ from $\D^k$ onto $U$, and  the study of the class of bounded strongly outer functions on $U$ (resp. the Smirnov class on $U$) reduces to the study of the class bounded strongly outer functions (resp. the Smirnov class) on $\D^k.$

Let $F\in H^{\infty}(D^k)$ be strongly outer, and let $(F_n)_{n\ge 1}$ be a sequence of invertible elements of $H^{\infty}(\D^k)$ satisfying the conditions of definition 12.6 with respect to $F.$ Denote by $\T=\partial D$ the unit circle. Then $H^{\infty}(\D^k)$ can be identified to a $w^*$-closed subspace of $L^{\infty}(\T^k)$ with respect to the $w^*$-topology $\sigma(L^{1}(\T^k), L^{\infty}(\T^k)).$ Let $L\in H^{\infty}(\D^k)$ be a $w^*$-cluster point of the sequence $(FF_n^{-1})_{n\ge 1}.$ Since the map $G\to G(\zeta)$ is $w^*$-continuous on $H^{\infty}(\D^k)$ for $\zeta \in \D^k,$ $L=1,$ and so $FH^{\infty}(D^k)$ is $w^*$-dense in $H^{\infty}.$ When $k=1,$ this implies as well-known that $F$ is outer, and the argument used for the half-plane shows that, conversely, every bounded outer function on $\D$ is strongly outer, and so ${\mathcal S}(\D)=\mathcal N^+(\D).$

Recall that a function $G \in H^{\infty}(\D^k)$ is said  to be outer if $log(\vert G(0,\dots,0)\vert) ={1\over (2\pi)^k}\int_{\T^k}log\vert G^*(e^{it_1},\dots,e^{it_k})\vert dt_1\dots dt_k,$ where $G^*(e^{it_1},\dots, e^{it_k})$ denotes a.e. the nontangential limit of $G$ at $(e^{it_1},\dots, e^{it_k}),$ see \cite{ru}, definition 4.4.3, and $G$ is outer if and only if almost every slice function $G_{\omega}$ is outer on $\D,$ where $G_{\omega}(\zeta)=G(\omega \zeta)$ for $\omega \in \T^k, \zeta \in \D,$ see \cite{ru}, lemma 4.4.4. If follows from definition 12.6 that every slice function $F_\omega$ is strongly outer on $\D$ if $F\in H^{\infty}(\D^k)$ is strongly outer on $\D^k,$ and so every strongly outer bounded function on $\D^k$ is outer. It follows from an example from \cite{ru} that the converse is false if $k\ge 2$.

\begin{prop} Let $k\ge 2,$ and set $F(\zeta_1, \dots, \zeta_k)=e^{\zeta_1+\zeta_2+2\over \zeta_1+\zeta_2-2}$ for $(\zeta_1,\dots, \zeta_k)\in \D^k.$ Then $F$ is outer on $\D^k,$ but $F$ is not strongly outer on $\D^k.$

\end{prop}

Proof: Set $f(\zeta)=e^{\zeta +1\over \zeta-1}$ for $\zeta \in \D.$ Then $f\in H^{\infty}(\D)$ is a singular inner function. Since $f(\zeta)\neq 0$ for $\zeta \in \D,$ it follows from \cite{ru}, lemma 4.4.4b that the function $\tilde f:(\zeta_1,\zeta_2)\to f({\zeta_1+\zeta_2\over 2},{\zeta_1+\zeta_2\over 2})=e^{\zeta_1+\zeta_2+2\over \zeta_1+\zeta_2-2}$ is outer on $\D^2.$ Hence we have

$$log\vert F(0,\dots,0)\vert =log\vert \tilde f (0,0)\vert={1\over (2\pi)^2}\int_{\T^2}\tilde f(e^{it_1},e^{it_2})dt_1dt_2$$ $$={1\over (2\pi)^k}\int_{\T^k} F(e^{it_1},\dots,e^{it_k})dt_1\dots dt_k,$$

and so $F$ is outer on $\D^k.$

Now set $\omega=(1, \dots, 1).$ Then $F_{\omega}=f$ is not outer on $\D,$ and so $F$ is not strongly outer on $\D^k.$ $\square$

The fact that some bounded outer functions on $\D$ are not strongly outer is not surprising: The Poisson intergal of a real valued integrable function on $\T^k$ is the real part of some
holomorphic function on $\D^k$ if an only if its Fourier coefficients vanish on $\Z^k\setminus (\Z^+)^k\cup (\Z^-)^k,$ see \cite{ru}, theorem 2.4.1, and so the construction of the sequence $(F_n)_{n\ge 1}$ satisfying the conditions of definition 12.6 with respect to a bounded outer function $F$ on $\D^k$ breaks down when $k\ge 2.$
We conclude this appendix with the following trivial observations.

\begin{prop} Let $U=\Pi_{j\le k}U_j \subset \C^k$ be an admissible open set with respect to some $(\alpha,\beta) \in M_{a,b}.$ 

(i) Let $\theta_j:U_j \to \D$ be a conformal map and let $\pi_j: (\zeta_1,\dots, \zeta_k) \to \zeta_j$ be the j-th cooordinate projection. If $f\in H^{\infty}(\D)$ is outer, then $f\circ \theta_j\circ \pi_j $ is strongly outer on $U.$

(ii) The Smirnov class $\mathcal S(U)$ contains all holomorphic functions on $U$ having polynomial growth at infinity.

\end{prop}

Proof: (i) Since $f$ is strongly outer on $D,$ there exists a sequence $(f_n)_{n\ge 1}$ of invertible elements of $H^{\infty}(D)$ satisfying the conditions of definition 12.6 with respect to $f.$ Then the sequence $(f_n\circ \theta_j \circ \pi_j)_{n\ge 1}$ satisfies the conditions  of definition 12.6 with respect to $f\circ \theta_j\circ \pi_j,$ and so $f\circ \theta_j\circ \pi_j$ is strongly outer on $U.$

(ii) For $j \le k$ there exists $\gamma_j \in [-\pi,\pi)$ and $m_j\in \R$ such that  open set $U_j$ is contained in the open half plane $P_j:=\{\zeta_j\in \C \ | \ Re(\zeta_je^{i\gamma_j})\ge m_j\}.$ The function $\sigma \to {1-\sigma\over 2}$ is outer on $\D,$ since $\vert {1-\sigma\over 2}\vert \le \left \vert {1+1/n-\sigma \over 2}\right \vert$ for $\sigma \in \D,$ and the function $\zeta_j \to {\zeta_je^{i\gamma_j}-m_j-1\over \zeta_je^{i\gamma_j}-m_j+1}$ maps conformally $U_j$ onto $\D.$ Set $F_j(\zeta_1, \dots, \zeta_k)= {1-{\zeta_je^{i\gamma_j}-m_j-1\over \zeta_je^{i\gamma_j}-m_j+1}\over 2}={1\over  \zeta_je^{i\gamma_j}-m_j+1}.$ It follows from (i) that $F_j$ is strongly outer on $\Pi_{j\le k}P_j,$ hence strongly outer on $U.$

Now assume that a function $F$ holomorphic on  $U$ has polynomial growth at infinity. Then there exists $p\ge 1$ such that $F\Pi_{j\le k}F_j^p$ is bounded on $U,$ and so $F\in \mathcal S(U).$ $\square$

IMB, UMR 5251

Universit\'e de Bordeaux

351, cours de la Lib\'eration

33405 - Talence

esterle@math.u-bordeaux.fr
  

\begin{thebibliography}{55}
\bibitem{ri} F. Ricci, Hardy spaces in one variable, 
 
 {\small{\color{blue}homepage.sns.it/fricci/papers/hardy.pdf}}
\end{thebibliography}
  \end{document}